\documentclass[preprint, 10pt]{elsarticle}

\usepackage[margin=2.5cm]{geometry}
\UseRawInputEncoding
\usepackage{setspace}
\usepackage{lmodern}
\usepackage{amsmath,amssymb}
\usepackage{amsfonts,amsthm}
\usepackage{amscd,amsxtra,latexsym}
\usepackage{graphicx}
\usepackage{subfigure}
\usepackage{epsf,epsfig}
\usepackage{bm,bbm}
\usepackage{algorithm}
\usepackage{physics}
\usepackage{color}
\usepackage[footnotesize,bf]{caption}
\usepackage{mathtools}
\usepackage{hhline}
\usepackage[T1]{fontenc}
\usepackage{placeins}
\newtheorem{rmk}{Remark}

\usepackage{multirow,dcolumn}
\usepackage{mathrsfs}
\usepackage{verbatim}
\usepackage{comment}
\usepackage{epstopdf}
\usepackage{lineno}

\graphicspath{{figures/}}
\journal{JCP}

\begin{document}
\begin{frontmatter}
    \title{Statistical higher-order multi-scale method for nonlinear thermo-mechanical simulation of random composite materials with temperature-dependent properties}
    \author[XDU]{Hao Dong\corref{cor}}
    \ead{donghaoxd@xidian.edu.cn}
    \cortext[cor]{Corresponding Author}
    \author[LESC]{Junzhi Cui}
    \address[XDU]{School of Mathematics and Statistics, Xidian University, Xi'an 710071, China}
    \address[LESC]{Academy of Mathematics and Systems Science, Chinese Academy of Sciences, Beijing 100190, China.}

\begin{abstract}
Stochastic multi-scale modeling and simulation for nonlinear thermo-mechanical problems of composite materials with complicated random microstructures remains a challenging issue. In this paper, we develop a novel statistical higher-order multi-scale (SHOMS) method for nonlinear thermo-mechanical simulation of random composite materials, which is designed to overcome limitations of prohibitive computation involving the macro-scale and micro-scale. By virtue of statistical multi-scale asymptotic analysis and Taylor series method, the SHOMS computational model is rigorously derived for accurately analyzing nonlinear thermo-mechanical responses of random composite materials both in the macro-scale and micro-scale. Moreover, the local error analysis of SHOMS solutions in the point-wise sense clearly illustrates the crucial indispensability of establishing the higher-order asymptotic corrected terms in SHOMS computational model for keeping the conservation of local energy and momentum. Then, the corresponding space-time multi-scale numerical algorithm with off-line and on-line stages is designed to efficiently simulate nonlinear thermo-mechanical behaviors of random composite materials. Finally, extensive numerical experiments are presented to gauge the efficiency and accuracy of the proposed SHOMS approach.
\end{abstract}


\end{frontmatter}
\section{Introduction}
In recent years, random composite materials have been extensively applied in a variety of engineering sectors, such as aviation, aerospace and civil construction, etc. By randomly distributing high-performance fibrous or particulate materials into ordinary matrix material, these synthetic composite materials exhibit high temperature resistance, high fatigue resistance and high fracture resistance, etc \cite{R1,R2}. Especially in aviation and aerospace industries, engineering structures manufactured by random composite materials often served under extreme heat environment while the thermal and mechanical properties of component materials exhibit significantly nonlinear temperature-dependent feature. These complicated nonlinear physical behaviors and randomly geometric heterogeneities of the considered structures raise a grand challenge for effective numerical simulation \cite{R3}.

To the best of our knowledge, traditional numerical methods including the finite element method (FEM) \cite{R33,R34}, boundary element method \cite{R35} and meshless method \cite{R36} have been adopted to the analysis and computation of nonlinear thermomechanical problems. Moreover, Abdoun et al. used homotopy and asymptotic numerical method to simulate and
analyze the thermal buckling and vibration of laminated composite plates with temperature-dependent properties
in \cite{R37}. In reference \cite{R38}, Najibi et al. employed higher-order graded finite element method to conduct transient thermal stress analysis for a hollow FGM cylinder with nonlinear temperature-dependent material properties. In reference \cite{R39}, the state space method and transfer-matrix method are adopted to obtain the displacements and stresses for the thick beams with temperature-dependent material properties under thermo-mechanical loads. However, it should be noted the equations, which govern the nonlinear thermo-mechanical behaviors for the composites, have rapidly varying and strongly discontinuous coefficients arising from the sharp variation between different constituents. As far as we know, the direct numerical simulation for composite materials needs a tremendous amount of computational resources or even ineffective to capture their microscopic behaviors due to the highly heterogeneous components.

To accomplish effective modeling and efficient simulation for inhomogeneous materials, scientists and engineers presented a variety of multi-scale methods, such as asymptotic homogenization method (AHM) \cite{R4}, multi-scale finite element method (MsFEM) \cite{R5}, heterogeneous multi-scale method (HMM) \cite{R6}, variational multi-scale method (VMS) \cite{R7}, multi-scale eigenelement method (MEM) \cite{R8}, localized orthogonal decomposition method (LOD) \cite{R9} and finite volume based asymptotic homogenization theory (FVBAHT) \cite{R10}, etc. However, numerical computation and theoretical analysis in \cite{R11,R12,R13} find that most of above-mentioned multi-scale methods are lower-order multi-scale method in essence, which can only capture macroscopic and inadequate microscopic information of heterogeneous materials, especially for high-contrast composite materials. To improve inadequate numerical accuracy of classical lower-order multi-scale approaches, Cui and his research team systematically developed a class of higher-order multi-scale methods, whose numerical accuracy is significantly improved for simulating authentic composite materials in practical engineering applications. Hence, these higher-order multi-scale approaches are extensively used in multi-physics coupling problems, stochastic multi-scale problem, structural mechanics problem and nonlinear multi-scale problem of heterogeneous materials, etc \cite{R14,R15,R16,R17,R18,R19,R20,R21}.
The reviews of above-mentioned multi-scale approaches show that these methods have a strong potential to encourage important advances in modeling and simulating a certain range of composites' behaviors. However, they still need to be improved for
composite materials with complex non-deterministic microstructure. The uncertainties in the microstructure prominently affect the mechanical properties of the composite materials. Some stochastic multi-scale computational schemes have been established in recent years based on perturbation-based stochastic finite element method \cite{R22,R23,R24}, spectral stochastic finite element method \cite{R25,R26,R27} and stochastic collocation method \cite{R28,R29} for specific problems. Furthermore, combining with Monte Carlo method, the higher-order multi-scale methods proposed by Cui and his research team have been applied to simulate a wide range of physical behaviors of random composite materials \cite{R15,R16,R21,R30,R31,R32}. However, it is worth noting that there are few works about multi-scale thermo-mechanical simulation of random composite materials with temperature-dependent properties. Hence, it is of great theoretical and engineering values to develop effective multi-scale approaches for nonlinear thermo-mechanical simulation of random composite materials with temperature-dependent properties.

The reminder of our article is organized as follows: In Section 2 the investigated random nonlinear governing equations with initial-boundary conditions are herein introduced for describing heat conduction and mechanical deformation of random composites with temperature-dependent properties. Furthermore, stochastic multi-scale asymptotic analysis and Taylor expansion approach are employed to establish statistical higher-order multi-scale computational model for nonlinear thermo-mechanical simulation of random composites. Also, the local error analysis of SHOMS computational model are derived.
In Section 3 we develop a space-time multi-scale numerical algorithm with off-line and on-line stages to efficiently simulate nonlinear thermo-mechanical behaviors of random composite materials. Illustrative examples are presented to validate the computational accuracy and efficiency of the proposed SHOMS computational model and corresponding numerical algorithm in Section 4. In the final section, we provide some meaningful conclusions and a brief outlook.
\section{The establishment of statistical higher-order multi-scale computational model}
\subsection{Microscopic computer representation of random composite materials}
In this study, the investigated composite materials are comprised of matrix and randomly distributed particles or fibers, as shown in Fig.\hspace{1mm}1.
\begin{figure}[!htb]
\centering
\begin{minipage}[c]{0.46\textwidth}
  \centering
  \includegraphics[width=0.9\linewidth,totalheight=1.6in]{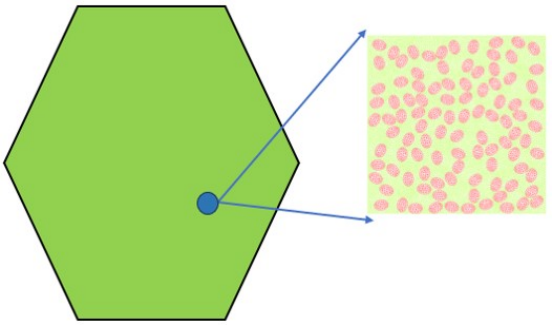} \\
  (a)
\end{minipage}
\begin{minipage}[c]{0.46\textwidth}
  \centering
  \includegraphics[width=0.9\linewidth,totalheight=1.6in]{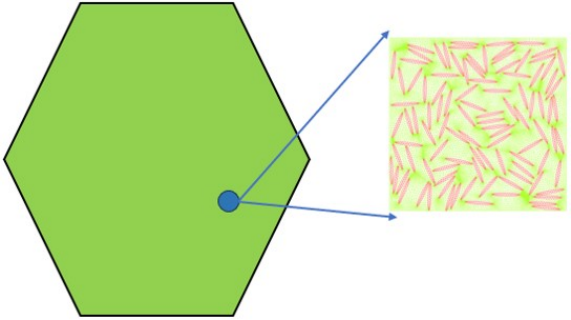} \\
  (b)
\end{minipage}
\caption{Random composite structures with different microscopic configurations. (a) Random particulate material; (b) Random fibrous material.}
\end{figure}

Based on computer representation idea and its improved algorithm devised by Li, Cui and Dong \cite{R40,R41}, we employ open-source Freefem++ software to establish the detailed computer representation algorithm for generating microscopic configurations of random composite materials as follows.
\begin{enumerate}[(S1)]
\item
Regarding random particulate and fibrous materials in Fig.\hspace{1mm}1(a) and Fig.\hspace{1mm}1(b), the probability distribution model is first employed to generate the random geometric parameters $\left(x_{1}, x_{2}, a, b, \theta_1\right)$ or $\left(x_{1}, x_{2}, x_{3}, a, b, c, \theta_1, \theta_2, \theta_3\right)$ of 2D or 3D randomly distributed configurations.
\item
Then, judge whether the newly generated configuration is located inside RVE and whether the newly generated configuration intersects with other previously generated configurations. For randomly distributed configurations, we use whether the distance between the centers of the previously generated configurations and newly generated configuration is greater than the sum of the radii of previously generated configurations and newly generated configuration as discriminate criterion. Additionally, to enhance the packing ratio of microscopic inclusions, we use whether there exist intersection points on previously generated configurations when connecting the centers of previously generated configurations and the points on the surface of newly generated configuration as revised discriminate criterion.
\item
When generating a sufficient amount of microscopic configurations, mesh generation algorithm based on Delaunay Refinement method (Freefem++ command: "buildmesh" or "tetg" for 2D or 3D geometrical configurations respectively) is adopted to create the microscopic configurations of the investigated random composites \cite{R44}.
\end{enumerate}
For geometric parameters $\left(x_{1}, x_{2}, a, b, \theta_1\right)$ in 2D case, $x_{1}$ and $x_{2}$ represent the central coordinates of the elliptical inclusion for the $x$-axis and $y$-axis. $a$ and $b$ denote the lengths of the long half-axis and short half-axis of the elliptical inclusion and $\theta_1$ represents the intersection angle between the long half-axis of the elliptical inclusion and the $x$-axis. For geometric parameters $\left(x_{1}, x_{2}, x_{3}, a, b, c, \theta_1, \theta_2, \theta_3\right)$ in 3D case, $x_{1}$, $x_{2}$ and $x_{3}$ are the central coordinates of the ellipsoidal inclusion for the $x$-axis, $y$-axis and $z$-axis. $a$, $b$ and $c$ denote the lengths of the long half-axis, middle half-axis and short half-axis of the ellipsoidal inclusion. $\theta_1$, $\theta_2$ and $\theta_3$ are the three Euler angles of the ellipsoidal inclusion, respectively. Moreover, by elongating the elliptical or ellipsoidal particles and increasing the ratio of their long half-axis to short half-axis, the elliptical or ellipsoidal particles can be changed as fibrous inclusions. To sum up, the above-mentioned methodologies accomplish the effective generation of finite element mesh for random composites we investigated at micro-scale.
\subsection{Setting of stochastic multi-scale nonlinear thermo-mechanical problems}
The primary challenge for solving random multi-scale problems pertains to their auxiliary cell problems defined on the entire space $\mathbb{R}^\mathcal{N}(\mathcal{N}=2,3)$. To tackle this challenge, by using "periodization" and "cutoff" techniques in previous studies \cite{R15,R16,R21}, the unit cell problems defined on an infinite domain are approximated by transforming them into unit cell problems on a finite domain with infinite random sampling, see Fig.\hspace{1mm}2 for a schematic explanation.
\begin{figure}[!htb]
\centering
  \includegraphics[width=0.7\linewidth,totalheight=2.6in]{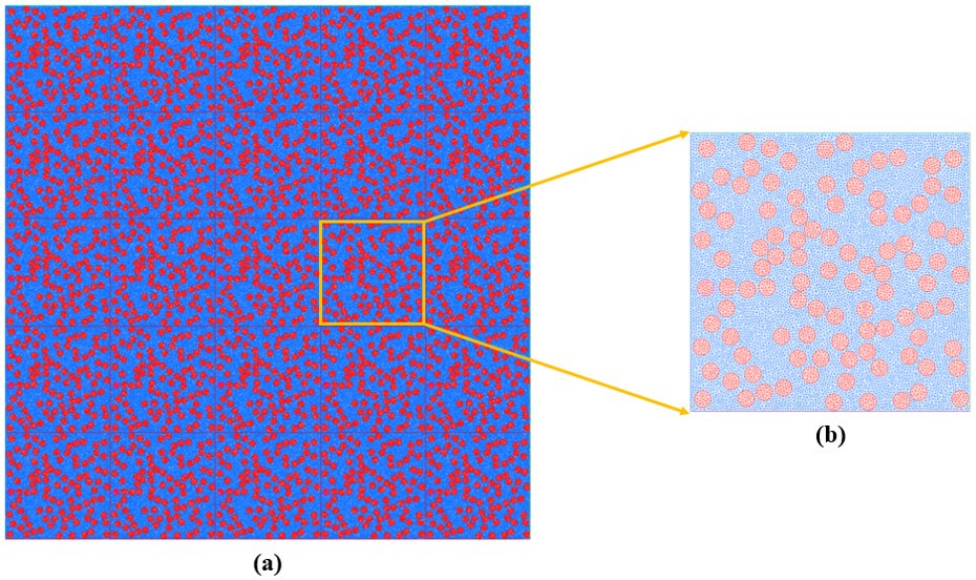}
\caption{Random composite structures with a statistical periodic layout. (a) Macroscopic composite structure $\Omega$; (b) Microscopic unit cell $Y^s$.}
\end{figure}

Based on the classical thermo-mechanical model \cite{R3}, the stochastic governing equations for describing the nonlinear thermo-mechanical problems of random composite materials are reported, whose material parameters all possess the temperature-dependent properties.
\begin{equation}
\left\{
\begin{aligned}
&\rho^{\varepsilon} ({\bm{x}},{T^\varepsilon },\bm{\omega})c^{\varepsilon}({\bm{x}},{T^\varepsilon },\bm{\omega})\frac{{\partial {T^\varepsilon (\bm{x},t,\bm{\omega}) }}}{{\partial t}}- \frac{\partial }{{\partial {x_i}}}\Big( {{k_{ij}^{\varepsilon}}({\bm{x}},{T^\varepsilon },\bm{\omega})\frac{{\partial {T^\varepsilon (\bm{x},t,\bm{\omega})}}}{{\partial {x_j}}}}\Big)= h(\bm{x},t),\;\;\text{in}\;\;\Omega\times(0,T^{*}),\\
&- \frac{\partial }{{\partial {x_j}}}\Big( {{C_{ijkl}^{\varepsilon}}({\bm{x}},{T^\varepsilon },\bm{\omega})\frac{{\partial u_k^\varepsilon (\bm{x},t,\bm{\omega}) }}{{\partial {x_l}}}}- {{\beta _{ij}^{\varepsilon}}({\bm{x}},{T^\varepsilon},\bm{\omega})( {{T^\varepsilon(\bm{x},t,\bm{\omega})} - \widetilde T} )}  \Big) = {f_i}(\bm{x},t),\;\;\text{in}\;\;\Omega\times(0,T^{*}),\\
&T^{{\varepsilon }}(\bm{x},t,\bm{\omega}) = \widehat T(\bm{x},t),\;\;\text{on}\;\;\partial {\Omega_T}\times(0,T^{*}),\\
&\bm{u}^{{\varepsilon}}(\bm{x},t,\bm{\omega})=\widehat{\bm{u}}(\bm{x},t),\;\;\text{on}\;\;\partial\Omega_{u}\times(0,T^{*}),\\
&{k_{ij}^{{\varepsilon}}({\bm{x}},{T^\varepsilon },\bm{\omega})\frac{\partial T^{{\varepsilon }}(\bm{x},t,\bm{\omega})}{\partial {x_j}}}{n_i} = \bar q(\bm{x},t),\;\;\text{on}\;\;\partial {\Omega_q}\times(0,T^{*}),\\
&\big[ {{C_{ijkl}^{{\varepsilon }}}({\bm{x}},{T^\varepsilon },\bm{\omega})\frac{{\partial u_k^{{\varepsilon }}(\bm{x},t,\bm{\omega})}}{{\partial {x_l}}}}- {\beta _{ij}^{{\varepsilon }}({\bm{x}},{T^\varepsilon },\bm{\omega})( {{T^{{\varepsilon }}(\bm{x},t,\bm{\omega})} - \widetilde T} )}  \big]n_j=\bar{\sigma}_i(\bm{x},t),\;\;\text{on}\;\;\partial\Omega_{\sigma}\times(0,T^{*}),\\
&T^{{\varepsilon }}({\bm{x}},0,\bm{\omega})=\widetilde T,\;\;\text{in}\;\;\Omega.
\end{aligned} \right.
\end{equation}
Here, $\Omega$ represents a bounded convex domain in $\mathbb{R}^\mathcal{N}(\mathcal{N}=2,3)$ with a boundary $\partial\Omega$. In the micro-scale, the domain $\Omega$ can be defined by a statistical periodic layout of microscopic unit cell $Y^s$ corresponding to random sample $\omega^s$ ($s=1,2,3,\dots$ denotes the index of random samples), as shown in Fig.\hspace{1mm}2. The characteristic size of microscopic cell $Y^s$ is characterized by a parameter $\varepsilon$. The $u_i^{\varepsilon}({\bm{x}},t,\bm{\omega})$ and $T^{\varepsilon}({\bm{x}},t,\bm{\omega})$ in equations (2.1) are the unresolved displacement and temperature fields. $\displaystyle\rho^{{\varepsilon }}({\bm{x}},{T^\varepsilon },\bm{\omega})$ is the mass density; $\displaystyle c^{{\varepsilon }}({\bm{x}},{T^\varepsilon },\bm{\omega})$ is specific heat; $\displaystyle k_{ij}^{{\varepsilon }}({\bm{x}},{T^\varepsilon },\bm{\omega})$ is the second order thermal conductivity tensor; $\displaystyle C_{ijkl}^{{\varepsilon }}({\bm{x}},{T^\varepsilon },\bm{\omega})$ is the fourth order elastic tensor; $\displaystyle\beta_{ij}^{{\varepsilon }}({\bm{x}},{T^\varepsilon },\bm{\omega})$ is the second order thermal modulus tensor. Furthermore, we assume that all material parameters satisfy Lipschitz continuous condition with respect to temperature variable $T^{\varepsilon}$ and are statistical periodic functions \cite{R14}. $\widehat{\bm{u}}({\bm{x}},t)$ is the prescribed displacements on the boundary $\partial\Omega_{u}$, $\widehat{T}({\bm{x}},t)$ is the prescribed temperature on the boundary $\partial\Omega_{T}$, $\bar{\sigma}_i({\bm{x}},t)$ is the prescribed traction on the boundary $\partial\Omega_{\sigma}$ with the normal vector $n_j$, and $\bar{q}({\bm{x}},t)$ is the prescribed heat flux normal to the boundary $\partial\Omega_{q}$ with the normal vector $n_i$; $\widetilde{T}$ is the initial temperature when the composites are in stress-free state; The body forces and internal heat source are represented by $f_i({\bm{x}},t)$ and $h({\bm{x}},t)$, respectively.

To begin with, let us set $\displaystyle\bm{y}={\bm{x}}/{\varepsilon}$ as microscopic coordinates of statistical periodic unit cell $Y^s=(0,1)^\mathcal{N}$. With this notation, the material parameters $\rho^{\varepsilon}({\bm{x}},{T^\varepsilon },\bm{\omega})$, $c^{\varepsilon}({\bm{x}},{T^\varepsilon },\bm{\omega})$, $k_{ij}^{\varepsilon}({\bm{x}},{T^\varepsilon },\bm{\omega})$, $C_{ijkl}^{\varepsilon}({\bm{x}},{T^\varepsilon },\bm{\omega})$ and $\beta_{ij}^{\varepsilon}({\bm{x}},{T^\varepsilon },\bm{\omega})$ can be rewritten as $\rho(\bm{y},{T^\varepsilon },\bm{\omega})$, $c(\bm{y},{T^\varepsilon },\bm{\omega})$, $k_{ij}(\bm{y},{T^\varepsilon },\bm{\omega})$, $C_{ijkl}(\bm{y},{T^\varepsilon },\bm{\omega})$ and $\beta_{ij}(\bm{y},{T^\varepsilon },\bm{\omega})$. Moreover, we have the chain rule for the performed spatial scales as follows.
\begin{equation}
\displaystyle\frac{\partial \Phi^\varepsilon (\bm{x},t,\bm{\omega})}{\partial x_i}=\frac{\partial \Phi (\bm{x},\bm{y},t,\bm{\omega})}{\partial x_i}+\frac{1}{\varepsilon}\frac{\partial \Phi (\bm{x},\bm{y},t,\bm{\omega})}{\partial y_i}.
\end{equation}
which will be extensively used in the sequel.

Moreover, in line with previous works such as \cite{R16,R17,R18,R19,R20}, we make certain assumptions for the governing equations (2.1).
\begin{enumerate}[(A)]
\item[(A)]
$k_{ij}^{\varepsilon}({\bm{x}},{T^\varepsilon },\bm{\omega})$, $C_{ijkl}^{\varepsilon}({\bm{x}},{T^\varepsilon },\bm{\omega})$ and $\beta_{ij}^{\varepsilon}({\bm{x}},{T^\varepsilon },\bm{\omega})$ are symmetric, and there exist two positive constant $\gamma_0$ and $\gamma_1$ independent of $\varepsilon$ such that
\begin{displaymath}
\begin{aligned}
&k_{ij}^{\varepsilon}({\bm{x}},{T^\varepsilon },\bm{\omega})=k_{ji}^{\varepsilon}({\bm{x}},{T^\varepsilon },\bm{\omega}),\gamma_0|\bm{\xi}|^2\leq k_{ij}^{\varepsilon}({\bm{x}},{T^\varepsilon },\bm{\omega})\xi_i\xi_j  \le\gamma_1|\bm{\xi}|^2,\\
&C_{ijkl}^{\varepsilon}({\bm{x}},{T^\varepsilon },\bm{\omega})=C_{ijlk}^{\varepsilon}({\bm{x}},{T^\varepsilon },\bm{\omega})=C_{klij}^{\varepsilon}({\bm{x}},{T^\varepsilon },\bm{\omega}),\gamma_0\eta_{ij}\eta_{ij}\leq C_{ijkl}^{\varepsilon}({\bm{x}},{T^\varepsilon },\bm{\omega})\eta_{ij}\eta_{kl} \le\gamma_1\eta_{ij}\eta_{ij},\\
&\beta_{ij}^{\varepsilon}({\bm{x}},{T^\varepsilon },\bm{\omega})=\beta_{ji}^{\varepsilon}({\bm{x}},{T^\varepsilon },\bm{\omega}),\gamma_0|\bm{\xi}|^2\leq \beta_{ij}^{\varepsilon}({\bm{x}},{T^\varepsilon },\bm{\omega})\xi_i\xi_j  \le\gamma_1|\bm{\xi}|^2.
\end{aligned}
\end{displaymath}
where $\{\eta_{ij}\}$ is an arbitrary symmetric matrix in $R^{\mathcal{N}\times\mathcal{N}}$, $\bm{\xi}=(\xi_1,\xi_2\cdots\xi_\mathcal{N})$ is an arbitrary vector with
real elements in $R^\mathcal{N}$, and $({\bm{x}},{T^\varepsilon })$ is an arbitrary point in $\Omega\times[T_{min},T_{max}+C_*]$.
\item[(B)]
$\rho^{\varepsilon}({\bm{x}},{T^\varepsilon },\bm{\omega})$, $c^{\varepsilon}({\bm{x}},{T^\varepsilon },\bm{\omega})$, $k_{ij}^{\varepsilon}({\bm{x}},{T^\varepsilon },\bm{\omega})$, $\displaystyle C_{ijkl}^{\varepsilon}({\bm{x}},{T^\varepsilon },\bm{\omega})$ and $\displaystyle\beta_{ij}^{\varepsilon}({\bm{x}},{T^\varepsilon },\bm{\omega})$ are scalar functions belonging to $L^\infty (\Omega)$, and there are constants $\rho^0$ and $c^0$ such that
\begin{displaymath}
0<\rho^0\leq \rho^{\varepsilon}({\bm{x}},{T^\varepsilon },\bm{\omega}),\;\;0<c^0\leq c^{\varepsilon}({\bm{x}},{T^\varepsilon },\bm{\omega}).
\end{displaymath}
functions $\rho^{\varepsilon}({\bm{x}},{T^\varepsilon },\bm{\omega})$, $c^{\varepsilon}({\bm{x}},{T^\varepsilon },\bm{\omega})$, $k_{ij}^{\varepsilon}({\bm{x}},{T^\varepsilon },\bm{\omega})$, $\displaystyle C_{ijkl}^{\varepsilon}({\bm{x}},{T^\varepsilon },\bm{\omega})$ and $\displaystyle\beta_{ij}^{\varepsilon}({\bm{x}},{T^\varepsilon },\bm{\omega})$ are $1-$periodic of microscopic variable ${\bm{y}}$ for stationary ${T^\varepsilon }\in [T_{min},T_{max}+C_*]$.
\item[(C)]
$h(\bm{x},t)\in L^2(\Omega\times(0,T^{*}))$, $f_{i}(\bm{x},t)\in L^2(\Omega\times(0,T^{*}))$, $\widehat{T}(\bm{x},t)\in L^2(0,T^{*};H^{1}(\Omega))$, $\widehat{\bm{u}}(\bm{x},t)\in L^2(0,T^{*};H^{1}(\Omega))^{3}$, $\bar{q}(\bm{x},t)\in L^2(\Omega\times(0,T^{*}))$, $\bar{\sigma}_i(\bm{x},t)\in L^2(\Omega\times(0,T^{*}))$.
\end{enumerate}
\subsection{Statistical higher-order multi-scale computational model of the governing equations}
The proposed second-order two-scale computational model is established in this subsection. To the multi-scale problem (2.1), we suppose that ${T^\varepsilon }(\bm{x},t,\bm{\omega})$ and $u_i^\varepsilon(\bm{x},t,\bm{\omega})$ can be expanded as the following power series representations.
\begin{equation}
\left\{\begin{array}{l}
{T^\varepsilon }(\bm{x},t,\bm{\omega}) ={T_{0}}(\bm{x},\bm{y},t,\bm{\omega}) + \varepsilon {T_{1}}(\bm{x},\bm{y},t,\bm{\omega}) + {\varepsilon ^2}{T_{2}}(\bm{x},\bm{y},t,\bm{\omega}) + {\rm O}({\varepsilon ^3}),\\
u_i^\varepsilon(\bm{x},t,\bm{\omega})= u_{i0}(\bm{x},\bm{y},t,\bm{\omega}) + \varepsilon u_{i1}(\bm{x},\bm{y},t,\bm{\omega}) + {\varepsilon ^2}u_{i2}(\bm{x},\bm{y},t,\bm{\omega}) + {\rm O}({\varepsilon ^3}).
\end{array}\right.
\end{equation}

Below we focus on the implementation of solving the specific expanding forms of $u _{i}^\varepsilon(\bm{x},t,\bm{\omega})$ and $T^\varepsilon(\bm{x},t,\bm{\omega})$. Firstly, the Taylor's formula of multivariate function is given as below \cite{R20}
\begin{equation}
f({x_0},{y_0} + h,z_0) = f({x_0},{y_0},z_0) + {f_y}({x_0},{y_0},z_0)h + \frac{1}{2}{f_{yy}}({x_0},{y_0},z_0){h^2} + O({h^3}).
\end{equation}
Then, the Taylor's formula (2.4) can be further rewritten by multi-index notation as follows \cite{R20}
\begin{equation}
f({x_0},{y_0} + h,z_0) = f({x_0},{y_0},z_0) + \mathbf{D}^{(0,1,0)}{f}({x_0},{y_0},z_0)h + \frac{1}{2}\mathbf{D}^{(0,2,0)}{f}({x_0},{y_0},z_0){h^2} + O({h^3}).
\end{equation}
By aid of the above-mentioned Taylor's formula (2.4) and multi-index notation (2.5), the material parameters depending on temperature $T^\varepsilon$ thus can be expanded as \cite{R20}
\begin{equation}
\begin{aligned}
&\quad\; k_{ij}^\varepsilon ({\bm{x}},{T^\varepsilon },\bm{\omega})= {k_{ij}}({\bm{y}},{T^{\varepsilon}},\bm{\omega}) = {C_{ijkl}}({\bm{y}},{T_{0}} + \varepsilon {T_{1}} + {\varepsilon ^2}{T_{2}} + {\rm O}({\varepsilon ^3}),\bm{\omega})\\
& = {C_{ijkl}}({\bm{y}},{T_{0}},\bm{\omega}) + \mathbf{D}^{(0,1,0)}{C_{ijkl}}({\bm{y}},{T_{0}},\bm{\omega})\big[ {\varepsilon {T_{1}} + {\varepsilon ^2}{T_{2}} + {\rm O}({\varepsilon ^3})} \big]\\
&+ \frac{1}{2}\mathbf{D}^{(0,2,0)}{C_{ijkl}}({\bm{y}},{T_{0}},\bm{\omega}){\big[ {\varepsilon {T_{1}} + {\varepsilon ^2}{T_{2}} + {\rm O}({\varepsilon ^3})} \big]^2} + O\Big( \big[ {\varepsilon {T_{1}} + {\varepsilon ^2}{T_{2}} + {\rm O}({\varepsilon ^3})} \big]^3 \Big)\\
& = {C_{ijkl}}({\bm{y}},{T_{0}},\bm{\omega}) + \varepsilon {T_{1}}\mathbf{D}^{(0,1,0)}{C_{ijkl}}({\bm{y}},{T_{0}},\bm{\omega})\\
&+ {\varepsilon ^2}\big[ {{T_{2}}\mathbf{D}^{(0,1,0)}{C_{ijkl}}({\bm{y}},{T_{0}},\bm{\omega}) + \frac{1}{2}{{( {{T_{1}}} )}^2}\mathbf{D}^{(0,2,0)}{C_{ijkl}}({\bm{y}},{T_{0}},\bm{\omega})} \big] + {\rm O}({\varepsilon ^3})\\
&= C_{ijkl}^{(0)} + \varepsilon C_{ijkl}^{(1)} + {\varepsilon ^2}C_{ijkl}^{(2)} + {\rm O}({\varepsilon ^3}).
\end{aligned}
\end{equation}
Using the similar expansion method as (2.6), other material parameters ${\rho^\varepsilon }({\bm{x}},{T^\varepsilon })$, $c^\varepsilon({\bm{x}},{T^\varepsilon })$, $C_{ijkl}^\varepsilon ({\bm{x}},{T^\varepsilon })$ and $\beta _{ij}^\varepsilon ({\bm{x}},{T^\varepsilon })$ are also expanded as the following forms
\begin{equation}
\begin{aligned}
{\rho^\varepsilon }({\bm{x}},{T^\varepsilon },\bm{\omega})& = \rho({\bm{y}},{T_{0}},\bm{\omega}) + \varepsilon {T_{1}}\mathbf{D}^{(0,1,0)}\rho({\bm{y}},{T_{0}},\bm{\omega})\\
& + {\varepsilon ^2}\big[ {{T_{2}}\mathbf{D}^{(0,1,0)}\rho({\bm{y}},{T_{0}},\bm{\omega}) + \frac{1}{2}{( {{T_{1}}} )}^2\mathbf{D}^{(0,2,0)}\rho({\bm{y}},{T_{0}},\bm{\omega})} \big] + {\rm O}({\varepsilon ^3})\\
& = {\rho^{(0)}} + \varepsilon {\rho^{(1)}} + {\varepsilon ^2}{\rho^{(2)}} + {\rm O}({\varepsilon ^3}),\\
{c^\varepsilon }({\bm{x}},{T^\varepsilon },\bm{\omega})& = c({\bm{y}},{T_{0}},\bm{\omega}) + \varepsilon {T_{1}}\mathbf{D}^{(0,1,0)}c({\bm{y}},{T_{0}},\bm{\omega})\\
& + {\varepsilon ^2}\big[ {{T_{2}}\mathbf{D}^{(0,1,0)}c({\bm{y}},{T_{0}},\bm{\omega}) + \frac{1}{2}{( {{T_{1}}} )}^2\mathbf{D}^{(0,2,0)}c({\bm{y}},{T_{0}},\bm{\omega})} \big] + {\rm O}({\varepsilon ^3})\\
& = {c^{(0)}} + \varepsilon {c^{(1)}} + {\varepsilon ^2}{c^{(2)}} + {\rm O}({\varepsilon ^3}),\\
C_{ijkl}^\varepsilon ({\bm{x}},{T^\varepsilon },\bm{\omega})& = {C_{ijkl}}({\bm{y}},{T_{0}},\bm{\omega}) + \varepsilon {T_{1}}\mathbf{D}^{(0,1,0)}C_{ijkl}({\bm{y}},{T_{0}},\bm{\omega})\\
& + {\varepsilon ^2}\big[ {{T_{2}}\mathbf{D}^{(0,1,0)}{C_{ijkl}}({\bm{y}},{T_{0}},\bm{\omega})+ \frac{1}{2}{{( {{T_{1}}} )}^2}\mathbf{D}^{(0,2,0)}{C_{ijkl}}({\bm{y}},{T_{0}},\bm{\omega}_{\bm{x}^{'}})} \big] + {\rm O}({\varepsilon ^3})\\
& = C_{ijkl}^{(0)} + \varepsilon C_{ijkl}^{(1)} + {\varepsilon ^2}C_{ijkl}^{(2)} + {\rm O}({\varepsilon ^3}),\\
\beta _{ij}^\varepsilon ({\bm{x}},{T^\varepsilon },\bm{\omega})&= {\beta _{ij}}({\bm{y}},{T_{0}},\bm{\omega}) + \varepsilon {T_{1}}\mathbf{D}^{(0,1,0)}\beta_{ij}({\bm{y}},{T_{0}},\bm{\omega})\\
& + {\varepsilon ^2}\big[ {{T_{2}}\mathbf{D}^{(0,1,0)}{\beta _{ij}}({\bm{y}},{T_{0}},\bm{\omega}) + \frac{1}{2}{{( {{T_{1}}} )}^2}\mathbf{D}^{(0,2,0)}{\beta _{ij}}({\bm{y}},{T_{0}},\bm{\omega})} \big] + {\rm O}({\varepsilon ^3})\\
&= \beta _{ij}^{(0)} + \varepsilon \beta _{ij}^{(1)} + {\varepsilon ^2}\beta _{ij}^{(2)} + {\rm O}({\varepsilon ^3}).
\end{aligned}
\end{equation}

By substituting equations (2.3), (2.6), and (2.7) into the multi-scale problem (2.1) and utilizing the chain rule (2.2), we can expand the derivatives and match terms with the same order of the small parameter $\varepsilon$. This process allows us to reasonably obtain the following equations
\begin{equation}
\left\{\begin{aligned}
&{\varepsilon ^{ - 2}}\frac{\partial }{{\partial {y_i}}}\Big( {{k_{ij}^{(0)}}\frac{{\partial T_{0} }}{{\partial {y_j}}}} \Big)+{\varepsilon ^{ - 1}}\frac{\partial }{{\partial {y_i}}}\Big( {{k_{ij}^{(1)}}\frac{{\partial T_{0} }}{{\partial {y_j}}}} \Big)+{\varepsilon ^{ - 1}}\frac{\partial }{{\partial {y_i}}}\Big( {{k_{ij}^{(0)}}(\frac{{\partial T_{0} }}{{\partial {x_j}}}+\frac{{\partial T_{1} }}{{\partial {y_j}}})} \Big)\\
&+{\varepsilon ^{ -1}}\frac{\partial }{{\partial {x_i}}}\Big( {{k_{ij}^{(0)}}\frac{{\partial T_{0} }}{{\partial {y_j}}}} \Big)+{\varepsilon ^{0}}\frac{\partial }{{\partial {y_i}}}\Big( {{k_{ij}^{(2)}}\frac{{\partial T_{0} }}{{\partial {y_j}}}} \Big)+{\varepsilon ^{0}}\frac{\partial }{{\partial {x_i}}}\Big( {{k_{ij}^{(1)}}\frac{{\partial T_{0} }}{{\partial {y_j}}}} \Big)+{\varepsilon ^{0}}\frac{\partial }{{\partial {y_i}}}\Big( {{k_{ij}^{(1)}}(\frac{{\partial T_{0} }}{{\partial {x_j}}}+\frac{{\partial T_{1} }}{{\partial {y_j}}})} \Big)\\
&+{\varepsilon ^{0}}\frac{\partial }{{\partial {y_i}}}\Big( {{k_{ij}^{(0)}}(\frac{{\partial T_{1} }}{{\partial {x_j}}}+\frac{{\partial T_{2} }}{{\partial {y_j}}})} \Big)+{\varepsilon ^{0}}\frac{\partial }{{\partial {x_i}}}\Big( {{k_{ij}^{(0)}}(\frac{{\partial T_{0} }}{{\partial {x_j}}}+\frac{{\partial T_{1} }}{{\partial {y_j}}})} \Big)+{\rm O}(\varepsilon )={\varepsilon ^{0}}\rho^{(0)}c^{(0)}\frac{{\partial {T_{0}}}}{{\partial t}}-h,\\
&{\varepsilon ^{ - 2}}\frac{\partial }{{\partial {y_j}}}\Big( {{C_{ijkl}^{(0)}}\frac{{\partial u_{k0} }}{{\partial {y_l}}}} \Big)+{\varepsilon ^{ - 1}}\frac{\partial }{{\partial {y_j}}}\Big( {{C_{ijkl}^{(1)}}\frac{{\partial u_{k0} }}{{\partial {y_l}}}} \Big)+{\varepsilon ^{ - 1}}\frac{\partial }{{\partial {y_j}}}\Big( {C_{ijkl}^{(0)}(\frac{\partial u_{k0} }{\partial {x_l}}}+\frac{\partial u_{k1} }{\partial {y_l}}) \Big)\\
&+{\varepsilon ^{ - 1}}\frac{\partial }{{\partial {x_j}}}\Big( {C_{ijkl}^{(0)}\frac{\partial u_{k0} }{\partial {y_l}}} \Big)-{\varepsilon ^{ - 1}} \frac{\partial }{{\partial {y_j}}}\Big({\beta _{ij}^{(0)}}(T_{0}-\widetilde T)\Big)+{\varepsilon ^{0}}\frac{\partial }{{\partial {y_j}}}\Big( {{C_{ijkl}^{(2)}}\frac{{\partial u_{k0} }}{{\partial {y_l}}}} \Big)\\
&+{\varepsilon ^{0}}\frac{\partial }{{\partial {y_j}}}\Big( {C_{ijkl}^{(1)}(\frac{\partial u_{k0} }{\partial {x_l}}}+\frac{\partial u_{k1} }{\partial {y_l}}) \Big)+{\varepsilon ^{0}}\frac{\partial }{{\partial {x_j}}}\Big( {C_{ijkl}^{(1)}\frac{\partial u_{k0} }{\partial {y_l}}} \Big)-{\varepsilon ^{0}} \frac{\partial }{{\partial {y_j}}}\Big({\beta _{ij}^{(1)}}(T_{0}-\widetilde T)\Big)\\
&+{\varepsilon ^{0}}\frac{\partial }{{\partial {y_j}}}\Big( {C_{ijkl}^{(0)}(\frac{\partial u_{k1} }{\partial {x_l}}}+\frac{\partial u_{k2} }{\partial {y_l}}) \Big)-{\varepsilon ^{0}}\frac{\partial }{{\partial {x_j}}}\Big({\beta _{ij}^{(0)}}(T_{0}-\widetilde T)\Big)\\
&-{\varepsilon ^{0}}\frac{\partial }{{\partial {y_j}}}\Big({\beta _{ij}^{(0)}}T_{1}\Big)+{\varepsilon ^{0}}\frac{\partial }{{\partial {x_j}}}\Big( {C_{ijkl}^{(0)}(\frac{\partial u_{k0} }{\partial {x_l}}}+\frac{\partial u_{k1} }{\partial {y_l}}) \Big)+ {\rm O}(\varepsilon )=-{f_i}.
\end{aligned}\right.
\end{equation}
As a consequence of equations (2.8), a series of equations are derived by matching terms of the same order of $\varepsilon$, following the classical procedure of AHM
\begin{equation}
{\rm O}({\varepsilon ^{ - 2}}):\left\{ \begin{aligned}
&\frac{\partial }{{\partial {y_i}}}\Big( {{k_{ij}^{(0)}}({\bm{y}},{T_{0}},\bm{\omega})\frac{\partial T_{0}}{\partial {y_j}}}\Big)=0,\\
&{\frac{\partial }{{\partial {y_j}}}\Big( {{C_{ijkl}^{(0)}}({\bm{y}},{T_{0}},\bm{\omega})\frac{{\partial u_{k0} }}{{\partial {y_l}}}} \Big)}=0.
\end{aligned} \right.
\end{equation}
\begin{equation}
{\rm O}({\varepsilon ^{ - 1}}):\left\{ \begin{aligned}
&\frac{\partial }{{\partial {y_i}}}\Big( {{k_{ij}^{(1)}}({\bm{y}},{T_{0}},\bm{\omega})\frac{{\partial T_{0} }}{{\partial {y_j}}}}\Big)+\frac{\partial }{{\partial {y_i}}}\Big( {{k_{ij}^{(0)}}({\bm{y}},{T_{0}},\bm{\omega})(\frac{{\partial T_{0} }}{{\partial {x_j}}}+\frac{{\partial T_{1}}}{{\partial {y_j}}})}\Big)\\
&+\frac{\partial }{{\partial {x_i}}}\Big({{k_{ij}^{(0)}}({\bm{y}},{T_{0}},\bm{\omega})\frac{{\partial T_{0} }}{{\partial {y_j}}}} \Big)=0,\\
&\frac{\partial }{{\partial {y_j}}}\Big( {{C_{ijkl}^{(1)}}({\bm{y}},{T_{0}},\bm{\omega})\frac{{\partial u_{k0} }}{{\partial {y_l}}}}\Big)+\frac{\partial }{{\partial {y_j}}}\Big( {C_{ijkl}^{(0)}({\bm{y}},{T_{0}},\bm{\omega})(\frac{\partial u_{k0} }{\partial {x_l}}}+\frac{\partial u_{k1} }{\partial {y_l}}) \Big)\\
&+\frac{\partial }{{\partial {x_j}}}\Big( {C_{ijkl}^{(0)}({\bm{y}},{T_{0}},\bm{\omega})\frac{\partial u_{k0} }{\partial {y_l}}} \Big)-\frac{\partial }{{\partial {y_j}}}\Big({\beta _{ij}^{(0)}}({\bm{y}},{T_{0}},\bm{\omega})(T_{0}-\widetilde T)\Big)=0.
\end{aligned} \right.
\end{equation}
\begin{equation}
{\rm O}({\varepsilon ^0}):\left\{ \begin{aligned}
&\frac{\partial }{{\partial {y_i}}}\Big( {{k_{ij}^{(2)}}({\bm{y}},{T_{0}},\bm{\omega})\frac{{\partial T_{0} }}{{\partial {y_j}}}} \Big)+\frac{\partial }{{\partial {x_i}}}\Big( {{k_{ij}^{(1)}}({\bm{y}},{T_{0}},\bm{\omega})\frac{{\partial T_{0} }}{{\partial {y_j}}}} \Big)\\
&+\frac{\partial }{{\partial {y_i}}}\Big( {{k_{ij}^{(1)}}({\bm{y}},{T_{0}},\bm{\omega})(\frac{{\partial T_{0} }}{{\partial {x_j}}}+\frac{{\partial T_{1} }}{{\partial {y_j}}})} \Big)+\frac{\partial }{{\partial {y_i}}}\Big( {{k_{ij}^{(0)}}({\bm{y}},{T_{0}},\bm{\omega})(\frac{{\partial T_{1} }}{{\partial {x_j}}}+\frac{{\partial T_{2}}}{{\partial {y_j}}})} \Big)\\
&+\frac{\partial }{{\partial {x_i}}}\Big( {{k_{ij}^{(0)}}({\bm{y}},{T_{0}},\bm{\omega})(\frac{{\partial T_{0} }}{{\partial {x_j}}}+\frac{{\partial T_{1} }}{{\partial {y_j}}})} \Big)=\rho^{(0)}({\bm{y}},{T_{0}},\bm{\omega})c^{(0)}({\bm{y}},{T_{0}},\bm{\omega})\frac{{\partial {T_{0}}}}{{\partial t}}-h
,\\
&\frac{\partial }{{\partial {y_j}}}\Big( {{C_{ijkl}^{(2)}}({\bm{y}},{T_{0}},\bm{\omega})\frac{{\partial u_{k0} }}{{\partial {y_l}}}} \Big)+\frac{\partial }{{\partial {y_j}}}\Big( {C_{ijkl}^{(1)}({\bm{y}},{T_{0}},\bm{\omega})(\frac{\partial u_{k0} }{\partial {x_l}}}+\frac{\partial u_{k1} }{\partial {y_l}}) \Big)\\
&+\frac{\partial }{{\partial {x_j}}}\Big( {C_{ijkl}^{(1)}({\bm{y}},{T_{0}},\bm{\omega})\frac{\partial u_{k0} }{\partial {y_l}}} \Big)- \frac{\partial }{{\partial {y_j}}}\Big({\beta _{ij}^{(1)}}({\bm{y}},{T_{0}},\bm{\omega})(T_{0}-\widetilde T)\Big)\\
&+\frac{\partial }{{\partial {y_j}}}\Big( {C_{ijkl}^{(0)}({\bm{y}},{T_{0}},\bm{\omega})(\frac{\partial u_{k1} }{\partial {x_l}}}+\frac{\partial u_{k2} }{\partial {y_l}}) \Big)-\frac{\partial }{{\partial {x_j}}}\Big({\beta _{ij}^{(0)}}({\bm{y}},{T_{0}},\bm{\omega})(T_{0}-\widetilde T)\Big)\\
&-\frac{\partial }{{\partial {y_j}}}\Big({\beta _{ij}^{(0)}}({\bm{y}},{T_{0}},\bm{\omega})T_{1}\Big)+\frac{\partial }{{\partial {x_j}}}\Big( {C_{ijkl}^{(0)}({\bm{y}},{T_{0}},\bm{\omega})(\frac{\partial u_{k0} }{\partial {x_l}}}+\frac{\partial u_{k1} }{\partial {y_l}}) \Big)=-{f_i}.
\end{aligned} \right.
\end{equation}

From (2.9), we next can determine that
\begin{equation}
\begin{aligned}
T_{0}(\bm{x},\bm{y},t,\bm{\omega}) =T_{0}(\bm{x},t),\;\;u_{i0}(\bm{x},\bm{y},t,\bm{\omega})= u_{i0}(\bm{x},t).
\end{aligned}
\end{equation}
Following this, taking advantage of (2.12), the terms $\displaystyle\frac{\partial u_{k0}}{\partial y_l}$ and $\displaystyle\frac{\partial T_{0}}{\partial y_j}$ both equate to zero. Subsequently, equations (2.10) can be further simplified as the subsequent equations.
\begin{equation}
\left\{ \begin{aligned}
&\frac{\partial }{{\partial {y_i}}}\Big( {{k_{ij}^{(0)}}({\bm{y}},{T_{0}},\bm{\omega})\frac{{\partial T_{1} }}{{\partial {y_j}}}} \Big)=-\frac{\partial }{{\partial {y_i}}}\Big( {{k_{ij}^{(0)}}({\bm{y}},{T_{0}},\bm{\omega})\frac{{\partial T_{0} }}{{\partial {x_j}}}} \Big),\\
&\frac{\partial }{{\partial {y_j}}}\Big( {C_{ijkl}^{(0)}({\bm{y}},{T_{0}},\bm{\omega})}\frac{\partial u_{k1} }{\partial {y_l}} \Big)=-\frac{\partial }{{\partial {y_j}}}\Big( {C_{ijkl}^{(0)}({\bm{y}},{T_{0}},\bm{\omega})\frac{\partial u_{k0} }{\partial {x_l}}} \Big)+\frac{\partial }{{\partial {y_j}}}\Big({\beta _{ij}^{(0)}}({\bm{y}},{T_{0}},\bm{\omega})(T_{0}-\widetilde T)\Big).
\end{aligned} \right.
\end{equation}
According to equations (2.13), we construct the separation forms of the first-order correctors $u_{i1}$ and $T_{1}$
\begin{equation}
\left\{
\begin{aligned}
&{T_{1}}(\bm{x},\bm{y},t,\bm{\omega}) = {M_{\alpha_1}}(\bm{y},T_{0},\bm{\omega})\frac{\partial T_{0}}{\partial x_{\alpha_1}},\\
&u_{i1}(\bm{x},\bm{y},t,\bm{\omega}) = N_{im}^{\alpha_1}(\bm{y},T_{0},\bm{\omega})\frac{\partial u_{m0}}{\partial x_{\alpha_1}} - {P_i}(\bm{y},T_{0},\bm{\omega})({T_{0}} - \widetilde T)
,\;\;\alpha_1,\alpha_2,m,n=1,2,3.
\end{aligned}\right.
\end{equation}
where ${M_{\alpha_1}}$, $N_{im}^{\alpha_1}$ and $P_i$ are $1-$periodic functions defined in microscopic unit cell $Y^s$ for any random sample $\omega^s$, which are defined as the first-order auxiliary cell functions. Now, substituting (2.14) into (2.13), the following equations with homogeneous Dirichlet boundary condition are obtained after simplification and calculation
\begin{equation}
\left\{
\begin{aligned}
&\frac{\partial}{\partial y_i}\big[ { k_{ij}^{(0)}(\bm{y},T_{0},\bm{\omega}^s){\frac{\partial M_{\alpha_1}}{\partial y_j}}} \big]= -\frac{\partial k_{i{\alpha_1}}^{(0)}(\bm{y},T_{0},\bm{\omega}^s)}{\partial y_i},&\;\;\;&\bm{y}\in Y^s, \\
&M_{\alpha_1}(\bm{y},T_{0},\bm{\omega}^s)= 0,&\;\;\;&\bm{y}\in\partial Y^s.
\end{aligned} \right.
\end{equation}
\begin{equation}
\left\{
\begin{aligned}
&\frac{\partial}{\partial y_j}\big[ { C_{ijkl}^{(0)}({\bm{y}},{T_{0}},\bm{\omega}^s){\frac{\partial N_{km}^{\alpha_1}}{\partial y_l}}} \big]= -\frac{\partial C_{ijm{\alpha_1}}^{(0)}({\bm{y}},{T_{0}},\bm{\omega}^s)}{\partial y_j},&\;\;\;&\bm{y}\in Y^s, \\
&N_{km}^{\alpha_1}(\bm{y},T_{0},\bm{\omega}^s)= 0,&\;\;\;&\bm{y}\in\partial Y^s.
\end{aligned} \right.
\end{equation}
\begin{equation}
\left\{
\begin{aligned}
&\frac{\partial}{\partial y_j}\big[ { C_{ijkl}^{(0)}(\bm{y},T_{0},\bm{\omega}^s){\frac{\partial P_k}{\partial y_l}}} \big]= -\frac{\partial \beta _{ij}^{(0)}(\bm{y},T_{0},\bm{\omega}^s)}{\partial y_j},&\;\;\;&\bm{y}\in Y^s, \\
&P_k(\bm{y},T_{0},\bm{\omega}^s)= 0,&\;\;\;&\bm{y}\in\partial Y^s.
\end{aligned} \right.
\end{equation}
\begin{rmk}
It should be mentioned that the first-order auxiliary cell functions are quasi-periodic functions which all depend on the expansion term $T_{0}$. Each expansion term $T_{0}$ acts as a varying parameter. This is a significant difference compared to linear periodic composites.
\end{rmk}
Subsequently, we perform a volume integral on both sides of equations (2.11) on microscopic unit cell $Y^s$ and apply the Gauss theorem on equations (2.11). By further applying the Kolmogorov's strong law of large numbers, these procedures allow to derive the macroscopic homogenized equations associated with multi-scale problem (2.1) as presented below
\begin{equation}
\left\{ \begin{aligned}
&\widehat S(T_{0})\frac{{\partial {T_{0}}}}{{\partial t}}- \frac{\partial }{{\partial {x_i}}}\Big( {{\widehat k_{ij}}(T_{0})\frac{{\partial {T_{0}}}}{{\partial {x_j}}}}\Big) = h,\;\;\text{in}\;\;\Omega\times(0,T^{*}),\\
&- \frac{\partial }{{\partial {x_j}}}\Big( {{\widehat C_{ijkl}}(T_{0})\frac{{\partial u_{k0} }}{{\partial {x_l}}}}- {{\widehat \beta _{ij}}(T_{0})( {{T_{0} } - \widetilde T} )}  \Big)= {f_i},\;\;\text{in}\;\;\Omega\times(0,T^{*}),\\
&T_{{0}}(\bm{x},t) = \widehat T(\bm{x},t),\;\;\text{on}\;\;\partial {\Omega_T}\times(0,T^{*}),\\
&\bm{u}_{0}(\bm{x},t)=\widehat{\bm{u}}(\bm{x},t),\;\;\text{on}\;\;\partial\Omega_{u}\times(0,T^{*}),\\
&{{\widehat k_{ij}}(T_{0})\frac{\partial T_{0}(\bm{x},t)}{\partial {x_j}}}{n_i} = \bar q(\bm{x},t),\;\;\text{on}\;\;\partial {\Omega_q}\times(0,T^{*}),\\
&\big[ {{\widehat C_{ijkl}}(T_{0})\frac{{\partial u_{{k0}}(\bm{x},t)}}{{\partial {x_l}}}}- {{\widehat \beta _{ij}}(T_{0})({{T_{{0}}(\bm{x},t)} - \widetilde T} )}  \big]n_j=\bar{\sigma}_i(\bm{x},t),\;\;\text{on}\;\;\partial\Omega_{\sigma}\times(0,T^{*}),\\
&T_{{0}}({\bm{x}},0)=\widetilde T,\;\;\text{in}\;\;\Omega.
\end{aligned} \right.
\end{equation}
Here, the macroscopic homogenized material parameters in (2.18) are defined as follows
\begin{equation}
\begin{aligned}
&\widehat S(T_{0})=\lim_{M\rightarrow +\infty}\frac{\displaystyle\sum_{s=1}^M \hat S(T_{0},\bm{\omega}^s)}{M},\\
&\widehat k_{ij}(T_{0})=\lim_{M\rightarrow +\infty}\frac{\displaystyle\sum_{s=1}^M \hat k_{ij}(T_{0},\bm{\omega}^s)}{M},\\
&\widehat C_{ijkl}(T_{0})=\lim_{M\rightarrow +\infty}\frac{\displaystyle\sum_{s=1}^M \hat C_{ijkl}(T_{0},\bm{\omega}^s)}{M},\\
&\widehat \beta_{ij}(T_{0})=\lim_{M\rightarrow +\infty}\frac{\displaystyle\sum_{s=1}^M \hat \beta_{ij}(T_{0},\bm{\omega}^s)}{M},
\end{aligned}
\end{equation}
where the effective material coefficients in (2.19) are evaluated using the following formulas, which correspond to microscopic unit cell $Y^s$ with random sample $\bm{\omega}^s$
\begin{equation}
\begin{aligned}
&\hat S(T_{0},\bm{\omega}^s) = \frac{1}{|Y^s|}{\int_{Y^s}}{\rho^{(0)}({\bm{y}},{T_{0}},\bm{\omega}^s){c}^{(0)}({\bm{y}},{T_{0}},\bm{\omega}^s)}dY^s,\\
&{\hat k_{ij}}(T_{0},\bm{\omega}^s) = \frac{1}{|Y^s|}{\int_{Y^s}}\big({k_{ij}^{(0)}({\bm{y}},{T_{0}},\bm{\omega}^s) + k_{ik}^{(0)}({\bm{y}},{T_{0}},\bm{\omega}^s){\frac{\partial M_j}{\partial y_k}}}\big)dY^s,\\
&{\hat C_{ijkl}}(T_{0},\bm{\omega}^s) = \frac{1}{|Y^s|}{\int_{Y^s}}\big({C_{ijkl}^{(0)}({\bm{y}},{T_{0}},\bm{\omega}^s)+ C_{ijmn}^{(0)}({\bm{y}},{T_{0}},\bm{\omega}^s){\frac{\partial N_{mk}^{l}}{\partial y_n}}}\big)dY^s,\\
&{\hat \beta _{ij}}(T_{0},\bm{\omega}^s) = \frac{1}{|Y^s|}{\int_{Y^s}}\big( {\beta _{ij}^{(0)}({\bm{y}},{T_{0}},\bm{\omega}^s)+ C_{ijkl}^{(0)}({\bm{y}},{T_{0}},\bm{\omega}^s){\frac{\partial P_k}{\partial y_l}}}\big)dY^s.
\end{aligned}
\end{equation}
\begin{rmk}
It is important to emphasize that all homogenized material parameters vary with the macroscopic homogenized solution $T_{0}$ due to the quasi-periodic properties of first-order cell functions. This discrepancy is significant when compared to linear composites.
\end{rmk}
Now, we proceed to establish the vital second-order correctors $u_{i2}$ and $T_{2}$. By substituting (2.12) and (2.14) into (2.11), and then subtracting (2.11) from (2.18), the following equations are obtained after simplification and computation
\begin{equation}
\left\{\begin{aligned}
&\frac{\partial}{\partial y_i}\Big( {k_{ij}^{(0)}(\bm{y},T_{0},\bm{\omega}){\frac{\partial T_{2}}{\partial y_j}}} \Big)=\Big[ { \rho^{(0)}(\bm{y},T_{0},\bm{\omega}) {c}^{(0)}(\bm{y},T_{0},\bm{\omega})-\widehat S(T_{0}) } \Big]\frac{{\partial {T_{0}}}}{{\partial t}}\\
&+ \Big[ { {\widehat k}_{\alpha_1\alpha_2}(T_{0}) -  {k_{\alpha_1\alpha_2}^{(0)}}(\bm{y},T_{0},\bm{\omega})}{- \frac{\partial}{\partial y_i}\big( {k_{i\alpha_1}^{(0)}(\bm{y},T_{0},\bm{\omega}){M_{\alpha_2}}} \big)-{k_{\alpha_1j}(\bm{y},T_{0},\bm{\omega})\frac{\partial M_{\alpha_2}}{\partial y_j}}} \Big]\frac{\partial^2 T_{0}}{\partial x_{\alpha_1}\partial x_{\alpha_2}}\\
&+ \Big[ { \frac{\partial{\widehat k}_{i\alpha_1}(T_{0})}{\partial x_i} -  {\frac{\partial{k}_{i\alpha_1}^{(0)}(\bm{y},T_{0},\bm{\omega})}{\partial x_i}}- \frac{\partial}{\partial y_i}\big( {k_{ij}^{(0)}(\bm{y},T_{0},\bm{\omega})\frac{\partial M_{\alpha_1}}{\partial x_{j}}} \big)}{-\frac{\partial}{\partial x_i}\big( {k_{ij}^{(0)}(\bm{y},T_{0},\bm{\omega})\frac{\partial M_{\alpha_1}}{\partial y_{j}}} \big)} \Big]\frac{\partial T_{0}}{\partial x_{\alpha_1}}\\
&- \frac{\partial}{\partial y_i}\Big({{M_{\alpha_1}}\mathbf{D}^{(0,1)}{k_{i\alpha_2}^{(0)}(\bm{y},T_{0},\bm{\omega})}}+ {{M_{\alpha_1}}\mathbf{D}^{(0,1)}{k_{ij}^{(0)}(\bm{y},T_{0},\bm{\omega})}\frac{\partial M_{\alpha_2}}{\partial y_j}}\Big)\frac{\partial T_{0}}{\partial x_{\alpha_1}}\frac{\partial T_{0}}{\partial x_{\alpha_2}},\\
&\frac{\partial}{\partial y_j}\Big( {C_{ijkl}^{(0)}(\bm{y},T_{0},\bm{\omega}){\frac{\partial u_{k2}}{\partial y_l}}} \Big)= \Big[ { {\widehat C}_{i\alpha_1m\alpha_2}(T_{0}) -  {C_{i\alpha_1m\alpha_2}^{(0)}}(\bm{y},T_{0},\bm{\omega})- \frac{\partial}{\partial y_j}\big( {C_{ijk\alpha_1}^{(0)}(\bm{y},T_{0},\bm{\omega}){N^{\alpha_2}_{km}}} \big)}\\
&{-{C_{i\alpha_1kj}^{(0)}(\bm{y},T_{0},\bm{\omega})\frac{\partial N^{\alpha_2}_{km}}{\partial y_j}}} \Big]\frac{\partial^2 u_{m0}}{\partial x_{\alpha_1}\partial x_{\alpha_2}}+ \Big[ { \frac{\partial{\widehat C}_{ijm\alpha_1}(T_{0})}{\partial x_j} -  {\frac{\partial{C}_{ijm\alpha_1}^{(0)}(\bm{y},T_{0},\bm{\omega})}{\partial x_j}}}\\
&{- \frac{\partial}{\partial y_j}\big( {C_{ijkl}^{(0)}(\bm{y},T_{0},\bm{\omega})\frac{\partial N^{\alpha_1}_{km}}{\partial x_{l}}} \big)}{-\frac{\partial}{\partial x_j}\big( {C_{ijkl}^{(0)}(\bm{y},T_{0},\bm{\omega})\frac{\partial N^{\alpha_1}_{km}}{\partial y_{l}}} \big)} \Big]\frac{\partial u_{m0}}{\partial x_{\alpha_1}}\\
&-\Big[ { \frac{\partial{\widehat \beta}_{ij}(T_{0})}{\partial x_j} -  {\frac{\partial{\beta}_{ij}^{(0)}(\bm{y},T_{0},\bm{\omega})}{\partial x_j}}- \frac{\partial}{\partial y_j}\big( {C_{ijkl}^{(0)}(\bm{y},T_{0},\bm{\omega})\frac{\partial P_{k}}{\partial x_{l}}} \big)}\\
&{-\frac{\partial}{\partial x_j}\big( {C_{ijkl}^{(0)}(\bm{y},T_{0},\bm{\omega})\frac{\partial P_{k}}{\partial y_{l}}} \big)} \Big]({T_{0}} - \widetilde T)-\big[ {{\widehat \beta }_{i\alpha_1}(T_{0})} - {\beta _{i\alpha_1}^{(0)}(\bm{y},T_{0},\bm{\omega})}\\
&{ -C_{i\alpha_1kl}^{(0)}(\bm{y},T_{0},\bm{\omega}){\frac{\partial P_k}{\partial y_l}}}{- \frac{\partial}{\partial y_j}( {C_{ijk\alpha_1}^{(0)}(\bm{y},T_{0},\bm{\omega}){P_k}})-{\frac{\partial}{\partial y_j}}( {\beta _{ij}^{(0)}(\bm{y},T_{0},\bm{\omega}){M_{\alpha_1}}} )} \big]\frac{\partial T_{0}}{\partial x_{\alpha_1}}\\
&+ \frac{\partial}{\partial y_j}\Big( {{M_{\alpha_1}}\mathbf{D}^{(0,1)}\beta _{ij}^{(0)}(\bm{y},T_{0},\bm{\omega}) + {M_{\alpha_1}}\mathbf{D}^{(0,1)}C_{i jkl}^{(0)}(\bm{y},T_{0},\bm{\omega})\frac{\partial P_k}{\partial y_l}} \Big)\frac{\partial T_{0}}{\partial x_{\alpha_1}}({{T_{0}}-\widetilde T})\\
&- \frac{\partial}{\partial y_j}\Big( {{M_{\alpha_1}}\mathbf{D}^{(0,1)}C_{ijm\alpha_2}^{(0)}(\bm{y},T_{0},\bm{\omega}) + {M_{\alpha_1}}\mathbf{D}^{(0,1)}C_{ijkl}^{(0)}(\bm{y},T_{0},\bm{\omega})\frac{\partial N^{\alpha_2}_{km}}{\partial y_l}} \Big)\frac{\partial T_{0}}{\partial x_{\alpha_1}}\frac{\partial u_{m0}}{\partial x_{\alpha_2}}.
\end{aligned}\right.
\end{equation}
According to equations (2.20), then we construct the concrete expressions for $u_{i2}$ and ${T_{2}}$ as follows
\begin{equation}
\left\{
\begin{aligned}
{T_{2}}(\bm{x},\bm{y},t,\bm{\omega})&= S({\bm{y}},{T_{0}},\bm{\omega})\frac{{\partial {T_{0}}}}{{\partial t}}+{M_{\alpha_1\alpha_2}}({\bm{y}},{T_{0}},\bm{\omega})\frac{\partial^2 T_{0}}{\partial x_{\alpha_1}\partial x_{\alpha_2}}+R_{\alpha_1}({\bm{y}},{T_{0}},\bm{\omega})\frac{\partial T_{0}}{\partial x_{\alpha_1}}\\&-{B_{\alpha_1\alpha_2}}({\bm{y}},{T_{0}},\bm{\omega})\frac{\partial T_{0}}{\partial x_{\alpha_1}}\frac{\partial T_{0}}{\partial x_{\alpha_2}},\\
u_{i2}(\bm{x},\bm{y},t,\bm{\omega})&= N_{im}^{\alpha_1\alpha_2}({\bm{y}},{T_{0}},\bm{\omega})\frac{\partial^2 u_{m0}}{\partial x_{\alpha_1}\partial x_{\alpha_2}}+ Z_{im}^{\alpha_1}({\bm{y}},{T_{0}},\bm{\omega})\frac{\partial u_{m0}}{\partial x_{\alpha_1}} \\
&-{Q_i}({\bm{y}},{T_{0}},\bm{\omega})( {{T_{0}} - \widetilde T} ) - H_i^{\alpha_1}({\bm{y}},{T_{0}},\bm{\omega})\frac{\partial T_{0}}{\partial x_{\alpha_1}}\\
&+W_i^{\alpha_1}({\bm{y}},{T_{0}},\bm{\omega})\frac{\partial T_{0}}{\partial x_{\alpha_1}}({{T_{0}}-\widetilde T})-A_{im}^{\alpha_1\alpha_2}({\bm{y}},{T_{0}},\bm{\omega})\frac{\partial T_{0}}{\partial x_{\alpha_1}}\frac{\partial u_{m0}}{\partial x_{\alpha_2}},\;\;\alpha_1,\alpha_2,m,n=1,2,3.
\end{aligned} \right.
\end{equation}
In the above expressions, $S$, $M_{\alpha_1\alpha_2}$, $R_{\alpha_1}$, ${B_{\alpha_1\alpha_2}}$, $N_{im}^{\alpha_1\alpha_2}$, $Z_{im}^{\alpha_1}$, ${Q_i}$, $H_i^{\alpha_1}$, $W_i^{\alpha_1}$ and $A_{im}^{\alpha_1\alpha_2}$ are $1-$periodic functions defined in microscopic unit cell $Y^s$ for any random sample $\omega^s$, which are referred to as second-order auxiliary cell functions. By substituting (2.21) into (2.20), a series of equations, which are attached with the homogeneous Dirichlet boundary condition, are derived as follows
\begin{equation}
\left\{
\begin{aligned}
&\frac{\partial}{\partial y_i}\big[ { k_{ij}^{(0)}({\bm{y}},{T_{0}},\bm{\omega}^s){\frac{\partial S}{\partial y_j}}} \big] = { \rho^{(0)} {c}^{(0)}-\widehat S },&\;\;\;&\bm{y}\in Y^s,\\
&S({\bm{y}},{T_{0}},\bm{\omega}^s)=0,&\;\;\;&\bm{y}\in\partial Y^s.
\end{aligned} \right.
\end{equation}
\begin{equation}
\left\{
\begin{aligned}
&\frac{\partial}{\partial y_i}\big[ { k_{ij}^{(0)}({\bm{y}},{T_{0}},\bm{\omega}^s){\frac{\partial {M_{\alpha_1\alpha_2}}}{\partial y_j}}} \big] = { {\widehat k}_{\alpha_1\alpha_2} -  {k_{\alpha_1\alpha_2}^{(0)}}- \frac{\partial}{\partial y_i}\big( {k_{i\alpha_1}^{(0)}{M_{\alpha_2}}} \big)-{k_{\alpha_1j}^{(0)}\frac{\partial M_{\alpha_2}}{\partial y_j}}},&\;\;\;&\bm{y}\in Y^s,\\
&{M_{\alpha_1\alpha_2}}({\bm{y}},{T_{0}},\bm{\omega}^s)=0,&\;\;\;&\bm{y}\in\partial Y^s.
\end{aligned} \right.
\end{equation}
\begin{equation}
\left\{
\begin{aligned}
&\frac{\partial}{\partial y_i}\big[ { k_{ij}^{(0)}({\bm{y}},{T_{0}},\bm{\omega}^s){\frac{\partial R_{\alpha_1}}{\partial y_j}}} \big] ={ \frac{\partial{\widehat k}_{i\alpha_1}}{\partial x_i} -  {\frac{\partial{k}_{i\alpha_1}^{(0)}}{\partial x_i}}- \frac{\partial}{\partial y_i}\big( {k_{ij}^{(0)}\frac{\partial M_{\alpha_1}}{\partial x_{j}}} \big)-\frac{\partial}{\partial x_i}\big( {k_{ij}^{(0)}\frac{\partial M_{\alpha_1}}{\partial y_{j}}} \big)},&\;\;\;&\bm{y}\in Y^s,\\
&R_{\alpha_1}({\bm{y}},{T_{0}},\bm{\omega}^s)=0,&\;\;\;&\bm{y}\in\partial Y^s.
\end{aligned} \right.
\end{equation}
\begin{equation}
\left\{
\begin{aligned}
&\frac{\partial}{\partial y_i}\big[ { k_{ij}^{(0)}({\bm{y}},{T_{0}},\bm{\omega}^s){\frac{\partial B_{\alpha_1\alpha_2}}{\partial y_j}}} \big] =\frac{\partial}{\partial y_i}\Big({{M_{\alpha_1}}\mathbf{D}^{(0,1)}{k_{i\alpha_2}^{(0)}}}+ {{M_{\alpha_1}}\mathbf{D}^{(0,1)}{k_{ij}^{(0)}}\frac{\partial M_{\alpha_2}}{\partial y_j}}\Big),&\;\;\;&\bm{y}\in Y^s,\\
&B_{\alpha_1\alpha_2}({\bm{y}},{T_{0}},\bm{\omega}^s)=0,&\;\;\;&\bm{y}\in\partial Y^s.
\end{aligned} \right.
\end{equation}
\begin{equation}
\left\{
\begin{aligned}
&\frac{\partial}{\partial y_j}\big[ C_{ijkl}^{(0)}({\bm{y}},{T_{0}},\bm{\omega}^s){\frac{\partial N_{km}^{\alpha_1\alpha_2}}{\partial y_l}}} \big] = { {\widehat C}_{i\alpha_1m\alpha_2} -  {C_{i\alpha_1m\alpha_2}^{(0)}}\\
&- \frac{\partial}{\partial y_j}\big( {C_{ijk\alpha_1}^{(0)}{N^{\alpha_2}_{km}}} \big)-{C_{i\alpha_1kj}^{(0)}\frac{\partial N^{\alpha_2}_{km}}{\partial y_j}},&&\bm{y}\in Y^s,\\
&N_{km}^{\alpha_1\alpha_2}({\bm{y}},{T_{0}},\bm{\omega}^s)= 0,&&\bm{y}\in\partial Y^s.
\end{aligned} \right.
\end{equation}
\begin{equation}
\left\{
\begin{aligned}
&\frac{\partial}{\partial y_j}\big[ { C_{ijkl}^{(0)}({\bm{y}},{T_{0}},\bm{\omega}^s){\frac{\partial Z_{km}^{\alpha_1}}{\partial y_l}}} \big] ={ \frac{\partial{\widehat C}_{ijm\alpha_1}}{\partial x_j} -  {\frac{\partial{C}_{ijm\alpha_1}^{(0)}}{\partial x_j}}}\\
&{- \frac{\partial}{\partial y_j}\big( {C_{ijkl}^{(0)}\frac{\partial N^{\alpha_1}_{km}}{\partial x_{l}}} \big)-\frac{\partial}{\partial x_j}\big( {C_{ijkl}^{(0)}\frac{\partial N^{\alpha_1}_{km}}{\partial y_{l}}} \big)},&&\bm{y}\in Y^s,\\
&Z_{km}^{\alpha_1}({\bm{y}},{T_{0}},\bm{\omega}^s)= 0,&&\bm{y}\in\partial Y^s.
\end{aligned} \right.
\end{equation}
\begin{equation}
\left\{
\begin{aligned}
&\frac{\partial}{\partial y_j}\big[ { C_{ijkl}^{(0)}({\bm{y}},{T_{0}},\bm{\omega}^s){\frac{\partial Q_{k}}{\partial y_l}}} \big] ={ \frac{\partial{\widehat \beta}_{ij}}{\partial x_j} -  {\frac{\partial{\beta}_{ij}^{(0)}}{\partial x_j}}- \frac{\partial}{\partial y_j}\big( {C_{ijkl}^{(0)}\frac{\partial P_{k}}{\partial x_{l}}} \big)-\frac{\partial}{\partial x_j}\big( {C_{ijkl}^{(0)}\frac{\partial P_{k}}{\partial y_{l}}} \big)},&&\bm{y}\in Y^s,\\
&Q_{k}({\bm{y}},{T_{0}},\bm{\omega}^s)= 0,&&\bm{y}\in\partial Y^s.
\end{aligned} \right.
\end{equation}
\begin{equation}
\left\{
\begin{aligned}
&\frac{\partial}{\partial y_j}\big[ { C_{ijkl}^{(0)}({\bm{y}},{T_{0}},\bm{\omega}^s){\frac{\partial H_{k}^{\alpha_1}}{\partial y_l}}} \big] = {{\widehat \beta }_{i\alpha_1}} - {\beta _{i\alpha_1}^{(0)} -C_{i\alpha_1kl}^{(0)}{\frac{\partial P_k}{\partial y_l}}}\\
&{- \frac{\partial}{\partial y_j}( {C_{ijk\alpha_1}^{(0)}{P_k}})-{\frac{\partial}{\partial y_j}}( {\beta _{ij}^{(0)}{M_{\alpha_1}}} )},&&\bm{y}\in Y^s,\\
&{H_{k}^{\alpha_1}({\bm{y}},{T_{0}},\bm{\omega}^s)= 0},&&\bm{y}\in\partial Y^s.
\end{aligned} \right.
\end{equation}
\begin{equation}
\left\{
\begin{aligned}
&\frac{\partial}{\partial y_j}\big[ { C_{ijkl}^{(0)}({\bm{y}},{T_{0}},\bm{\omega}^s){\frac{\partial W_{k}^{\alpha_1}}{\partial y_l}}} \big]  = \frac{\partial}{\partial y_j}\Big( {{M_{\alpha_1}}\mathbf{D}^{(0,1)}\beta _{i j}^{(0)} + {M_{\alpha_1}}\mathbf{D}^{(0,1)}C_{i jkl}^{(0)}\frac{\partial P_k}{\partial y_l}} \Big),&\;&\bm{y}\in Y^s,\\
&{W_{k}^{\alpha_1}({\bm{y}},{T_{0}},\bm{\omega}^s)= 0},&\;&\bm{y}\in\partial Y^s.
\end{aligned} \right.
\end{equation}
\begin{equation}
\left\{
\begin{aligned}
&\frac{\partial}{\partial y_j}\big[ { C_{ijkl}^{(0)}({\bm{y}},{T_{0}},\bm{\omega}^s){\frac{\partial A_{km}^{\alpha_1\alpha_2}}{\partial y_l}}} \big]  = \frac{\partial}{\partial y_j}\Big( {{M_{\alpha_1}}\mathbf{D}^{(0,1)}C_{ijm\alpha_2}^{(0)} + {M_{\alpha_1}}\mathbf{D}^{(0,1)}C_{ijkl}^{(0)}\frac{\partial N^{\alpha_2}_{km}}{\partial y_l}} \Big),&\;&\bm{y}\in Y^s,\\
&{A_{km}^{\alpha_1\alpha_2}({\bm{y}},{T_{0}},\bm{\omega}^s)= 0},&\;&\bm{y}\in\partial Y^s.
\end{aligned} \right.
\end{equation}

Summing up, the macro-micro coupled SHOMS asymptotic solutions of the multi-scale nonlinear dynamic thermo-mechanical problem (2.1) are established as follows
\begin{equation}
{\begin{aligned}
{\hat T^\varepsilon }(\bm{x},t,\bm{\omega})&=
T_{0}(\bm{x},t)+\varepsilon{M_{\alpha_1}}({\bm{y}},{T_{0}},\bm{\omega})\frac{\partial T_{0}(\bm{x},t)}{\partial x_{\alpha_1}}\\
&+ \varepsilon^2\big[S({\bm{y}},{T_{0}},\bm{\omega})\frac{{\partial {T_{0}(\bm{x},t)}}}{{\partial t}}+{M_{\alpha_1\alpha_2}}({\bm{y}},{T_{0}},\bm{\omega})\frac{\partial^2 T_{0}(\bm{x},t)}{\partial x_{\alpha_1}\partial x_{\alpha_2}}\\
&+R_{\alpha_1}({\bm{y}},{T_{0}},\bm{\omega})\frac{\partial T_{0}(\bm{x},t)}{\partial x_{\alpha_1}}-{B_{\alpha_1\alpha_2}}({\bm{y}},{T_{0}},\bm{\omega})\frac{\partial T_{0}(\bm{x},t)}{\partial x_{\alpha_1}}\frac{\partial T_{0}(\bm{x},t)}{\partial x_{\alpha_2}}\big].
\end{aligned}}
\end{equation}
\begin{equation}
{\begin{aligned}
\hat u_i^\varepsilon(\bm{x},t,\bm{\omega})&=
u_{i0}(\bm{x},t)+\varepsilon\big[N_{im}^{\alpha_1}({\bm{y}},{T_{0}},\bm{\omega})\frac{\partial u_{m0}(\bm{x},t)}{\partial x_{\alpha_1}} - {P_i}({\bm{y}},{T_{0}},\bm{\omega})({T_{0}(\bm{x},t)} - \widetilde T)\big]\\
&+ \varepsilon^2\big[ N_{im}^{\alpha_1\alpha_2}({\bm{y}},{T_{0}},\bm{\omega})\frac{\partial^2 u_{m0}(\bm{x},t)}{\partial x_{\alpha_1}\partial x_{\alpha_2}}+ Z_{im}^{\alpha_1}({\bm{y}},{T_{0}},\bm{\omega})\frac{\partial u_{m0}(\bm{x},t)}{\partial x_{\alpha_1}}\\
&-{Q_i}({\bm{y}},{T_{0}},\bm{\omega})( {{T_{0}(\bm{x},t)} - \widetilde T}) - H_i^{\alpha_1}({\bm{y}},{T_{0}},\bm{\omega})\frac{\partial T_{0}(\bm{x},t)}{\partial x_{\alpha_1}}\\
&+W_i^{\alpha_1}({\bm{y}},{T_{0}},\bm{\omega})\frac{\partial T_{0}(\bm{x},t)}{\partial x_{\alpha_1}}({{T_{0}(\bm{x},t)}-\widetilde T})-A_{im}^{\alpha_1\alpha_2}({\bm{y}},{T_{0}},\bm{\omega})\frac{\partial T_{0}(\bm{x},t)}{\partial x_{\alpha_1}}\frac{\partial u_{m0}(\bm{x},t)}{\partial x_{\alpha_2}}\big].
\end{aligned}}
\end{equation}
\begin{rmk}
It is worth emphasizing that the proposed macro-micro coupled SHOMS asymptotic solutions of the multi-scale nonlinear dynamic thermo-mechanical problem (2.1) can be reasonably degenerated to the SHOMS asymptotic solutions of the multi-scale linear dynamic thermo-mechanical problem. Since the material parameters of multi-scale linear dynamic thermo-mechanical problem are independent with temperature variable $T^{\varepsilon}$, the second-order cell functions $R_{\alpha_1}$, ${B_{\alpha_1\alpha_2}}$ $Z_{im}^{\alpha_1}$, ${Q_i}$, $W_i^{\alpha_1}$ and $A_{im}^{\alpha_1\alpha_2}$ will be degenerated to zero and disappear in formulas (2.34) and (2.35).
\end{rmk}

Based on the theoretical analysis mentioned above, the macro-micro coupled SHOMS computational model has been established. This computational model includes first-order and second-order auxiliary cell model in micro-scale, homogenized model in macro-scale and macro-micro coupled SHOMS solutions. Furthermore, by respecting the chain rule (2.2) and the macro-micro asymptotic formulas (2.34) and (2.35), the heat fluxes, strains and stresses and inside $\Omega$ of random composite materials with temperature-dependent properties are approximately evaluated as follows
\begin{equation}
\hat q_i^\varepsilon({\bm{x}},t,\bm{\omega}) = -k_{ij}^\varepsilon({\bm{x}},{T^\varepsilon },\bm{\omega})\frac{\partial \hat T^{\varepsilon}({\bm{x}},t,\bm{\omega})}{\partial x_j}.
\end{equation}
\begin{equation}
\hat\varepsilon_{ij}^\varepsilon({\bm{x}},t,\bm{\omega})=\frac{1}{2}\big[\frac{\partial \hat u_i^{\varepsilon}({\bm{x}},t,\bm{\omega})}{\partial x_j}+\frac{\partial \hat u_j^{\varepsilon}({\bm{x}},t,\bm{\omega})}{\partial x_i}\big].
\end{equation}
and
\begin{equation}
\begin{aligned}
&\hat\sigma _{ij}^\varepsilon({\bm{x}},t,\bm{\omega}) = {C_{ijkl}^\varepsilon}({\bm{x}},{T^\varepsilon },\bm{\omega})\frac{\partial \hat u_k^{\varepsilon}({\bm{x}},t,\bm{\omega})}{\partial x_l}- {\beta _{ij}^\varepsilon}({\bm{x}},{T^\varepsilon },\bm{\omega})( {{\hat T^{{\varepsilon}}}({\bm{x}},t,\bm{\omega}) - \widetilde T} ).
\end{aligned}
\end{equation}
\subsection{Local error analysis of statistical multi-scale solutions in point-wise sense}
In this subsection, we begin by introducing the multi-scale approximate solutions of stochastic multi-scale problem (2.1) as follows
\begin{equation}
\begin{aligned}
&T^{(1\varepsilon)}({\bm{x}},t,\bm{\omega})=T_{0}+\varepsilon T_{1},\;\;u^{(1\varepsilon)}_{i}({\bm{x}},t,\bm{\omega})=u_{i0}+\varepsilon u_{i1},\\
&T^{(2\varepsilon)}({\bm{x}},t,\bm{\omega})=T_{0}+\varepsilon T_{1}+\varepsilon^2 T_{2},\;\;u^{(2\varepsilon)}_{i}({\bm{x}},t,\bm{\omega})=u_{i0}+\varepsilon u_{i1}+\varepsilon^2 u_{i2},
\end{aligned}
\end{equation}
where we defined $T_{0}$ and $u_{i0}$ as the macroscopic homogenized solutions, $T^{(1\varepsilon)}$ and $u^{(1\varepsilon)}_{i}$ as the statistical lower-order multi-scale (SLOMS) solutions, $T^{(2\varepsilon)}$ and $u^{(2\varepsilon)}_{i}$ as the SHOMS approximate solutions, which equate to $\hat T^\varepsilon$ in (2.33) and $\hat u_i^\varepsilon$ in (2.34), respectively.

Next, we introduce the following residual functions for the purpose of numerical accuracy analysis.
\begin{align}
{\begin{array}{*{20}{l}}
T^{(1\varepsilon)}_{\Delta}({\bm{x}},t,\bm{\omega})=T^\varepsilon-T^{(1\varepsilon)},\;\;u^{(1\varepsilon)}_{\Delta i}({\bm{x}},t,\bm{\omega})=u^\varepsilon_{i}-u^{(1\varepsilon)}_{i},\\
T^{(2\varepsilon)}_{\Delta}({\bm{x}},t,\bm{\omega})=T^\varepsilon-T^{(2\varepsilon)},\;\;u^{(2\varepsilon)}_{\Delta i}({\bm{x}},t,\bm{\omega})=u^\varepsilon_{ i}-u^{(2\varepsilon)}_{i}.
\end{array}}
\end{align}

After that, by inserting $T^{(1\varepsilon)}_{\Delta}$ and $u^{(1\varepsilon)}_{\Delta i}$ from (2.39) into (2.1), the following residual equations of SLOMS solutions are obtained, which hold in the distribution sense.
\begin{equation}
\left\{ \begin{aligned}
&\rho^{\varepsilon} ({\bm{x}},{T^\varepsilon },\bm{\omega})c^{\varepsilon}({\bm{x}},{T^\varepsilon },\bm{\omega})\frac{{\partial {T^{(1\varepsilon)}_{\Delta}}}}{{\partial t}}- \frac{\partial }{{\partial {x_i}}}\Big( {{k_{ij}^{\varepsilon}}({\bm{x}},{T^\varepsilon },\bm{\omega})\frac{{\partial {T^{(1\varepsilon)}_{\Delta}}}}{{\partial {x_j}}}}\Big) \\
&\quad\quad\quad=\mathbb{F}_{0}(\bm{x},\bm{y},t,\bm{\omega})+\varepsilon \mathbb{F}_1(\bm{x},\bm{y},t,\bm{\omega}),\;\;\text{in}\;\;\Omega\times(0,T^{*}),\\
&- \frac{\partial }{{\partial {x_j}}}\Big( {{C_{ijkl}^{\varepsilon}}({\bm{x}},{T^\varepsilon },\bm{\omega})\frac{{\partial u^{(1\varepsilon)}_{\Delta k}}}{{\partial {x_l}}}}- {{\beta _{ij}^{\varepsilon}}({\bm{x}},{T^\varepsilon },\bm{\omega})( {{T^{(1\varepsilon)}_{\Delta}} - \widetilde T} )}  \Big) \\
&\quad\quad\quad=\mathbb{S}_{0i}(\bm{x},\bm{y},t,\bm{\omega})+\varepsilon \mathbb{S}_{1i}(\bm{x},\bm{y},t,\bm{\omega}),\;\;\text{in}\;\;\Omega\times(0,T^{*}).
\end{aligned} \right.
\end{equation}
Here, the detailed forms of $\mathbb{S}_{0i}(\bm{x},\bm{y},t,\bm{\omega})$, $\mathbb{S}_{1i}(\bm{x},\bm{y},t,\bm{\omega})$, $\mathbb{F}_{0}(\bm{x},\bm{y},t,\bm{\omega})$ and $\mathbb{F}_{1}(\bm{x},\bm{y},t,\bm{\omega})$ can be obtained easily by referring to the results in \cite{R20}. However, due to their lengthy expressions, they are not exhibited in the current study.

Furthermore, by substituting $T^{(2\varepsilon)}_{\Delta}$ and $u^{(2\varepsilon)}_{\Delta i}$ from (2.39) into (2.1), the resulting residual equations of SHOMS solutions are derived, which also hold in the distribution sense.
\begin{equation}
\left\{\begin{aligned}
&\rho^{\varepsilon} ({\bm{x}},{T^\varepsilon },\bm{\omega})c^{\varepsilon}({\bm{x}},{T^\varepsilon },\bm{\omega})\frac{{\partial {T^{(2\varepsilon)}_{\Delta}}}}{{\partial t}}- \frac{\partial }{{\partial {x_i}}}\Big( {{k_{ij}^{\varepsilon}}({\bm{x}},{T^\varepsilon },\bm{\omega})\frac{{\partial {T^{(2\varepsilon)}_{\Delta}}}}{{\partial {x_j}}}}\Big)\\
&\quad\quad\quad=\varepsilon \mathbb{G}(\bm{x},\bm{y},t,\bm{\omega}),\;\;\text{in}\;\;\Omega\times(0,T^{*}),\\
&- \frac{\partial }{{\partial {x_j}}}\Big( {{C_{ijkl}^{\varepsilon}}({\bm{x}},{T^\varepsilon },\bm{\omega})\frac{{\partial u^{(2\varepsilon)}_{\Delta k}}}{{\partial {x_l}}}}- {{\beta _{ij}^{\varepsilon}}({\bm{x}},{T^\varepsilon },\bm{\omega})( {{T^{(2\varepsilon)}_{\Delta}} - \widetilde T} )}  \Big) \\
&\quad\quad\quad=\varepsilon \mathbb{H}_{i}(\bm{x},\bm{y},t,\bm{\omega}),\;\;\text{in}\;\;\Omega\times(0,T^{*}).
\end{aligned} \right.
\end{equation}
In the above residual equations, the specific forms of $\mathbb{H}_{i}(\bm{x},\bm{y},t,\bm{\omega})$ and $\mathbb{G}(\bm{x},\bm{y},t,\bm{\omega})$ are also not displayed in the present study owing to the same reason mentioned above.

In conclusion, the numerical accuracy analysis in the point-wise sense is presented, focusing on the residual equations (2.40) and (2.41). Noting the residual equations (2.40), it can be concluded that the residual error of SLOMS solutions is of order-$O(1)$ in the point-wise sense, primarily due to the presence of terms $\mathbb{S}_{0i}({\bm{x}},{\bm{y}},t,\bm{\omega})$ and $\mathbb{F}_0({\bm{x}},{\bm{y}},t,\bm{\omega})$. Comparing with the residual equations (2.41), the residual error of SHOMS solutions is of order-$O(\varepsilon)$ in the point-wise sense. This implies that the SHOMS solutions can satisfy the local energy conservation of thermal equation and local momentum conservation of mechanical equations of the original stochastic multi-scale equations (2.1) in the point-wise sense. Thus even $\varepsilon$ is a small constant, the SHOMS solutions can still provide the required accuracy for engineering computation and accurately capture the microscopic oscillating behaviors exhibited by random composite materials. This is the fundamental rationale and motivation for us to establish higher-order corrected terms and develop the SHOMS solutions.
\section{Space-time multi-scale numerical algorithm}
In this section, an innovative space-time multi-scale numerical algorithm with off-line and on-line stages is presented for the stochastic nonlinear governing equations (2.1), which is based on the FEM in spatial region and the finite difference method (FDM) in time direction. The detailed algorithm procedures for stochastic multi-scale nonlinear problem (2.1) are listed as follows.
\subsection{Off-line computation for microscopic cell problems}
\begin{enumerate}[(1)]
\item
Generate a random sample $\omega^s$ for microscopic unit cell $Y^s$ according to a given probability distribution model $\bar P$ and the volume fractions. And then, partition $Y^s$ into finite element meshes $J^{h_1}=\{K\}$ where $h_1=$max$_K\{h_K\}$ and define the linear conforming finite element space $V_{h_1}(Y^s)=\{\nu\in C^0(\bar{Y}^s):\nu\mid_{\partial Y^s}=0,\nu\mid_{K}\in P_1(K)\}\subset H_0^1(Y^s)$.
\item
Define service temperature range $[T_{min},T_{max}]$ of investigated random composite materials and choose a certain number of representative macroscopic parameters $\bar T_{s}$ in service temperature range where $s=1,2,\cdots,\bar{L}(\bar{L}\in \bf{N}^{*})$. Then, solve the first-order auxiliary cell problems (2.15)-(2.17) with given material properties on $V_{h_1}(Y^s)$ corresponding to distinct representative macroscopic parameters $\bar T_{s}$. The specific FEM scheme for auxiliary cell problem (2.15) is established as follows
\begin{equation}
\int_{Y^s}{k}_{ij}^{(0)}({\bm{y}},{T_{0}},\bm{\omega}^s)\frac{{\partial {H_{{\alpha _1}}}({\bm{y}},{T_{0}},\bm{\omega}^s)}}{{\partial {y_j}}}\frac{\partial \upsilon^{h_1}}{\partial y_i}dY^s=-\int_{Y^s} {k}_{i\alpha_1}^{(0)}({\bm{y}},{T_{0}},\bm{\omega}^s)\frac{\partial \upsilon^{h_1}}{\partial y_i}dY^s,\;\forall\upsilon^{h_1}\in V_{h_1}(Y^s).
\end{equation}
Other first-order auxiliary cell problems can be solved similarly.
\item
Repeat Steps 1 and 2 for different random samples $\omega^s$ $(s=1,\cdots,M)$ and compute macroscopic homogenized material parameters $\widehat{S}(T_{0})$, $\widehat{k}_{ij}(T_{0})$, $\widehat{C}_{ijkl}(T_{0})$ and $\widehat{\beta} _{ij}(T_{0})$ by using the formula (2.19) corresponding to different macroscopic parameters $\bar T_{s}$.
\item
Employing the same mesh as first-order auxiliary cell functions, the second-order auxiliary cell functions are evaluated by solving second-order auxiliary cell problems (2.23)-(2.32) with given material properties on $V_{h_1}(Y^s)$ respectively, which correspond to different macroscopic parameters $\bar T_{s}$.
\end{enumerate}
\subsection{On-line computation for macroscopic homogenized problem and statistical higher-order multi-scale solutions}
\begin{enumerate}[(1)]
\item
Define the homogenized macroscopic region $\Omega$ in $\mathbb{R}^\mathcal{N}$ and partition $Y^s$ into finite element meshes $J^{h_0}=\{e\}$ where $h_0=$max$_e\{h_e\}$. Then, define the linear conforming finite element spaces  $V_{h_0}(\Omega)=\{\nu\in C^0(\bar{\Omega}):\nu\mid_{\partial\Omega_{T}}=0,\nu\mid_{e}\in P_1(e)\}\subset H^1(\Omega)$ and $V_{h_0}^{*}(\Omega)=\{\nu\in C^0(\bar{\Omega}):\nu\mid_{\partial\Omega_{u}}=0,\nu\mid_{e}\in P_1(e)\}\subset H^1(\Omega)$ for simulating temperature and displacement fields, respectively.
\item
Evaluate macroscopic homogenized material parameters on each nodes of $V_{h_0}(\Omega)$ and $V_{h_0}^{*}(\Omega)$ by  interpolation method. After that, simulate the macroscopic homogenized problem (2.18) by coupled FDM-FEM method in a coarse mesh and with a large time step on the whole domain $\Omega\times(0,T^{*})$. Employing the uniform time step $\displaystyle\Delta t={T^{*}}/{M}$ to discretize time-domain $(0,T^{*})$ as $0=t_0<t_1<\cdots<t_M=T^{*}$ and $t_N=N\Delta t (N=0,\cdots,M)$, then we define $f_{i}^{N}=f_i(\bm{x},t_N)$. In order to preserve the unconditional stability of our multi-scale numerical algorithm, the implicit FDM scheme is adopted in time domain and FEM scheme is employed in spatial domain. The detailed FDM-FEM scheme for macroscopic homogenized problem (2.18) is presented as below
\begin{equation}
\left\{ \begin{aligned}
&\int_{\Omega}\widehat{S}(T_0^{N+1})\frac{{ T_0^{N+1}-T_0^{N}}}{{\Delta t}}\widetilde{\varphi}^{h_0}d\Omega+\int_{\Omega}{\widehat{k}_{ij}(T_0^{N+1})\frac{\partial T_0^{N+1}}{\partial x_j} }\frac{\partial \widetilde{\varphi}^{h_0}}{\partial x_i}d\Omega\\
&\quad\quad\quad=\int_{\Omega} h^{N+1}\widetilde{\varphi}^{h_0}d\Omega+\int_{\partial\Omega_q}\bar{q}^{N+1}\widetilde{\varphi}^{h_0}ds,\;\forall\widetilde{\varphi}^{h_0}\in V_{h_0}(\Omega),\\
&T_{0}^{N}=\widehat{T}(\bm{x},t_N),\;\;\;\text{on}\;\;\partial\Omega_{T}.
\end{aligned} \right.
\end{equation}
\begin{equation}
\left\{ \begin{aligned}
&\int_{\Omega}{\widehat{C}_{ijkl}(T_0^{N+1})\frac{\partial u_{k0}^{N+1}}{\partial x_l}}\frac{\partial \nu_i^{h_0} }{\partial x_j}d\Omega-\int_{\Omega} {\widehat{\beta}_{ij}}(T_0^{N+1}){ (T_0^{N+1}-\widetilde T)}\frac{\partial \nu_i^{h_0} }{\partial x_j}d\Omega\\
&\quad\quad\quad=\int_{\Omega}{ f_{i}^{N+1}}\nu_i^{h_0}d\Omega+{\int_{\partial\Omega_{\sigma}}}\bar{\sigma}_i^{N+1}\nu_i^{h_0}ds,\;\forall\bm{\nu}^{h_0}\in \big(V_{h_0}^{*}(\Omega)\big)^\mathcal{N},\\
&\bm{u}_{0}^{N}=\widehat{\bm{u}}(\bm{x},t_N),\;\;\;\text{on}\;\;\partial\Omega_{u},\\
\end{aligned} \right.
\end{equation}
As a result, we can solve the macroscopic temperature and displacement fields at each time step via (3.43) and (3.44) by turn.
\item
It should be highlighted that, decoupled thermal equation (3.43) is a nonlinear system which can not be computed directly. Herein, the following direct iterative method is presented for simulating the nonlinear system (3.43).
\begin{enumerate}[Step 1:]
\item
Let $\breve T_0(\bm{x})$ be the initial function, and $\breve T_\lambda(\bm{x})$ be the solution at $\lambda$-th iterative step, $\lambda\geq 1$. Set the iteration threshold as $E_{tol}$, and begin iterating.
\item
At the $\lambda$-th iteration step, apply $\breve T_{\lambda-1}(\bm{x})$ to linearize the nonlinear system (3.43) as follows.
\begin{equation}
\left\{ \begin{aligned}
&\int_{\Omega}\widehat{S}(\breve T_{\lambda-1})\frac{{\breve T_{\lambda}-T_0^{N}}}{{\Delta t}}\widetilde{\varphi}^{h_0}d\Omega+\int_{\Omega}{\widehat{k}_{ij}(\breve T_{\lambda-1})\frac{\partial \breve  T_{\lambda}}{\partial x_j}}\frac{\partial \widetilde{\varphi}^{h_0}}{\partial x_i}d\Omega\\
&=\int_{\Omega} h^{N+1}\widetilde{\varphi}^{h_0}d\Omega+\int_{\partial\Omega_q}\bar{q}^{N+1}\widetilde{\varphi}^{h_0}ds,\;\forall\widetilde{\varphi}^{h_0}\in V_{h_0}(\Omega),\\
&T_{0}^{N}=\widehat{T}(\bm{x},t_N),\;\;\;\text{on}\;\;\partial\Omega_{T}.
\end{aligned} \right.
\end{equation}
\item
If $||\breve T_\lambda(\bm{x})-\breve T_{\lambda-1}(\bm{x})||_{L^\infty(\Omega)}\leq E_{tol}$, stop; otherwise $\lambda=\lambda+1$, go back to Step 2.
\item
Set $T_0^{N+1}=T_{sat}$, where $T_{sat}$ is the solution of linear system (3.45) which satisfies the iteration threshold as $E_{tol}$.
\end{enumerate}
Next, substitute the obtained temperature field $T_0^{N+1}$ into mechanical equation (3.44) and solve equation (3.44) to obtain displacement field $\bm{u}_{0}^{N+1}$.
\item
For arbitrary point $(\bm{x},t)\in \Omega\times(0,T^{*})$, we apply the interpolation method to acquire the corresponding values of first-order auxiliary cell functions, second-order auxiliary cell functions and homogenized solutions. The spatial derivatives $\displaystyle\frac{\partial u_{m0}}{\partial x_{\alpha_1}}$, $\displaystyle\frac{\partial^2 u_{m0}}{\partial x_{\alpha_1}\partial x_{\alpha_2}}$, $\displaystyle\frac{\partial T_{0}}{\partial x_{\alpha_1}}$ and $\displaystyle\frac{\partial^2 T_{0}}{\partial x_{\alpha_1}\partial x_{\alpha_2}}$ in formulas (2.33) and (2.34)
are evaluated by the average technique on relative elements \cite{R45,R46}, and the temporal derivatives $\displaystyle\frac{\partial T_{0}}{\partial t}$ in formula (2.33) are evaluated by using the difference scheme (3.43) at every time steps. Then, the displacement field $\bm{u}^{(2\varepsilon)}(\bm{x},t)$ and temperature field $T^{(2\varepsilon)}(\bm{x},t)$ can be computed by the formulas (2.33) and (2.34). Moreover, we can further employ the higher-order interpolation method and post-processing technique in \cite{R42,R43} to obtain the high-accuracy SOTS solutions.
\end{enumerate}
\section{Numerical experiments}
In this section, the numerical accuracy and computational efficiency of the proposed SHOMS method are demonstrated by several numerical examples. Since obtaining the exact solutions for the multi-scale nonlinear problem (2.1) is extremely difficult or even impossible, direct numerical simulation (DNS) solutions on the fine grid for the multi-scale nonlinear problem (2.1) are taken as the reference solutions denoted as $\bm{u}_{\text{DNS}}(\bm{x},t)$ and $T_{\text{DNS}}(\bm{x},t)$. In the following numerical experiments, some error notations are introduced as follows.
\begin{equation}
\text{{T}err0}=\frac{||T_{0}-T_{\text{DNS}}||_{L^2(\Omega)}}{||T_{\text{DNS}}||_{L^2(\Omega)}},\;
\text{{T}err1}=\frac{||T^{(1\varepsilon)}-T_{\text{DNS}}||_{L^2(\Omega)}}{||T_{\text{DNS}}||_{L^2(\Omega)}},\;
\text{{T}err2}=\frac{||T^{(2\varepsilon)}-T_{\text{DNS}}||_{L^2(\Omega)}}{||T_{\text{DNS}}||_{L^2(\Omega)}}.
\end{equation}
\begin{equation}
\text{{T}Err0}=\frac{|T_{0}-T_{\text{DNS}}|_{H^1(\Omega)}}{|T_{\text{DNS}}|_{H^1(\Omega)}},\;
\text{{T}Err1}=\frac{|T^{(1\varepsilon)}-T_{\text{DNS}}|_{H^1(\Omega)}}{|T_{\text{DNS}}|_{H^1(\Omega)}},\;
\text{{T}Err2}=\frac{|T^{(2\varepsilon)}-T_{\text{DNS}}|_{H^1(\Omega)}}{|T_{\text{DNS}}|_{H^1(\Omega)}}.
\end{equation}
\begin{equation}
\text{Uerr0}=\frac{||\bm{u}_{0}-\bm{u}_{\text{DNS}}||_{L^2(\Omega)}}{||\bm{u}_{\text{DNS}}||_{L^2(\Omega)}},\;
\text{Uerr1}=\frac{||\bm{u}^{(1\varepsilon)}-\bm{u}_{\text{DNS}}||_{L^2(\Omega)}}{||\bm{u}_{\text{DNS}}||_{L^2(\Omega)}},\;
\text{Uerr2}=\frac{||\bm{u}^{(2\varepsilon)}-\bm{u}_{\text{DNS}}||_{L^2(\Omega)}}{||\bm{u}_{\text{DNS}}||_{L^2(\Omega)}}.
\end{equation}
\begin{equation}
\text{UErr0}=\frac{|\bm{u}_{0}-\bm{u}_{\text{DNS}}|_{H^1(\Omega)}}{|\bm{u}_{\text{DNS}}|_{H^1(\Omega)}},\;
\text{UErr1}=\frac{|\bm{u}^{(1\varepsilon)}-\bm{u}_{\text{DNS}}|_{H^1(\Omega)}}{\bm{u}_{\text{DNS}}|_{H^1(\Omega)}},\;
\text{UErr2}=\frac{|\bm{u}^{(2\varepsilon)}-\bm{u}_{\text{DNS}}|_{H^1(\Omega)}}{\bm{u}_{\text{DNS}}|_{H^1(\Omega)}}.
\end{equation}
In the above formulas, $|\bm{u}_{0}-\bm{u}_{\text{DNS}}|_{H^1(\Omega)}=\Big(\sum\limits_{i,j=1}^{\mathcal{N}}  \big|\big|\varepsilon_{ij}(\bm{u}_{0}-\bm{u}_{\text{DNS}})\big|\big|_{L^2(\Omega)}\Big)^{\frac{1}{2}}$, where the symbol $\varepsilon_{ij}$ represents the strain operator.
\subsection{Example 1: Validation of the SHOMS method for nonlinear thermo-mechanical simulation in periodic composite structure}
A 2D composite structure with periodically microscopic configurations is investigated here, whose macrostructure $\Omega$ and unit cell $\Theta$ are shown in Fig.\hspace{1mm}3, where $\Omega=(x_1,x_2)=[0,1]\times[0,1]cm^2$ and periodic unit cell size $\varepsilon=1/5$. This 2D composite structure is clamped on its four side surfaces, and the temperature at the side surfaces is kept at $373.15K$. In addition, we set internal heat source $h = 5000J/(cm^2\bm{\cdot}s)$ and body forces $(f_1,f_2) = (-2000,-2000)N/cm^2$. Furthermore, the material parameters of 2D composite structure are demonstrated in Table 1.
\begin{figure}[!htb]
\centering
\begin{minipage}[c]{0.32\textwidth}
  \centering
  \includegraphics[width=0.8\linewidth,totalheight=1.6in]{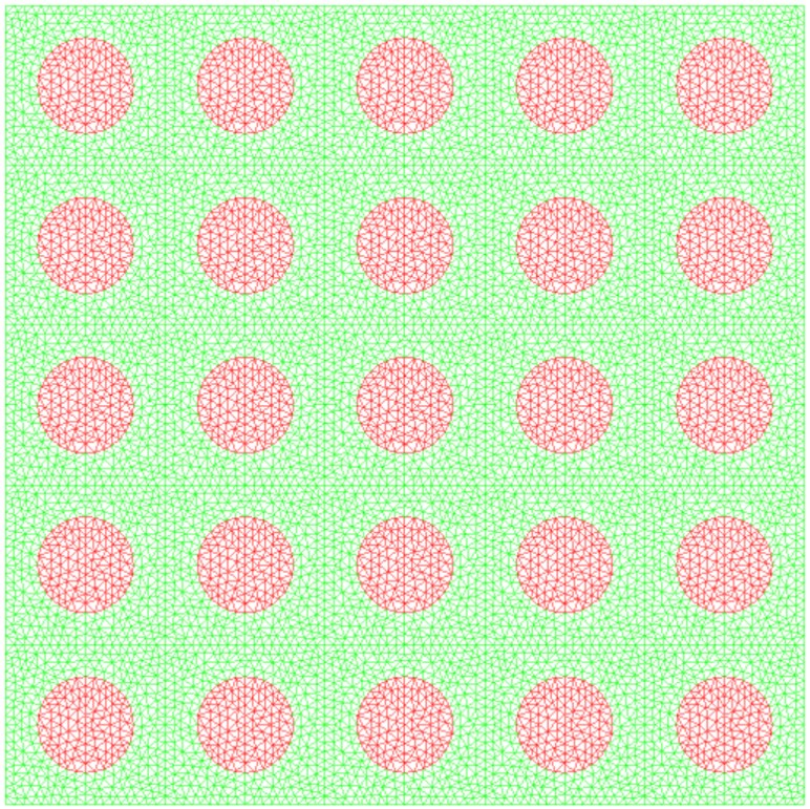} \\
  (a)
\end{minipage}
\begin{minipage}[c]{0.32\textwidth}
  \centering
  \includegraphics[width=0.8\linewidth,totalheight=1.6in]{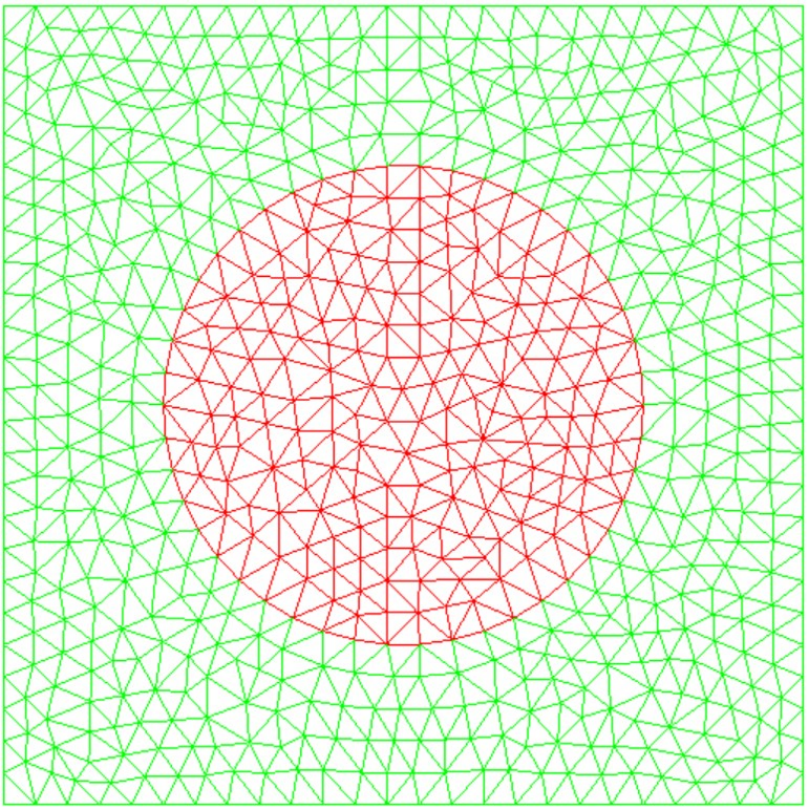} \\
  (b)
\end{minipage}
\begin{minipage}[c]{0.32\textwidth}
  \centering
  \includegraphics[width=0.8\linewidth,totalheight=1.6in]{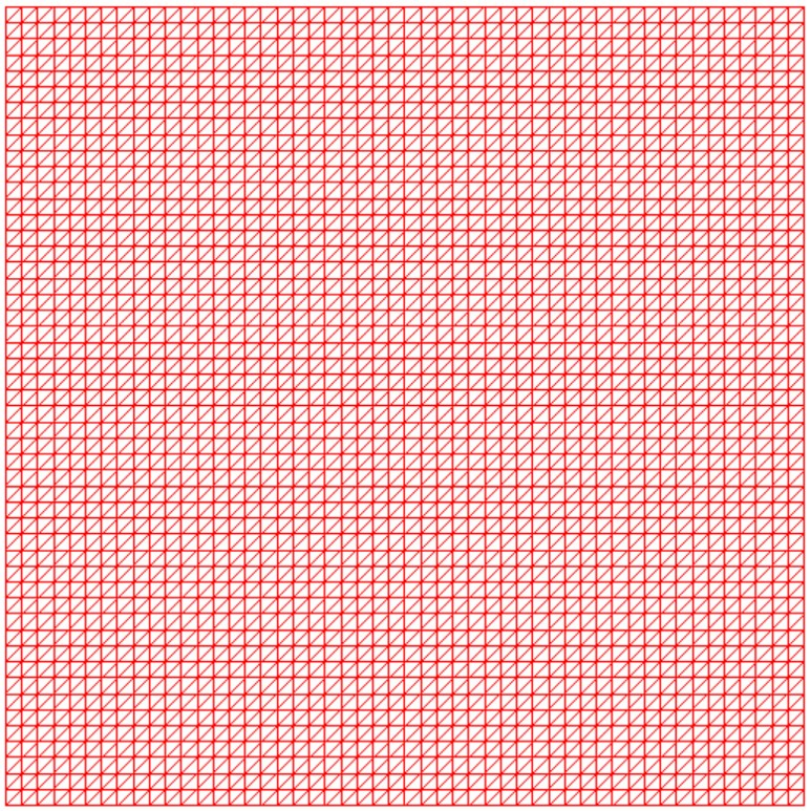} \\
  (a)
\end{minipage}
\caption{The illustration of investigated 2D composite structure: (a) FEM mesh of composite structure; (b) FEM mesh of microscopic unit cell; (c) FEM mesh of macroscopic homogenized structure.}
\end{figure}

\begin{table}[h]{\caption{Material parameters of 2D composite structure.}}
\centering
\begin{tabular}{ccc}
\hline
Material parameters & Matrix & Inclusion\\
\hline
$\rho^\varepsilon$  $(kg/m^3)$ & 3210&1760   \\
$c^\varepsilon$ $(J/(kg\bm{\cdot} K))$ & 660.0+1.915$T$-1.491$\times10^{-3}T^2$&710.0+$1.781T$-1.976$\times10^{-3}T^2$  \\
$k_{ij}^\varepsilon$ $(W/(m\bm{\cdot} K))$ & 250.0+0.02728$T$&8.0+0.02535$T$  \\
$E^\varepsilon$ $(GPa)$ & 350.0-3.04$\times10^{-2}T$&220.0-1.10$\times10^{-4}T$ \\
$\nu^\varepsilon$ & 0.25&0.20   \\
$\beta_{ij}^\varepsilon$ $(Pa/K)$ & 3.50&9273.0-57.53$\times T$+0.0817$\times T^2$  \\
\hline
\end{tabular}
\end{table}

Implementing the SHOMS method to multi-scale nonlinear coupling equations (2.1) in time interval $t=[0,1]s$ with temporal step $\Delta t=0.002s$, we define service temperature range $[273.15,873.15]K$ of investigated composite structure and 60 representative macroscopic parameters $\bar T_{s}$ in service temperature range. In this example, the total auxiliary cell problems need to be solved off-line 4380 times, in which the 13 first-order cell functions and 60 second-order cell functions are solved on 60 macroscopic temperature interpolation points. The information of FEM meshes is listed in Table 2. After off-line computation for microscopic cell problems and on-line computation for macroscopic homogenized problems and higher-order multi-scale solutions, we depict the computational results in Figs.\hspace{1mm}4-8.
\begin{table}[!htb]{\caption{Summary of computational cost.}}
\centering
\begin{tabular}{cccc}
\hline
 & Multi-scale equations. & Cell equations. & Homogenized equations. \\
\hline
FEM elements & 17700 & 1370 & 5000\\
FEM nodes    & 9031 & 736 & 2601\\
\hline
\end{tabular}
\end{table}

\begin{figure}[!htb]
\centering
\begin{minipage}[c]{0.4\textwidth}
  \centering
  \includegraphics[width=50mm]{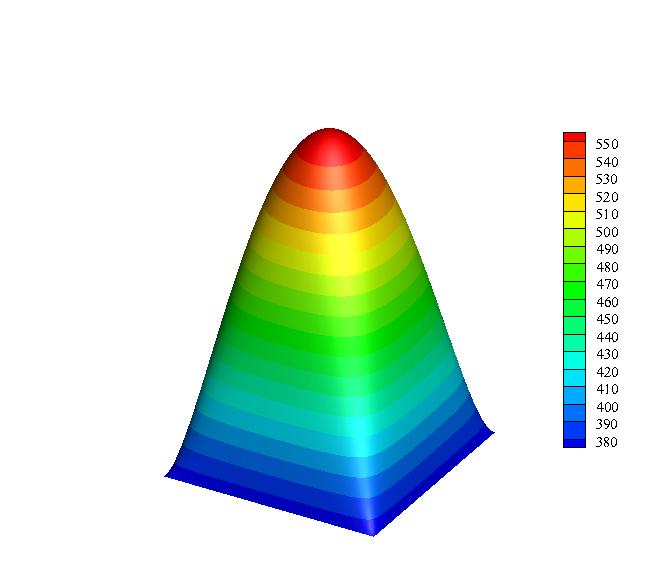}\\
  (a)
\end{minipage}
\begin{minipage}[c]{0.4\textwidth}
  \centering
  \includegraphics[width=50mm]{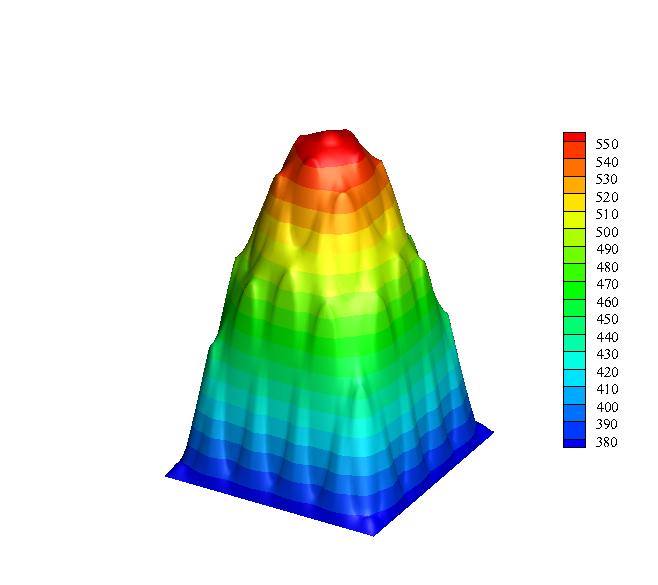}\\
  (b)
\end{minipage}
\begin{minipage}[c]{0.4\textwidth}
  \centering
  \includegraphics[width=50mm]{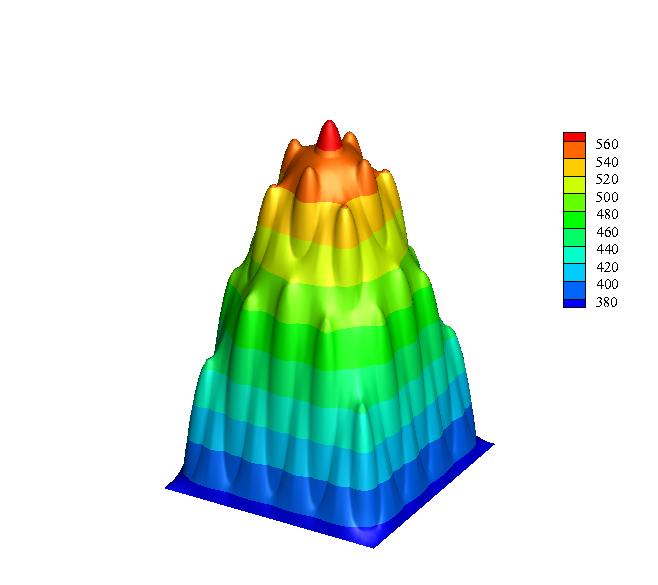}\\
  (c)
\end{minipage}
\begin{minipage}[c]{0.4\textwidth}
  \centering
  \includegraphics[width=50mm]{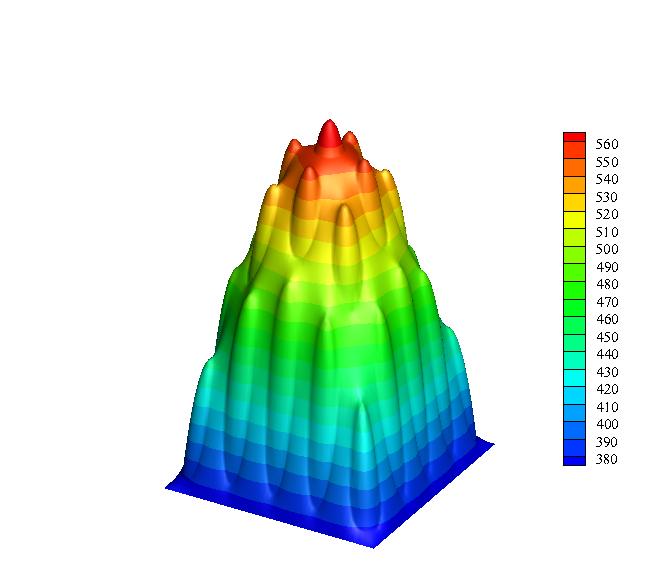}\\
  (d)
\end{minipage}
\caption{The numerical results of temperature field at $t=0.2s$: (a) $T_0$; (b) $T^{(1\varepsilon)}$; (c) $T^{(2\varepsilon)}$; (d) $T_{\text{DNS}}$.}
\end{figure}

\begin{figure}[!htb]
\centering
\begin{minipage}[c]{0.4\textwidth}
  \centering
  \includegraphics[width=50mm]{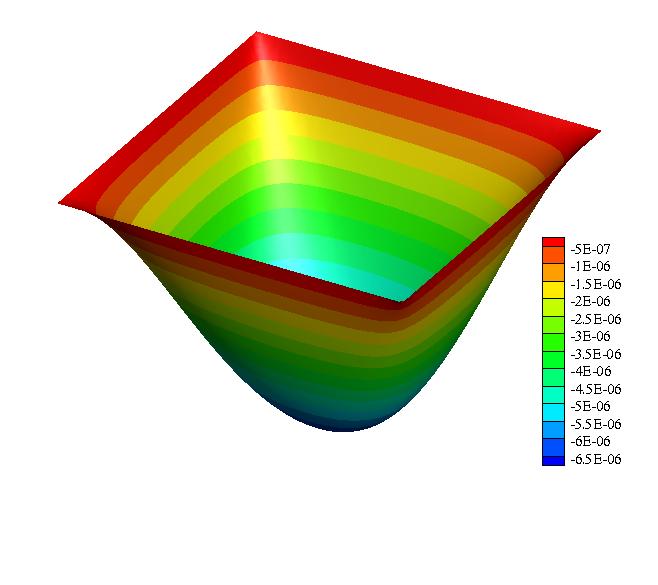}\\
  (a)
\end{minipage}
\begin{minipage}[c]{0.4\textwidth}
  \centering
  \includegraphics[width=50mm]{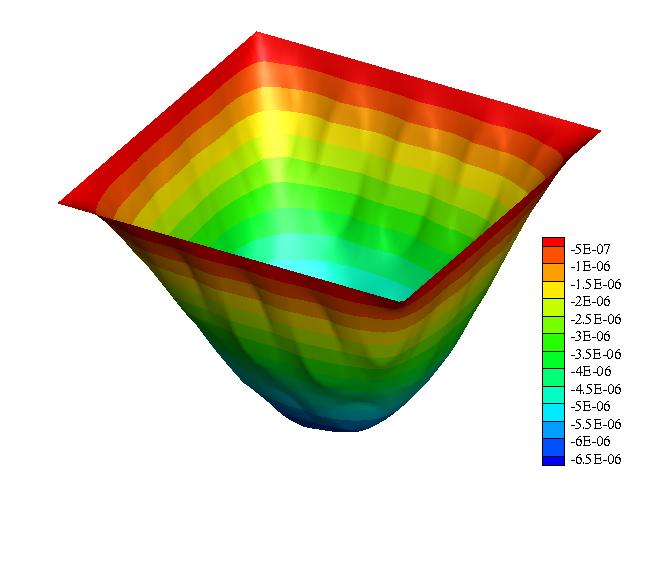}\\
  (b)
\end{minipage}
\begin{minipage}[c]{0.4\textwidth}
  \centering
  \includegraphics[width=50mm]{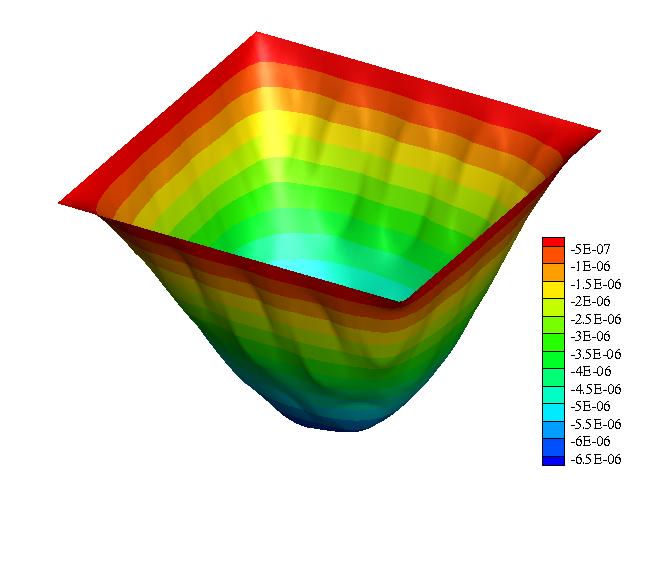}\\
  (c)
\end{minipage}
\begin{minipage}[c]{0.4\textwidth}
  \centering
  \includegraphics[width=50mm]{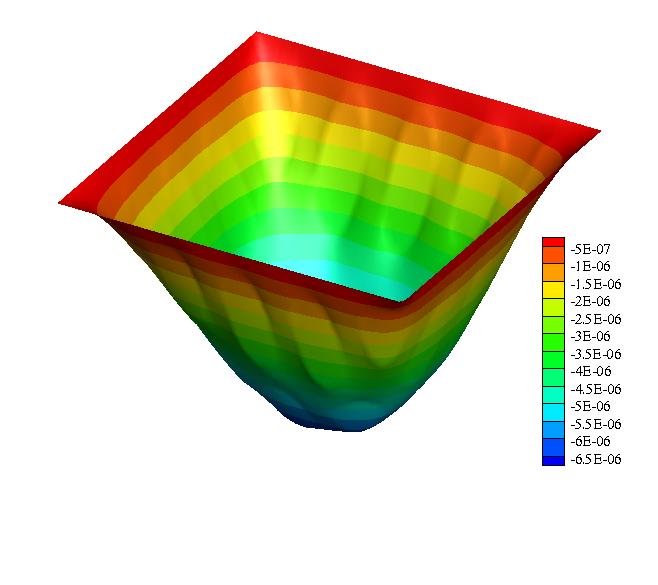}\\
  (d)
\end{minipage}
\caption{The numerical results of temperature field at $t=0.2s$: (a) $u_{10}$; (b) $u_1^{(1\varepsilon)}$; (c) $u_1^{(2\varepsilon)}$; (d) $u_{1\text{DNS}}$.}
\end{figure}

\begin{figure}[!htb]
\centering
\begin{minipage}[c]{0.4\textwidth}
  \centering
  \includegraphics[width=50mm]{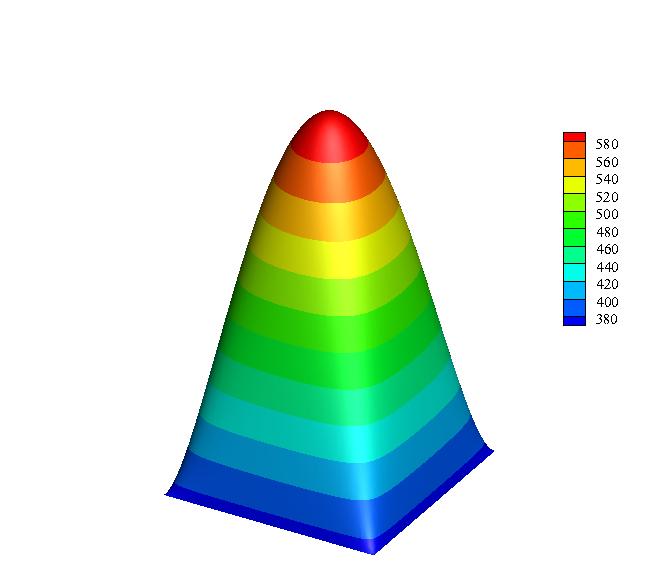}\\
  (a)
\end{minipage}
\begin{minipage}[c]{0.4\textwidth}
  \centering
  \includegraphics[width=50mm]{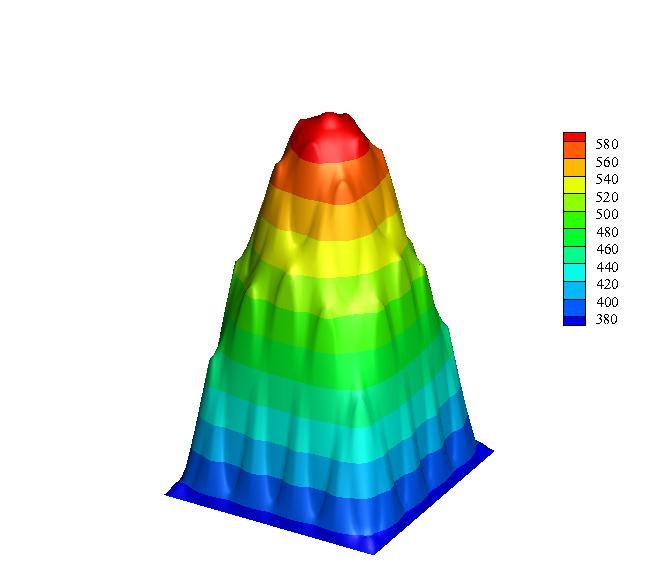}\\
  (b)
\end{minipage}
\begin{minipage}[c]{0.4\textwidth}
  \centering
  \includegraphics[width=50mm]{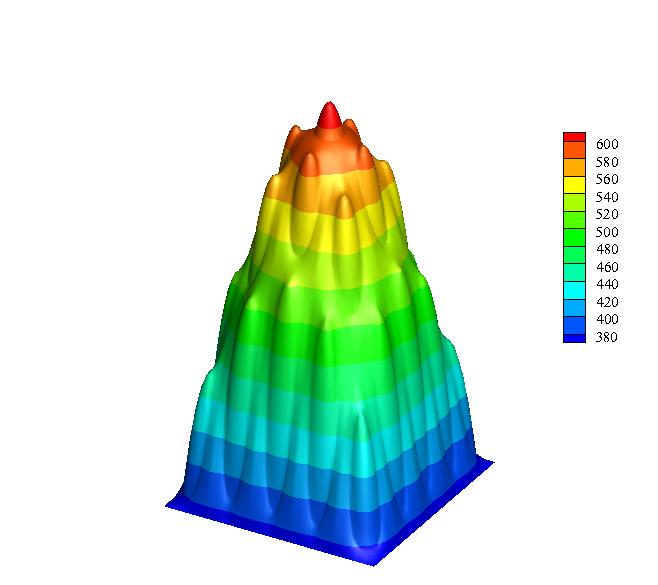}\\
  (c)
\end{minipage}
\begin{minipage}[c]{0.4\textwidth}
  \centering
  \includegraphics[width=50mm]{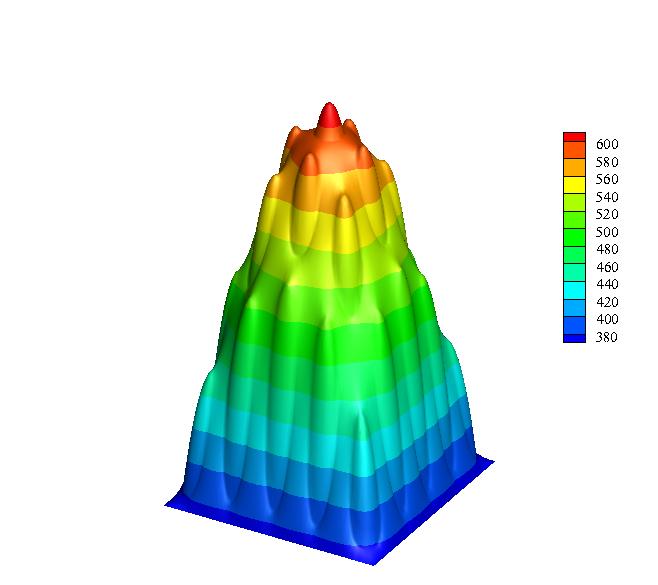}\\
  (d)
\end{minipage}
\caption{The numerical results of temperature field at $t=1.0s$: (a) $T_0$; (b) $T^{(1\varepsilon)}$; (c) $T^{(2\varepsilon)}$; (d) $T_{\text{DNS}}$.}
\end{figure}

\begin{figure}[!htb]
\centering
\begin{minipage}[c]{0.4\textwidth}
  \centering
  \includegraphics[width=50mm]{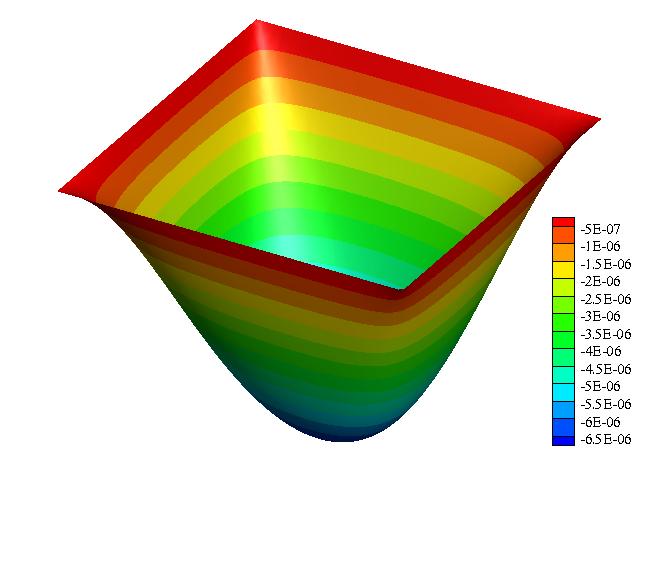}\\
  (a)
\end{minipage}
\begin{minipage}[c]{0.4\textwidth}
  \centering
  \includegraphics[width=50mm]{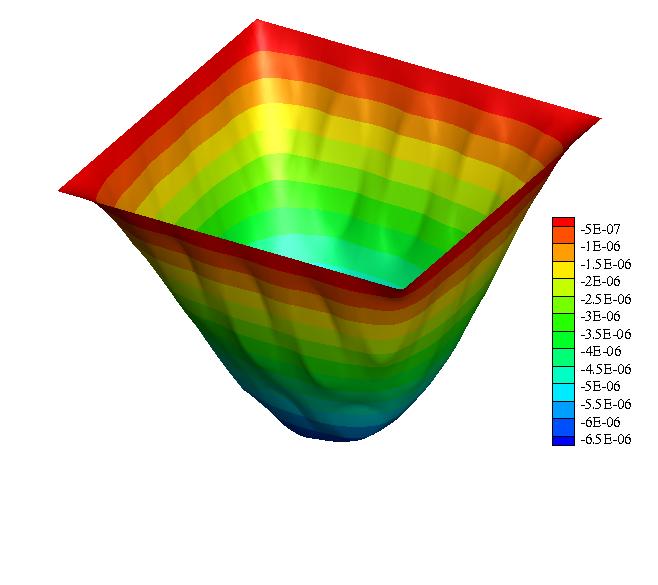}\\
  (b)
\end{minipage}
\begin{minipage}[c]{0.4\textwidth}
  \centering
  \includegraphics[width=50mm]{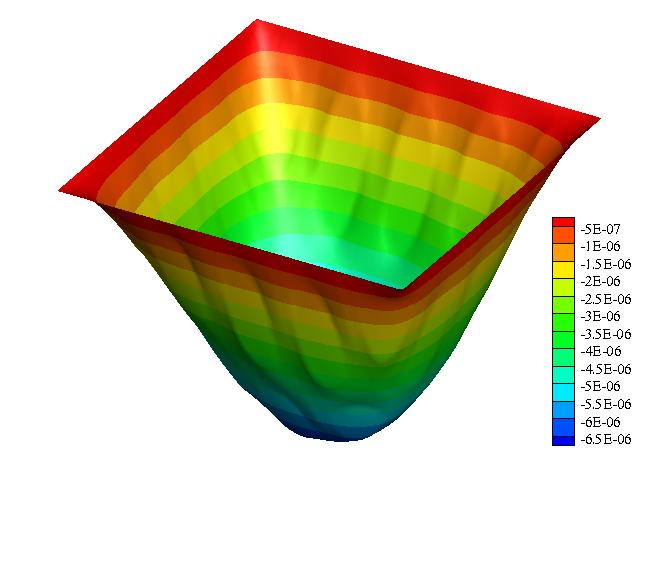}\\
  (c)
\end{minipage}
\begin{minipage}[c]{0.4\textwidth}
  \centering
  \includegraphics[width=50mm]{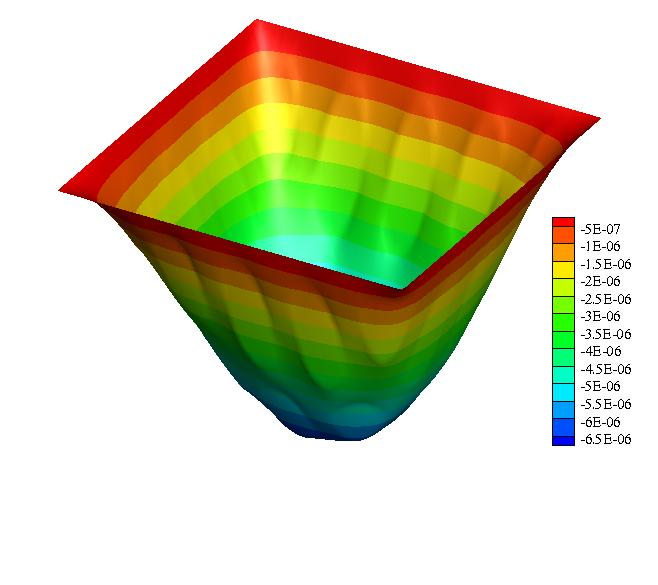}\\
  (d)
\end{minipage}
\caption{The numerical results of temperature field at $t=1.0s$: (a) $u_{10}$; (b) $u_1^{(1\varepsilon)}$; (c) $u_1^{(2\varepsilon)}$; (d) $u_{1\text{DNS}}$.}
\end{figure}

\begin{figure}[!htb]
\centering
\begin{minipage}[c]{0.4\textwidth}
  \centering
  \includegraphics[width=50mm]{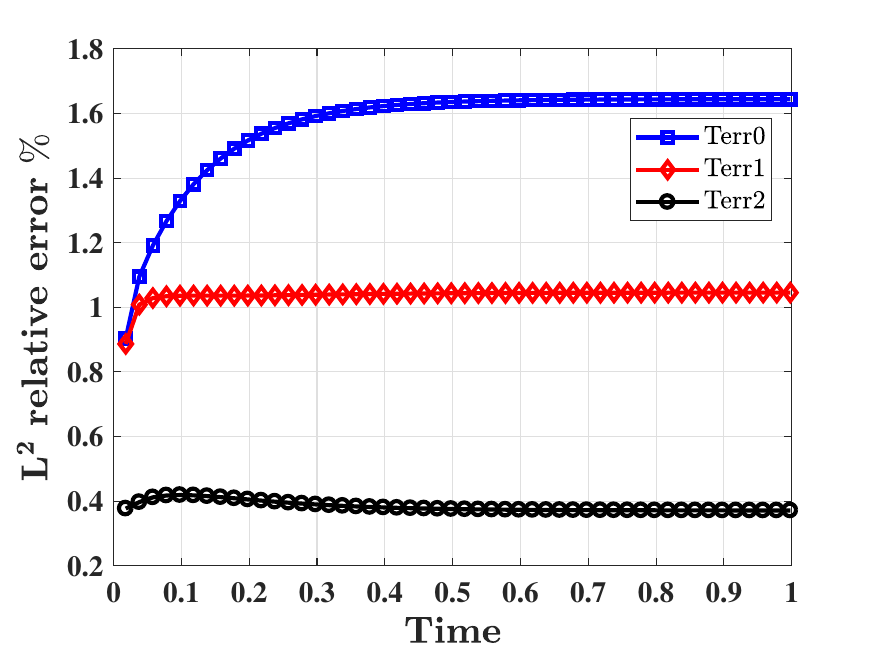}\\
  (a)
\end{minipage}
\begin{minipage}[c]{0.4\textwidth}
  \centering
  \includegraphics[width=50mm]{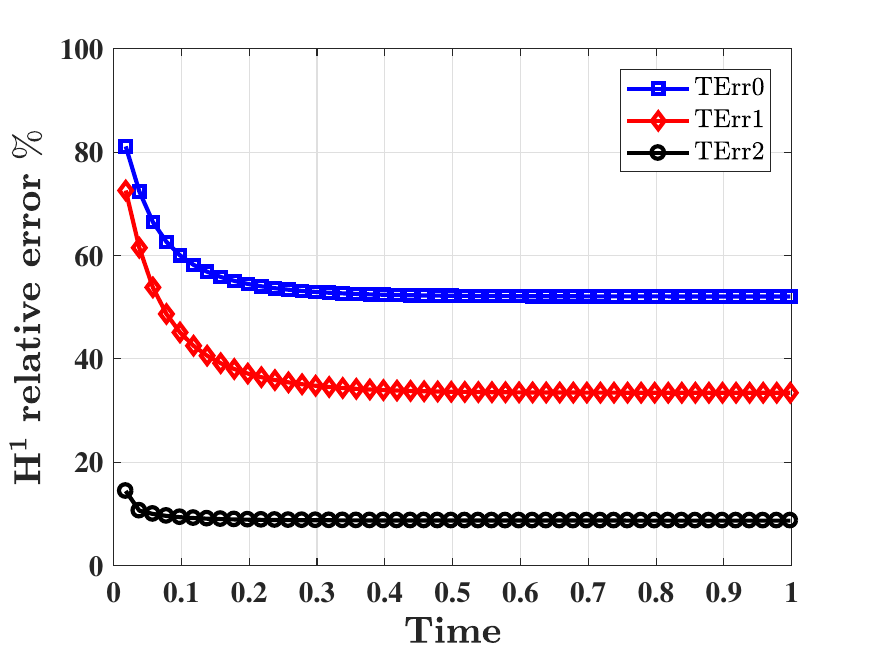}\\
  (b)
\end{minipage}
\begin{minipage}[c]{0.4\textwidth}
  \centering
  \includegraphics[width=50mm]{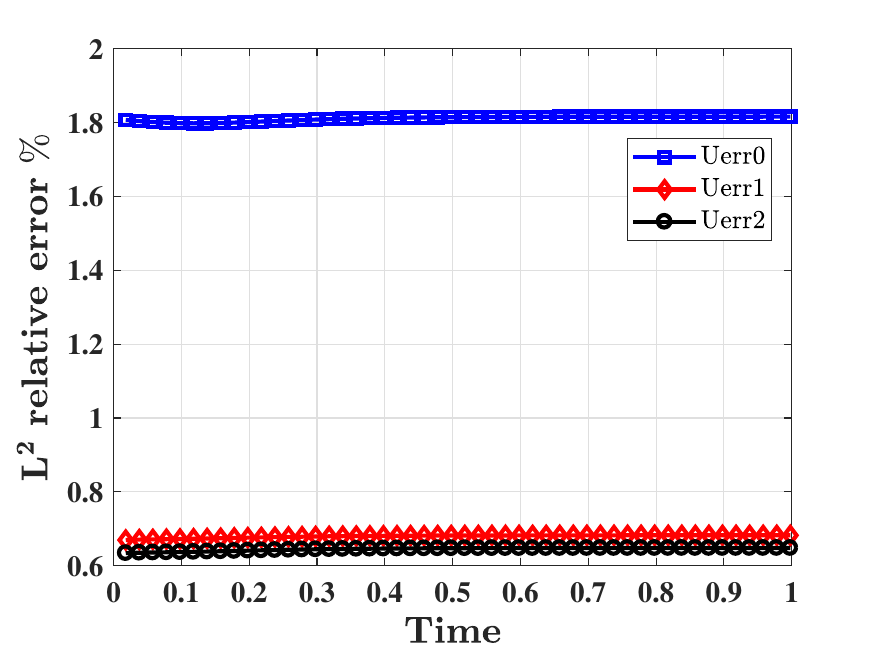}\\
  (c)
\end{minipage}
\begin{minipage}[c]{0.4\textwidth}
  \centering
  \includegraphics[width=50mm]{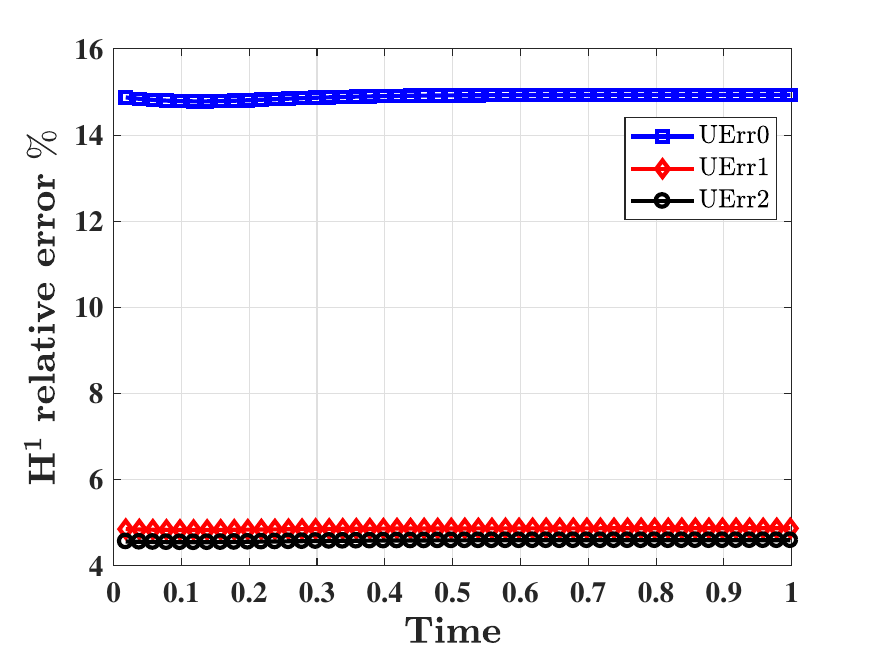}\\
  (d)
\end{minipage}
\caption{The evolutive relative errors of temperature and displacement fields: (a) \rm{Terr}; (b) \rm{TErr}; (c) \rm{Uerr}; (d) \rm{UErr}.}
\end{figure}

According to the computational resource cost in Table 2, the presented SHOMS method can greatly economize computer memory without losing precision. Actually, both the SHOMS method and direct numerical simulation are performed on a HP desktop workstation equipped with an Intel(R) Core(TM) i7-8750H processor (2.20 GHz) and 16.0 GB of internal memory. As illustrated in Figs.\hspace{1mm}4-7, we can conclude that the higher-order multi-scale solutions can accurately capture the microscopic oscillatory behaviors and preferably approximate the exact solutions of the investigated 2D composite structure compared with macroscopic homogenized solutions and lower-order multi-scale solutions, especially for temperature field. From the evolutive relative errors in Fig.\hspace{1mm}8, it can clearly demonstrate that the two-stages space-time multi-scale numerical algorithm is accurate and stable in the long-time numerical simulation. Furthermore, it is worth emphasizing that
the presented SHOMS approach remains effective even for a relatively small parameter $\varepsilon$, namely existing a great number of microscopic unit cells in inhomogeneous structures. At this time, the high-resolution DNS simulation can not guarantee the convergence for the investigated large-scale problems. This obvious advantage of the SHOMS approach is of great application values for engineering computation.
\subsection{Example 2: Application of the SHOMS method for equivalent material parameters computation of random composite structure}
In this example, two kinds of composite materials with matrix Ti-6Al-4V and random inclusion ZrO$_2$, and matrix SiC and random inclusion C are investigated by the SHOMS method, as exhibited in Figs.\hspace{1mm}9. The detailed material parameters for the investigated composite materials are presented in the following Tables 3 and 4.
\begin{table}[h]{\caption{Material property parameters of Ti-6Al-4V/ZrO$_2$ composite.}}
\centering
\begin{tabular}{ccc}
\hline
Material properties & Matrix Ti-6Al-4V & Inclusion ZrO$_2$ \\
\hline
$k_{ij}^\varepsilon$ $(W/(m\bm{\cdot} K))$ & 1.10+0.017$T$ & 1.71+2.1$\times10^{-4}T$+1.16$\times10^{-7}T^2$ \\
$E^\varepsilon$ $(GPa)$ & 122.7-5.65$\times10^{-3}T$ & 132.2-50.3$\times10^{-3}T$-8.1$\times10^{-6}T^2$ \\
$\nu^\varepsilon$ & 0.289 & 0.333  \\
\hline
\end{tabular}
\end{table}

\begin{table}[h]{\caption{Material property parameters of SiC/C composite.}}
\centering
\begin{tabular}{ccc}
\hline
Material properties & Matrix SiC & Inclusion C \\
\hline
$k_{ij}^\varepsilon$ $(W/(m\bm{\cdot} K))$ & 250.0+0.02728$T$ & 8.0+0.02535$T$  \\
$E^\varepsilon$ $(GPa)$ & 350.0-3.04$\times10^{-2}T$ & 220.0-1.10$\times10^{-4}T$ \\
$\nu^\varepsilon$ & 0.25 & 0.20  \\
\hline
\end{tabular}
\end{table}

\begin{figure}[!htb]
\centering
\begin{minipage}[c]{0.32\textwidth}
  \centering
  \includegraphics[width=35mm]{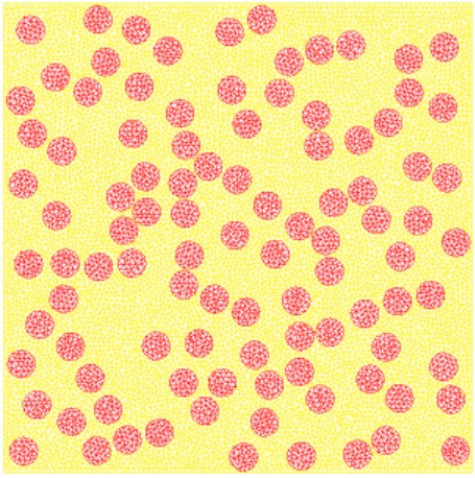}\\
  (a)
\end{minipage}
\begin{minipage}[c]{0.32\textwidth}
  \centering
  \includegraphics[width=35mm]{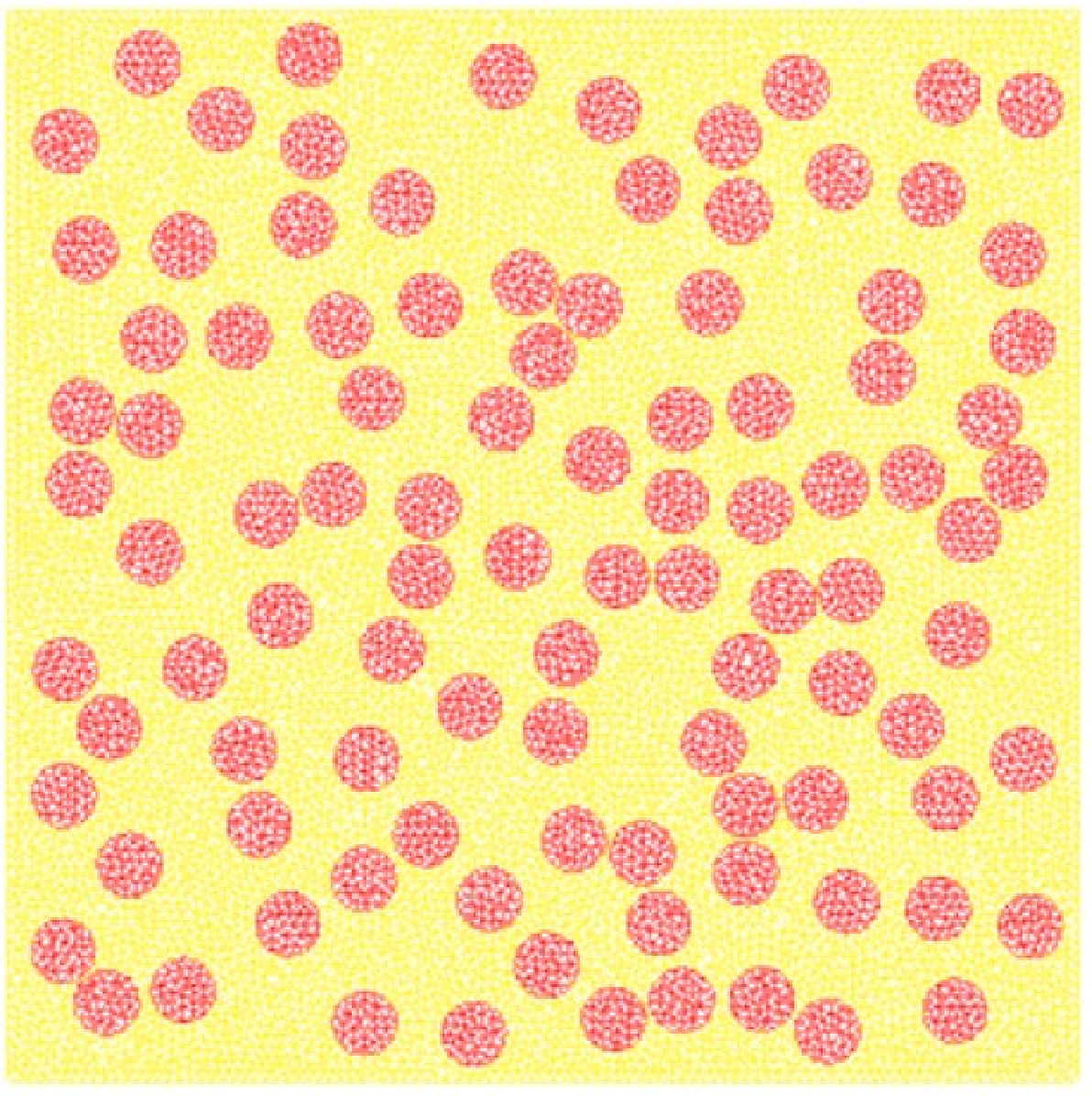}\\
  (b)
\end{minipage}
\begin{minipage}[c]{0.32\textwidth}
  \centering
  \includegraphics[width=35mm]{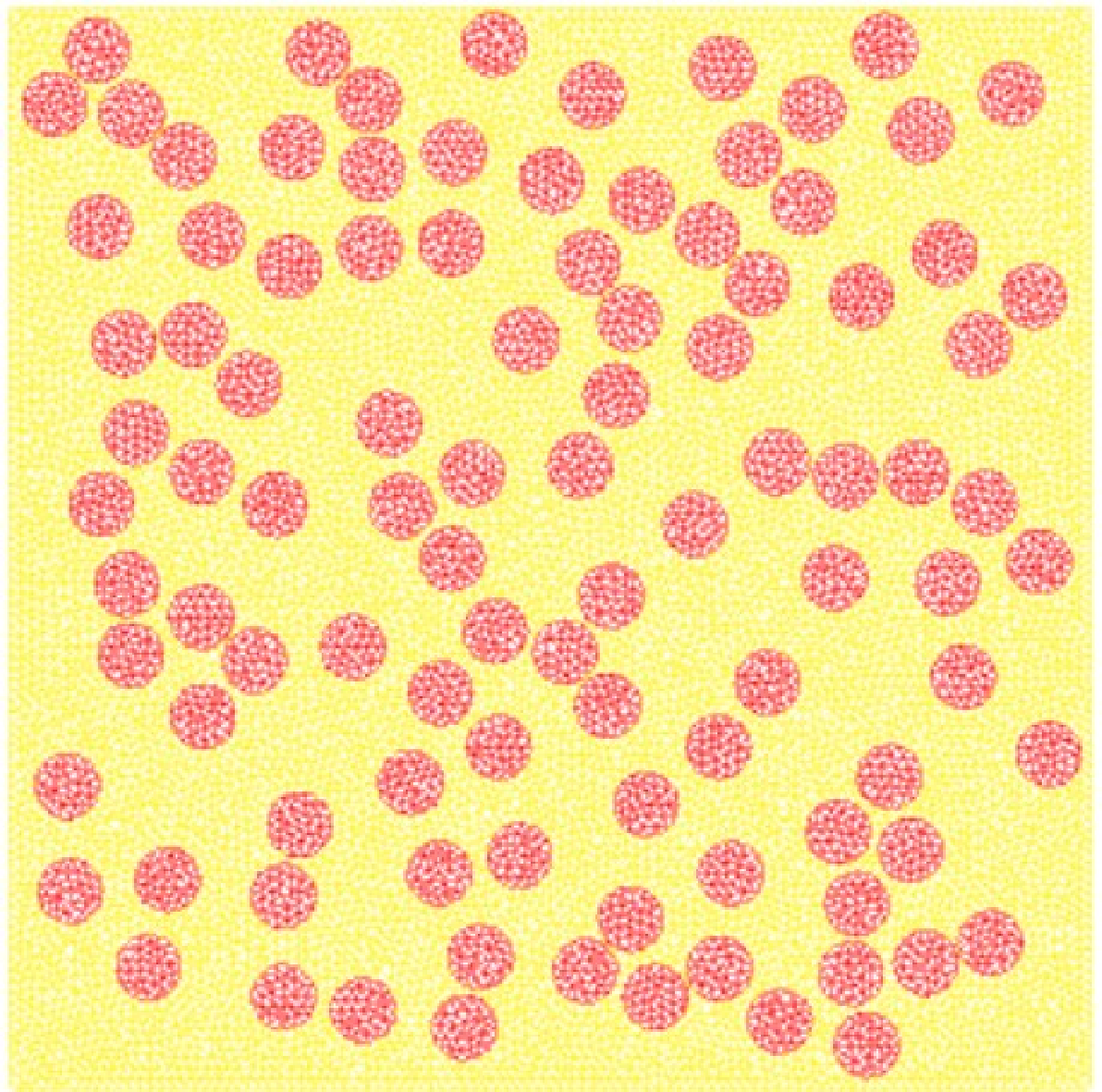}\\
  (c)
\end{minipage}
\begin{minipage}[c]{0.32\textwidth}
  \centering
  \includegraphics[width=35mm]{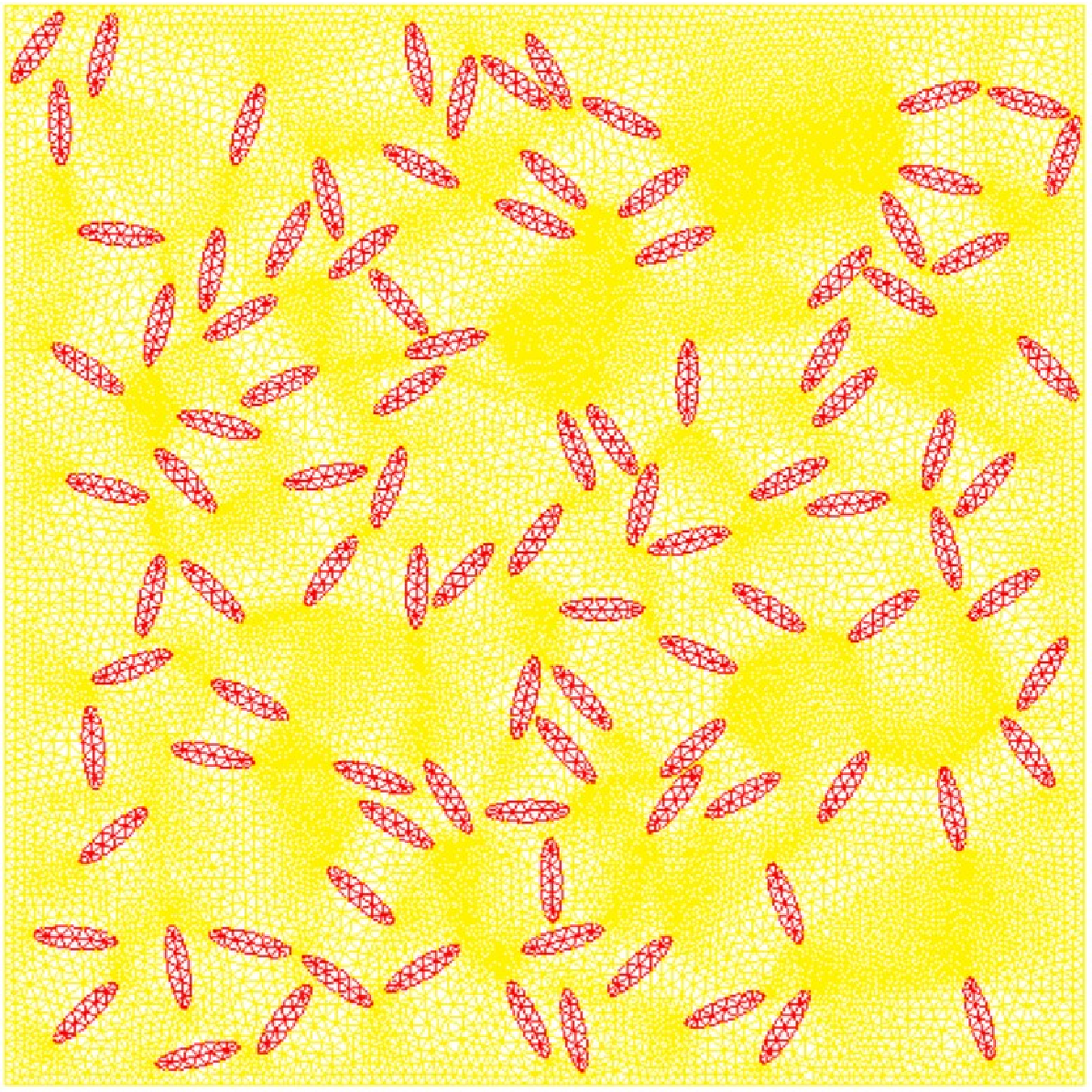}\\
  (d)
\end{minipage}
\begin{minipage}[c]{0.32\textwidth}
  \centering
  \includegraphics[width=35mm]{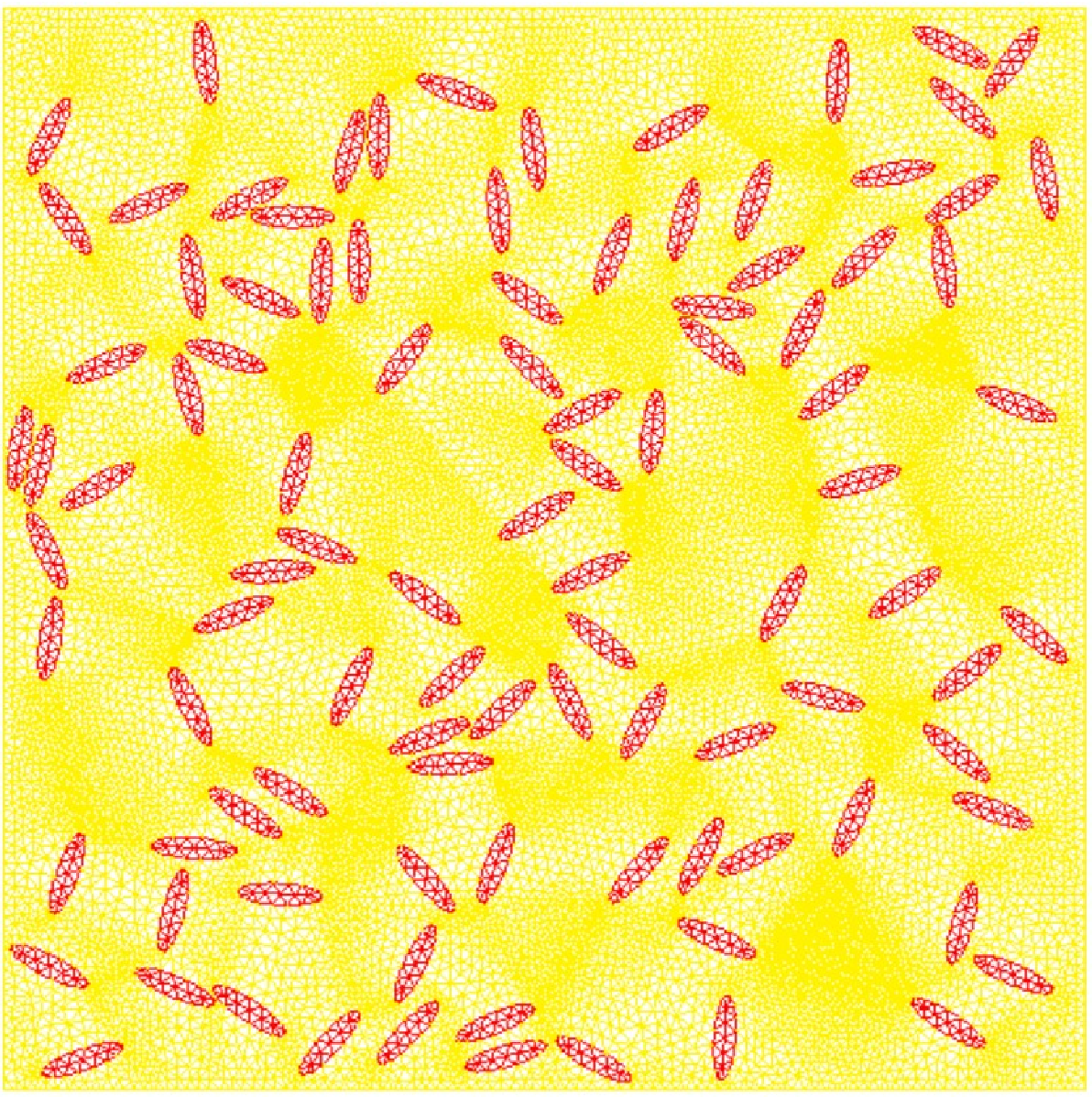}\\
  (e)
\end{minipage}
\begin{minipage}[c]{0.32\textwidth}
  \centering
  \includegraphics[width=35mm]{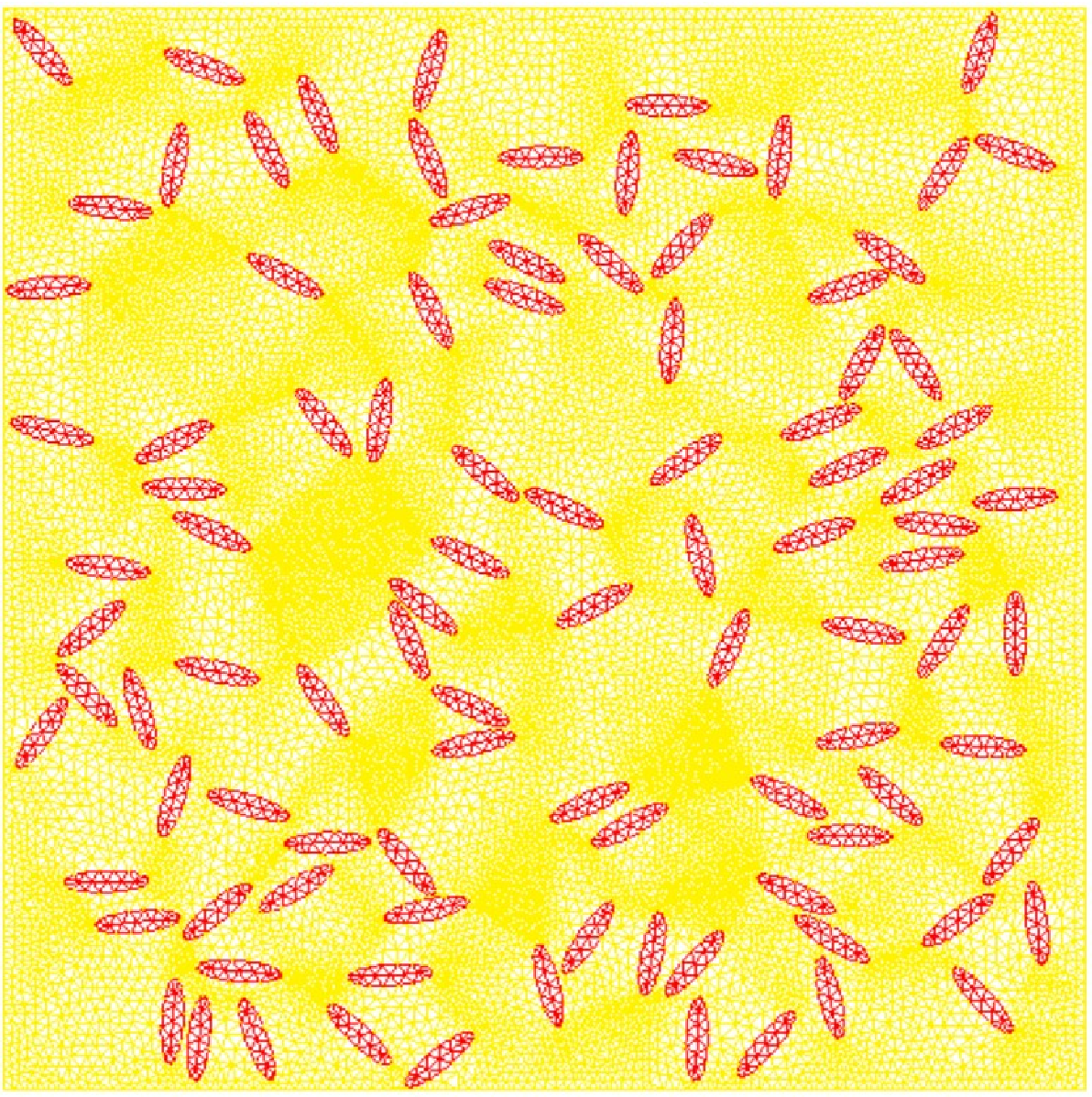}\\
  (f)
\end{minipage}
\caption{Several random RVEs of the composites employed for predicting equivalent material parameters: (a)-(c) particulate composites with volume fraction 28.2743\%; (d)-(f) fibrous composites with volume fraction 12.5664\%.}
\end{figure}

\begin{figure}[!htb]
\centering

\begin{minipage}[c]{0.35\textwidth}
  \centering
  \includegraphics[width=55mm]{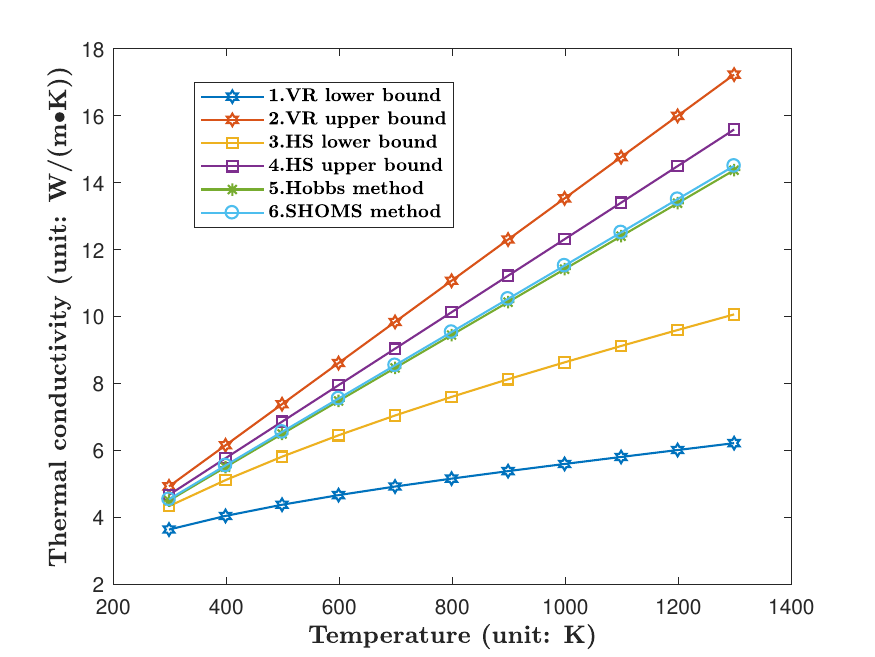}\\
  (a)
\end{minipage}
\begin{minipage}[c]{0.35\textwidth}
  \centering
  \includegraphics[width=55mm]{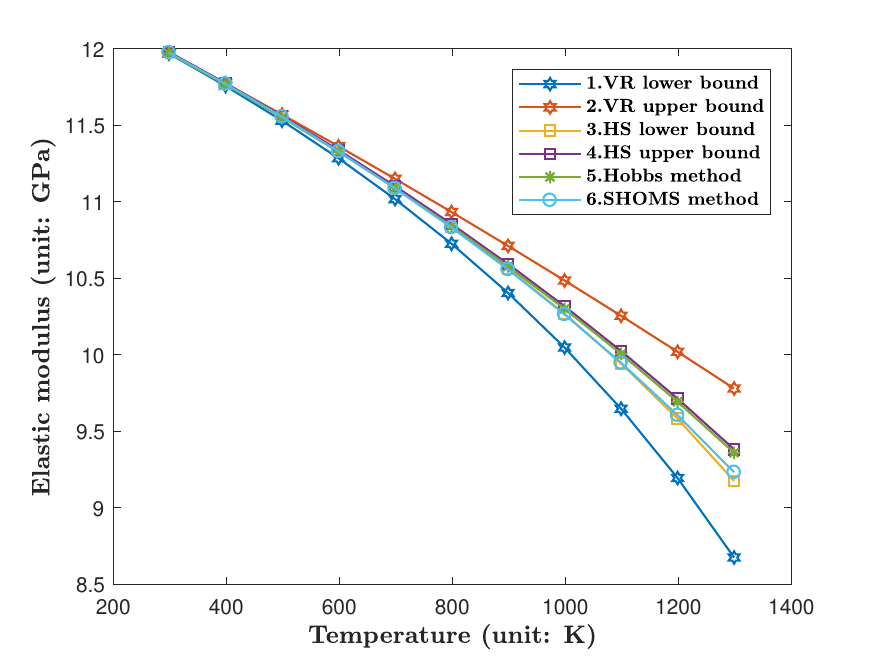}\\
  (b)
\end{minipage}
\begin{minipage}[c]{0.35\textwidth}
  \centering
  \includegraphics[width=55mm]{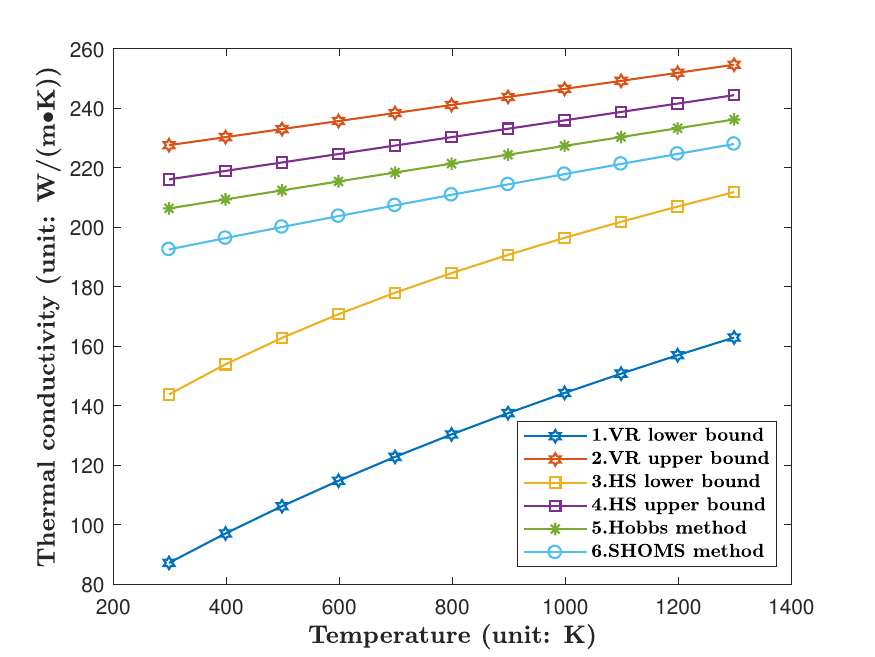}\\
  (c)
\end{minipage}
\begin{minipage}[c]{0.35\textwidth}
  \centering
  \includegraphics[width=55mm]{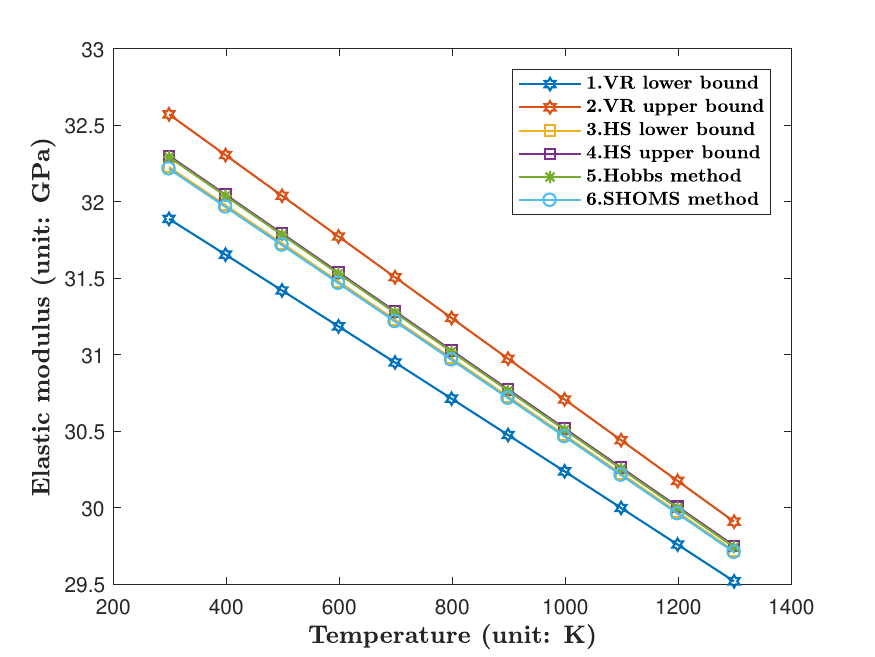}\\
  (d)
\end{minipage}
\caption{A comparison of the predictive results of equivalent material parameters: (a) nonlinear elastic modulus of particulate composites; (b) nonlinear thermal conductivity of particulate composites; (c) nonlinear elastic modulus of fibrous composites; (d) nonlinear thermal conductivity of fibrous composites.}
\end{figure}

By using the SHOMS method, the equivalent material parameters at macro-scale are obtained by the mean value of 50 randomly microscopic samples. The corresponding results are depicted in Fig.\hspace{1mm}10. According to the numerical results in Fig.\hspace{1mm}10, we can conclude that the predictive values of Ti-6Al-4V/ZrO$_2$ composite and SiC/C composite fall between lower and upper bounds of Voigt-Reuss method, Hashin¨CShtrikman method and also approximate the predicted values of Hobbs method. Hence, the proposed SHOMS can be employed to predict the temperature-dependent equivalent material properties of Ti-6Al-4V/ZrO$_2$ composite and SiC/C composite.
\subsection{Example 3: Application of the SHOMS method for nonlinear thermo-mechanical simulation in random composite structure}
This example study the nonlinear thermo-mechanical simulation of 2D composite structure with randomly microscopic configurations, as depicted in Fig.\hspace{1mm}11. In addition, the setting of initial-boundary conditions, heat source, body forces and material parameters in this example is the same as those of Example 1.
\begin{figure}[!htb]
\centering
\begin{minipage}[c]{0.32\textwidth}
  \centering
  \includegraphics[width=0.8\linewidth,totalheight=1.6in]{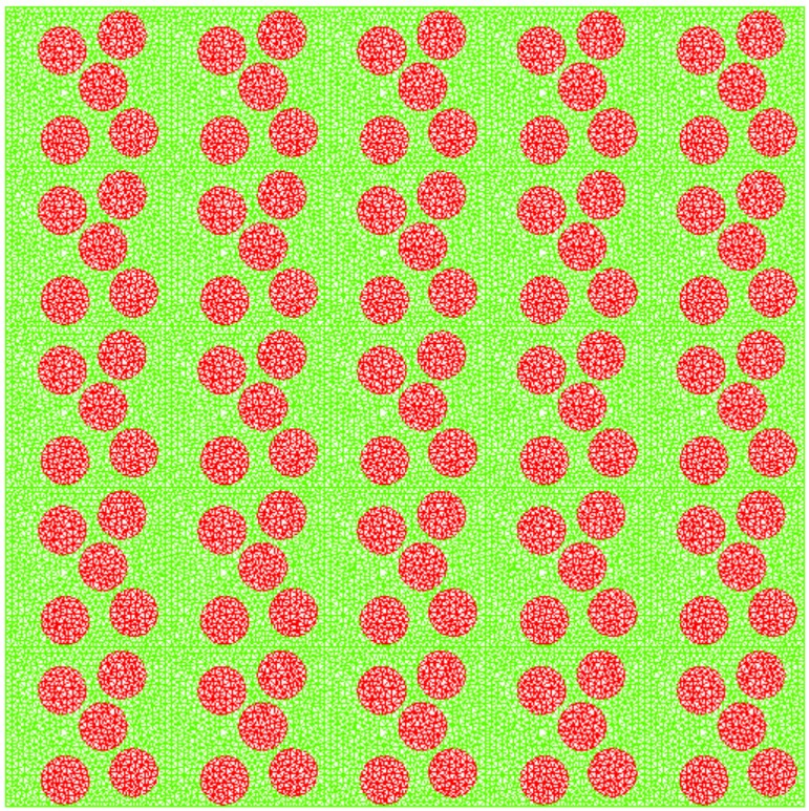} \\
  (a)
\end{minipage}
\begin{minipage}[c]{0.32\textwidth}
  \centering
  \includegraphics[width=0.8\linewidth,totalheight=1.6in]{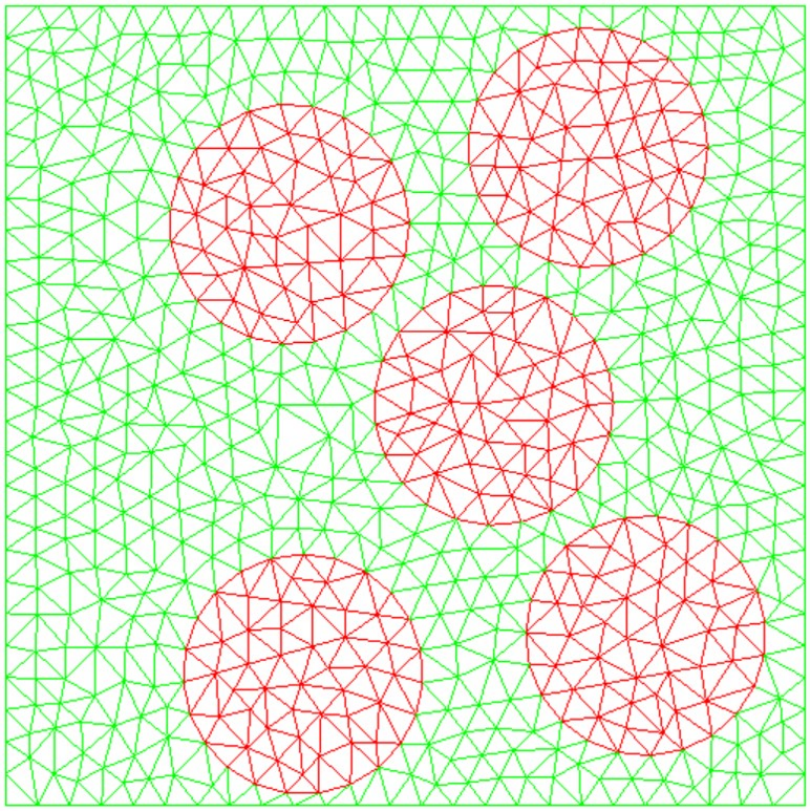} \\
  (b)
\end{minipage}
\begin{minipage}[c]{0.32\textwidth}
  \centering
  \includegraphics[width=0.8\linewidth,totalheight=1.6in]{homodomain1.pdf} \\
  (a)
\end{minipage}
\caption{The illustration of investigated 2D composite structure: (a) FEM mesh of composite structure; (b) FEM mesh of microscopic unit cell; (c) FEM mesh of macroscopic homogenized structure.}
\end{figure}

Applying the SHOMS method to multi-scale nonlinear coupling equations (2.1) within time interval $t=[0,1]s$ with temporal step $\Delta t=0.002s$, we establish the service temperature range of the investigated composite structure as $[273.15,873.15]K$. In this service temperature range, we distribute 60 representative macroscopic parameters $\bar T_{s}$. The detailed information of FEM meshes is listed in Table 5. After off-line computation for microscopic cell problems and on-line computation for macroscopic homogenized problems and higher-order multi-scale solutions, we present the computational results in Figs.\hspace{1mm}12-16.
\begin{table}[!htb]{\caption{Summary of computational cost.}}
\centering
\begin{tabular}{cccc}
\hline
 & Multi-scale equations. & Cell equations. & Homogenized equations. \\
\hline
FEM elements & 33800 & 1352 & 5000\\
FEM nodes    & 17151 & 727 & 2601\\
\hline
\end{tabular}
\end{table}

\begin{figure}[!htb]
\centering
\begin{minipage}[c]{0.4\textwidth}
  \centering
  \includegraphics[width=50mm]{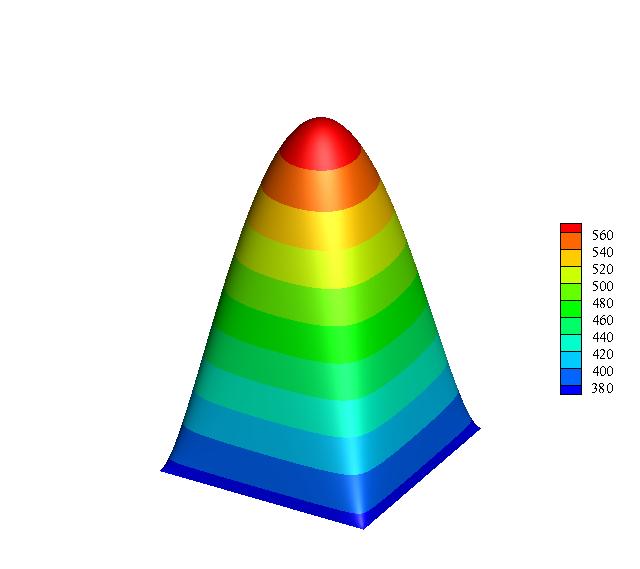}\\
  (a)
\end{minipage}
\begin{minipage}[c]{0.4\textwidth}
  \centering
  \includegraphics[width=50mm]{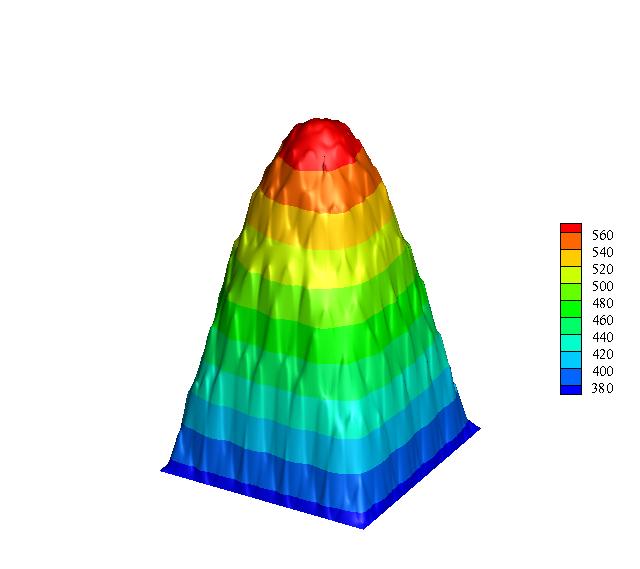}\\
  (b)
\end{minipage}
\begin{minipage}[c]{0.4\textwidth}
  \centering
  \includegraphics[width=50mm]{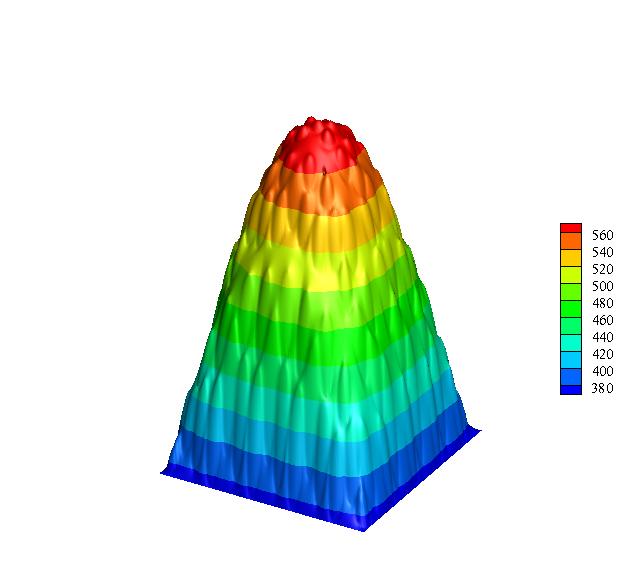}\\
  (c)
\end{minipage}
\begin{minipage}[c]{0.4\textwidth}
  \centering
  \includegraphics[width=50mm]{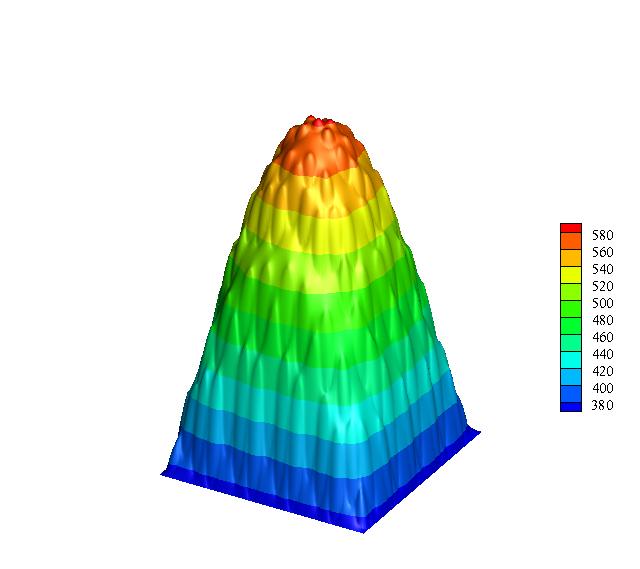}\\
  (d)
\end{minipage}
\caption{The numerical results of temperature field at $t=0.2s$: (a) $T_0$; (b) $T^{(1\varepsilon)}$; (c) $T^{(2\varepsilon)}$; (d) $T_{\text{DNS}}$.}
\end{figure}

\begin{figure}[!htb]
\centering
\begin{minipage}[c]{0.4\textwidth}
  \centering
  \includegraphics[width=50mm]{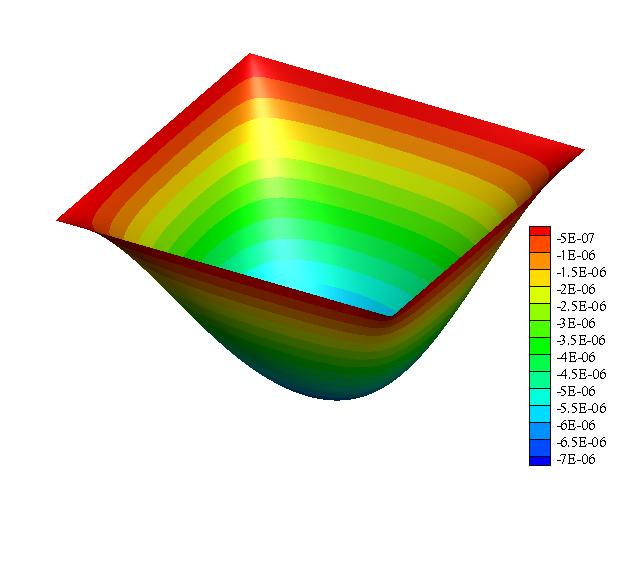}\\
  (a)
\end{minipage}
\begin{minipage}[c]{0.4\textwidth}
  \centering
  \includegraphics[width=50mm]{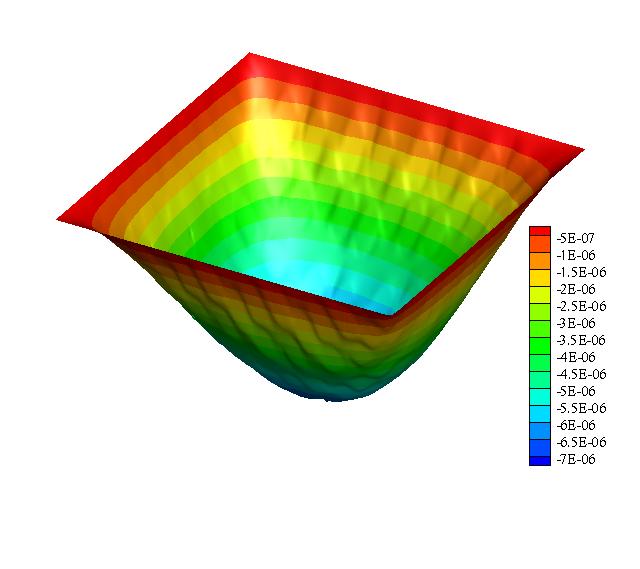}\\
  (b)
\end{minipage}
\begin{minipage}[c]{0.4\textwidth}
  \centering
  \includegraphics[width=50mm]{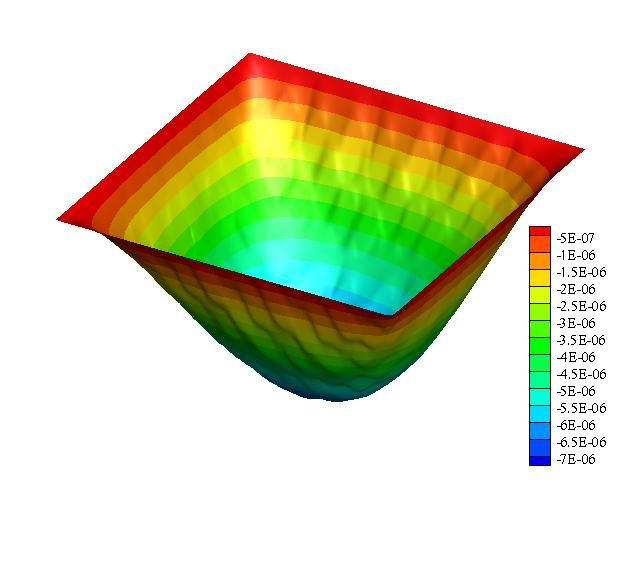}\\
  (c)
\end{minipage}
\begin{minipage}[c]{0.4\textwidth}
  \centering
  \includegraphics[width=50mm]{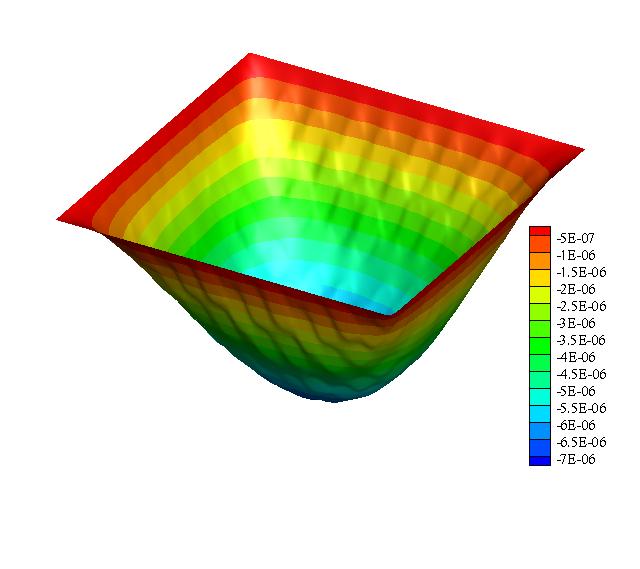}\\
  (d)
\end{minipage}
\caption{The numerical results of temperature field at $t=0.2s$: (a) $u_{10}$; (b) $u_1^{(1\varepsilon)}$; (c) $u_1^{(2\varepsilon)}$; (d) $u_{1\text{DNS}}$.}
\end{figure}

\begin{figure}[!htb]
\centering
\begin{minipage}[c]{0.4\textwidth}
  \centering
  \includegraphics[width=50mm]{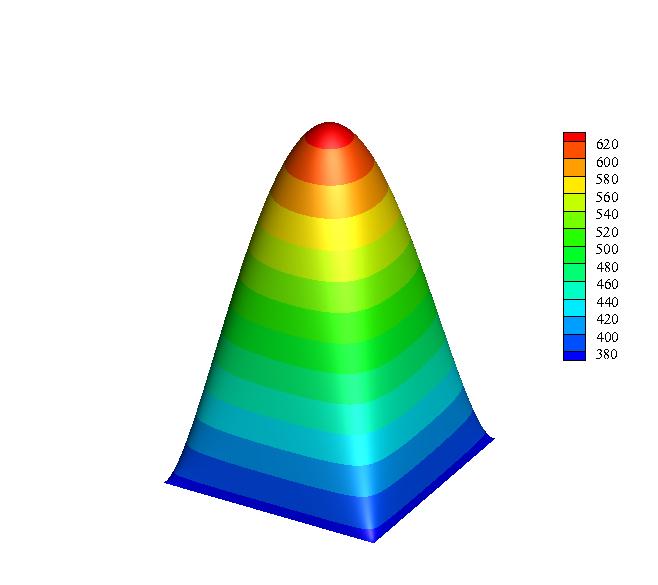}\\
  (a)
\end{minipage}
\begin{minipage}[c]{0.4\textwidth}
  \centering
  \includegraphics[width=50mm]{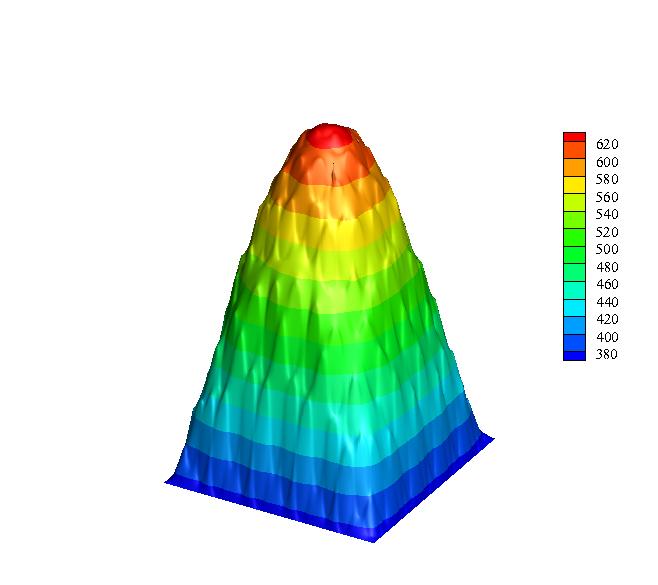}\\
  (b)
\end{minipage}
\begin{minipage}[c]{0.4\textwidth}
  \centering
  \includegraphics[width=50mm]{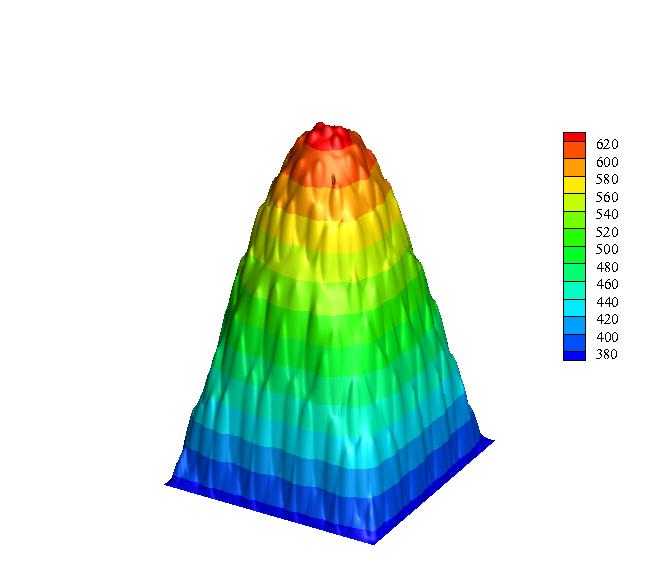}\\
  (c)
\end{minipage}
\begin{minipage}[c]{0.4\textwidth}
  \centering
  \includegraphics[width=50mm]{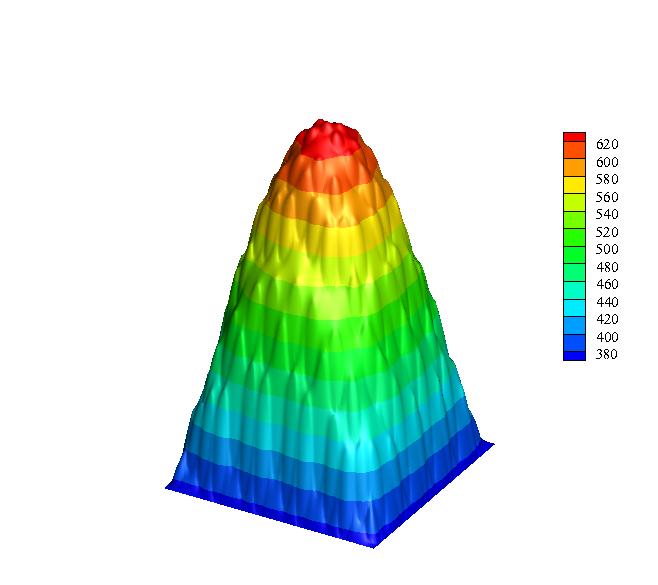}\\
  (d)
\end{minipage}
\caption{The numerical results of temperature field at $t=1.0s$: (a) $T_0$; (b) $T^{(1\varepsilon)}$; (c) $T^{(2\varepsilon)}$; (d) $T_{\text{DNS}}$.}
\end{figure}

\begin{figure}[!htb]
\centering
\begin{minipage}[c]{0.4\textwidth}
  \centering
  \includegraphics[width=50mm]{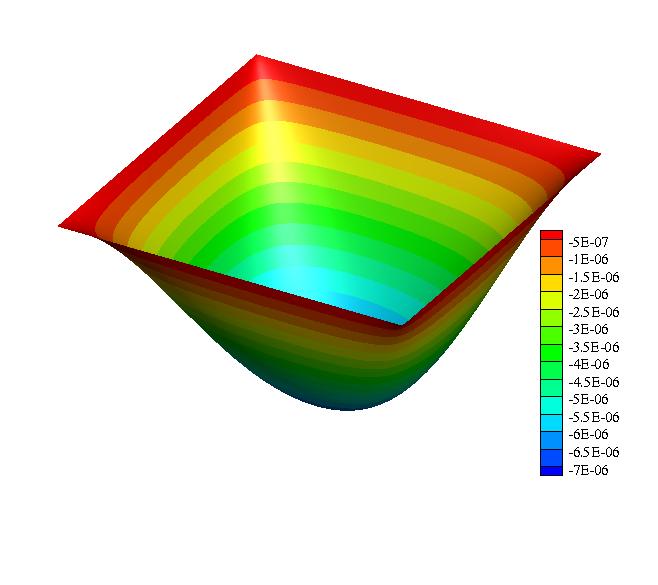}\\
  (a)
\end{minipage}
\begin{minipage}[c]{0.4\textwidth}
  \centering
  \includegraphics[width=50mm]{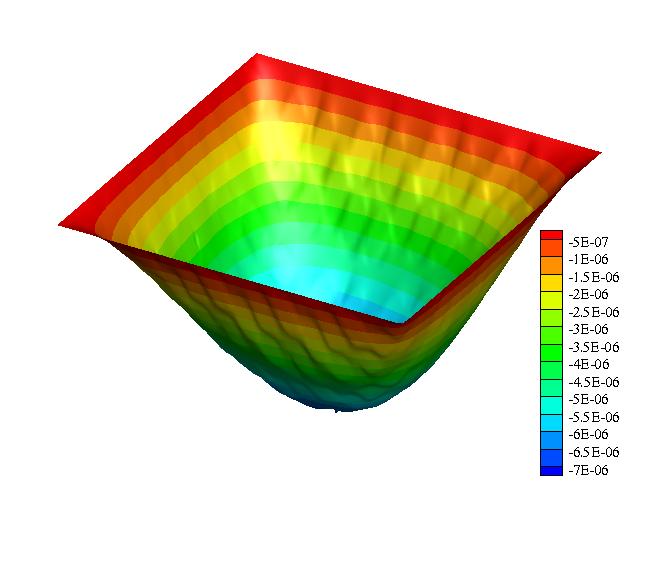}\\
  (b)
\end{minipage}
\begin{minipage}[c]{0.4\textwidth}
  \centering
  \includegraphics[width=50mm]{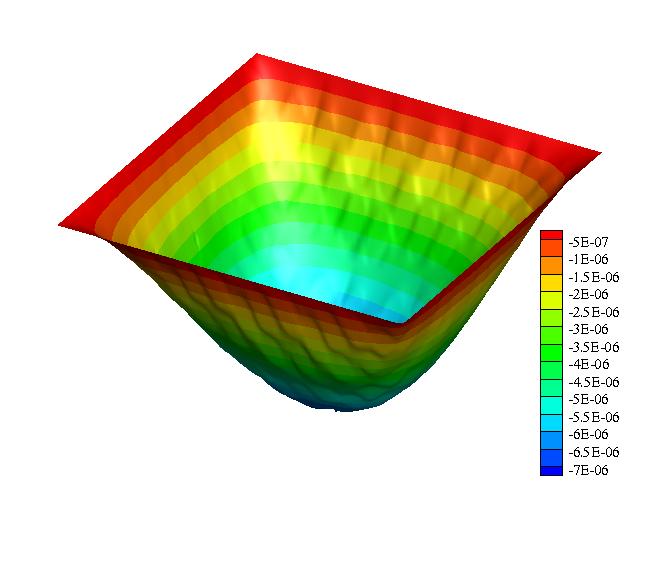}\\
  (c)
\end{minipage}
\begin{minipage}[c]{0.4\textwidth}
  \centering
  \includegraphics[width=50mm]{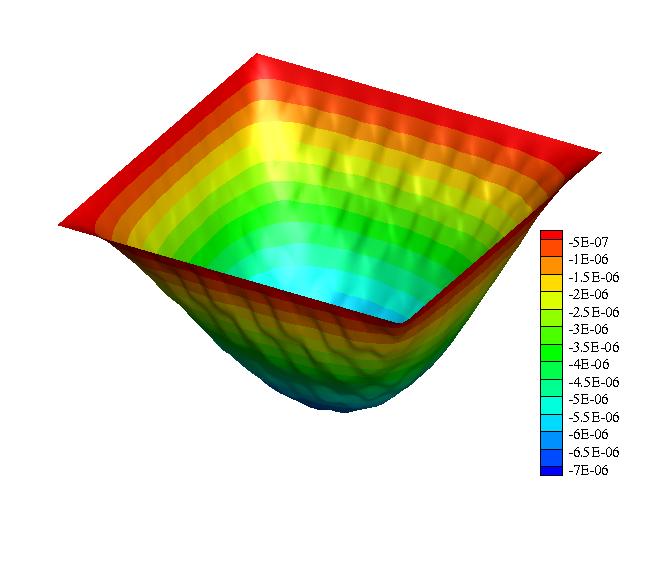}\\
  (d)
\end{minipage}
\caption{The numerical results of temperature field at $t=1.0s$: (a) $u_{10}$; (b) $u_1^{(1\varepsilon)}$; (c) $u_1^{(2\varepsilon)}$; (d) $u_{1\text{DNS}}$.}
\end{figure}

\begin{figure}[!htb]
\centering
\begin{minipage}[c]{0.4\textwidth}
  \centering
  \includegraphics[width=50mm]{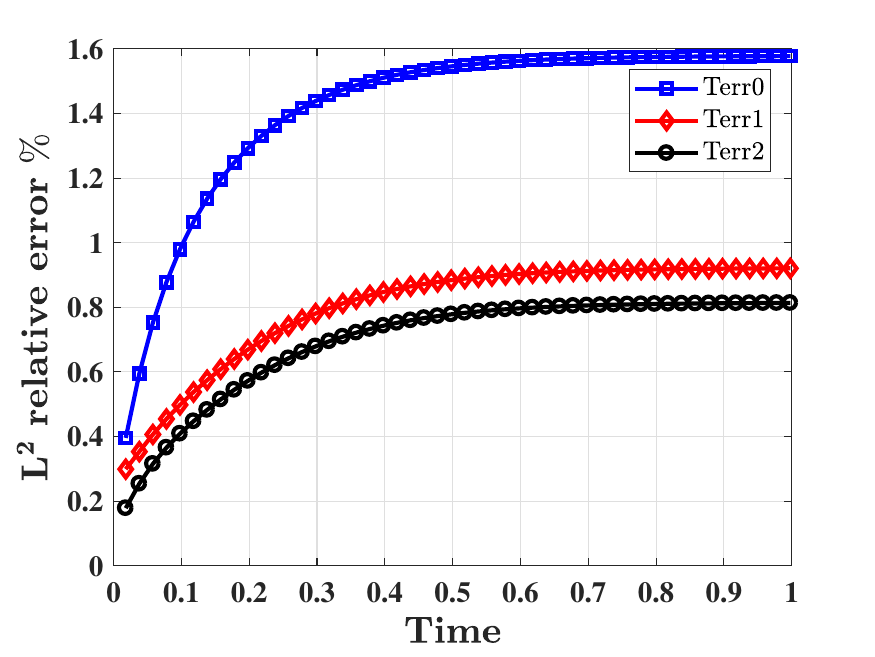}\\
  (a)
\end{minipage}
\begin{minipage}[c]{0.4\textwidth}
  \centering
  \includegraphics[width=50mm]{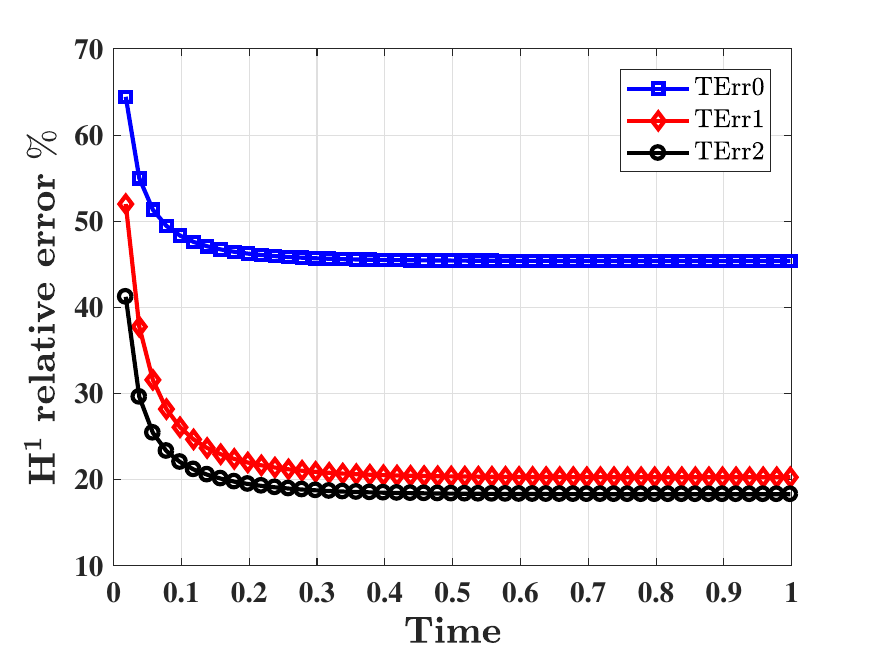}\\
  (b)
\end{minipage}
\begin{minipage}[c]{0.4\textwidth}
  \centering
  \includegraphics[width=50mm]{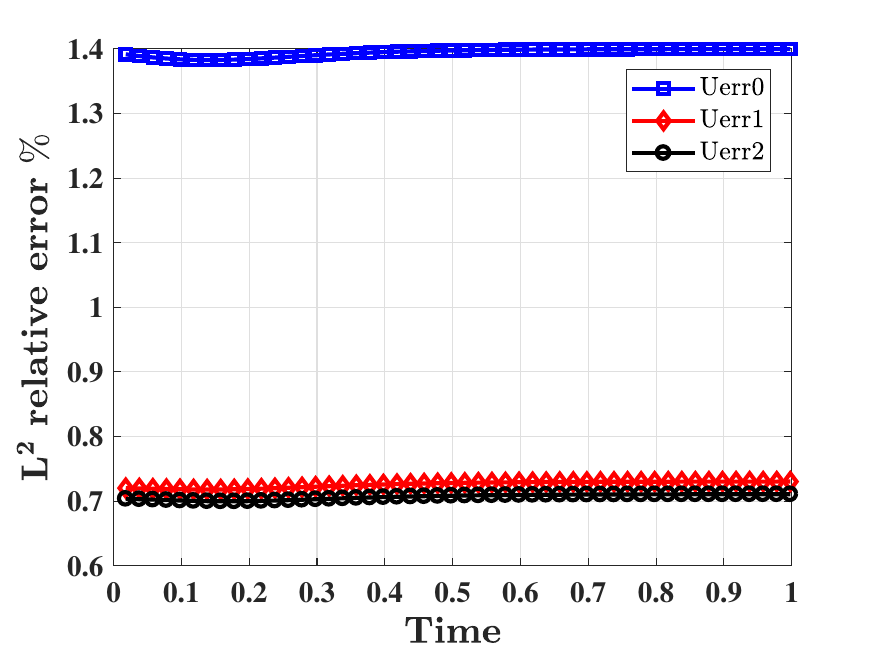}\\
  (c)
\end{minipage}
\begin{minipage}[c]{0.4\textwidth}
  \centering
  \includegraphics[width=50mm]{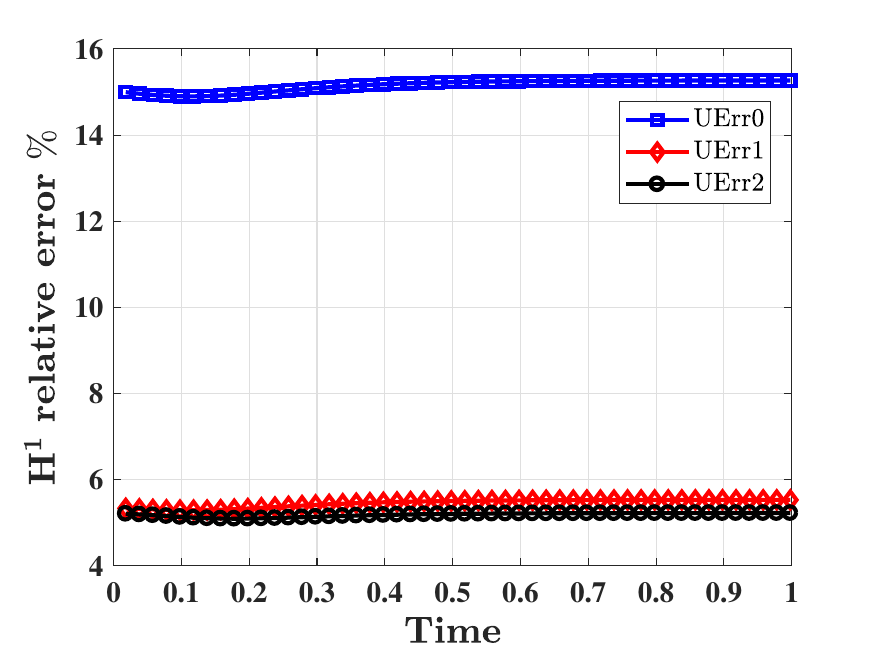}\\
  (d)
\end{minipage}
\caption{The evolutive relative errors of temperature and displacement fields: (a) \rm{Terr}; (b) \rm{TErr}; (c) \rm{Uerr}; (d) \rm{UErr}.}
\end{figure}

As indicated Table 5, the presented SHOMS method can significant reduce computer memory without losing precision. Moreover, the numerical results in Figs.\hspace{1mm}12-15 reveal that the higher-order multi-scale solutions can accurately capture the microscopic oscillatory behaviors and provide preferable approximations to the exact solutions of the investigated 2D composite structure compared with macroscopic homogenized solutions and lower-order multi-scale solutions, especially for temperature field. The evolutive relative errors shown in Fig.\hspace{1mm}16 clearly demonstrate the accuracy and stability of the two-stages space-time multi-scale numerical algorithm in the long-time numerical simulation. Furthermore, it is important to highlight that the presented SHOMS approach remains effective even for a relatively small parameter $\varepsilon$, which corresponds to a large number of microscopic unit cells in inhomogeneous structures. In contrast, the high-resolution DNS simulation fails to converge for the investigated large-scale problems. This prominent computational advantage of the proposed SHOMS framework is of significant practical value of the SHOMS approach in engineering computations.
\section{Conclusions and outlook}
In the present work, a novel statistical higher-order multi-scale method is developed with the aim of effectively simulating nonlinear thermo-mechanical simulation of random composite materials with temperature-dependent properties, which served
under extreme heat environment. The main contributions of this work are threefold: First, the statistical multi-scale formulations with the higher-order correction terms are established for random composites under statistical periodic configurations. Second, the local error estimations for the statistical multi-scale solutions of nonlinear thermo-mechanical systems are derived in detail. Third, a space-time numerical algorithm with off-line and on-line stages is designed to overcome the prohibitive computation of direct numerical simulation. Furthermore, numerical results demonstrate that the presented SHOMS approach can effectively simulate nonlinear thermo-mechanical coupling behaviors with less computational cost and accurately capture the microscopic oscillatory information caused by randomly heterogeneous configurations. Besides, the proposed SHOMS approach can accurately predict the equivalent material parameters of random composites compared with the predictive results of some theoretical models, which illustrate that high temperature field has a remarkable effect on macroscopic thermo-mechanical properties.

In the future, the SHOMS method will be extended to more complex nonlinear problems including thermal convection and radiation effects under extreme thermal environment. Additionally, machine learning approaches and parallel algorithm will be introduced in the off-line stage of SHOMS framework, in order to avoid repetitive statistical computation and improve computational efficiency.
\section*{Acknowledgments}
The authors gratefully acknowledge the support of the the National Natural Science Foundation of China (Nos.\hspace{1mm}12001414 and 51739007), Young Talent Fund of Association for Science and Technology in Shaanxi, China (No.\hspace{1mm}20220506), Young Talent Fund of Association for Science and Technology in Xi'an, China (No.\hspace{1mm}095920221338), the Fundamental Research Funds for the Central Universities (No.\hspace{1mm}KYFZ23020).

\FloatBarrier

\section*{References}
\bibliography{reference_SHTSP}

\newpage
\appendix
\end{document}